\documentclass[reqno,12pt]{amsart} 
\usepackage[margin=2.85cm,footskip=1cm]{geometry}

\usepackage{amsmath} 
\usepackage{amssymb}      
\usepackage{amscd}      
\usepackage{amsthm}      
\usepackage{epsfig}      
\usepackage{amstext}      
\usepackage[all]{xy} 
\usepackage{isolatin1}
\usepackage{multirow}

 \newcommand{\dist}{{\operatorname{dist}}}
   \newcommand{\Hom}{\operatorname{Hom}}

\newcommand{\Ad}{\operatorname{Ad}}

\newcommand{\id}{\operatorname{id}}

\newcommand{\inter}{\operatorname{int}}  
\newcommand{\Aff}{\operatorname{Aff}} 
\newcommand{\diam}{\operatorname{diam}}

\newcommand{\diag}{\operatorname{diag}}

\newcommand{\var}{\operatorname{var}}

\newcommand{\Per}{\operatorname{Per}}
\newcommand{\Gl}{\operatorname{Gl}}

   \theoremstyle{plain}%default
   \newtheorem{thm}{Theorem}[section]
   \newtheorem{prop}[thm]{Proposition}
   \newtheorem{lemma}[thm]{Lemma}  
   \newtheorem{cor}[thm]{Corollary}
   \theoremstyle{definition}
   
   \newtheorem{defn}[thm]{Definition}
   \newtheorem{example}[thm]{Example}
   \theoremstyle{remark}
   
   \newtheorem{remark}[thm]{Remark}

   \numberwithin{equation}{section}

\title[1-solenoids]{The homoclinic and heteroclinic $C^*$-algebras
  of a generalized one-dimensional solenoid}
  
\author{Klaus Thomsen}

        \date{\today}
\newcommand\xsquid[1]{\mathrel{\smash{\overset{#1}{\rightsquigarrow}}}}
\newcommand\soverline[1]{\smash{\overline{#1}}}

\DeclareMathOperator\coker{coker}

\date{}

\email{matkt@imf.au.dk}
\address{Institut for matematiske fag, Ny Munkegade, 8000 Aarhus C, Denmark}

\begin{document}

\maketitle

\section{Introduction}

The homoclinic and heteroclinic structure in dynamical systems was
first used to produce $C^*$-algebras in the way breaking work of Cuntz
and Krieger in \cite{Kr1}, \cite{Kr2} and \cite{CuK}. This work has
been generalized in many directions where the relation to dynamical
systems is either absent or appears very implicit, but Ian Putnam
described in \cite{Pu} a natural way to extend the constructions of
Cuntz and Krieger to higher dimensions such that the point of
departure is the heteroclinic structure in a Smale space, just as the
work of Cuntz and Krieger departed from the heteroclinic structure in
a shift of finite type, which is a zero-dimensional Smale space.
Putnam builds his approach on the work of D.~Ruelle, \cite{Ru1},
\cite{Ru2}, who introduced the notion of a Smale space in \cite{Ru1}
and constructed the so-called asymptotic algebra from the homoclinic
equivalence relation in \cite{Ru2}.

The work of Putnam and Ruelle was further generalized by the author in
\cite{Th1} where it was shown that Ruelle's approach can be adopted as
soon as there is enough expansiveness in the underlying dynamical
system; the local product structure in a Smale space is not
crucial for the construction. 
Furthermore, in \cite{Th1} the alternative approach was used to obtain
inductive limit decompositions for the algebras of Putnam arising from
particular classes of
Smale spaces, e.g. expansive group automorphisms and one-dimensional
generalized solenoids in the sense of R.F. Williams, \cite{Wi2}, and
I. Yi, \cite{Y1}. For expansive group automorphisms it was shown
that the $C^*$-algebras are all AT-algebras of real rank zero, and
hence are classified by their K-theory groups, thanks to the work of
G. Elliott, \cite{Ell1}. For one-dimensional generalized solenoids the exact nature
of the inductive limit decomposition was not determined and the
homoclinic algebra was not examined. In
particular, it was not decided if the $C^*$-algebras are classified
by K-theory. The main
purpose of the present paper is tie up this loose end by showing that
they are, although they turn out to be more general AH-algebras and exhibit
more complicated K-theory than the algebras arising from expansive
group automorphisms, at least in the sense that torsion
appears. Specifically, it is shown that the heteroclinic algebra of
both a
one-solenoid and its inverse, as well as the homoclinic algebra are
all AH-algebras of real rank zero with no dimension growth. They
are therefore classified by K-theory thanks to the work of Elliott and
Gong, \cite{EG}. This conclusion is obtained for the heteroclinic
algebra by combining a thorough study of the inductive limit
decomposition obtained in \cite{Th1} with results on the
classification of simple $C^*$-algebras, in particular results by
H. Lin on algebras of tracial rank zero, cf. e.g. \cite{Lin4}.

The study of the $K_1$-group
of the heteroclinic algebra reveals a fundamental dichotomy between
the orientable and the non-orientable case. For orientable
one-solenoids the $K_1$-group of the heteroclinic algebra is $\mathbb
Z$ while it collapses to $\mathbb Z_2$ when the one-solenoid is not
orientable. It follows from this that a one-solenoid is only conjugate
to an orientable one-solenoid when the graph map defining it is already
oriented.

The dichotomy between the oriented and non-oriented case
turns out to be
of pervasive importance for the structure of the three algebras
arising from the heteroclinic and homoclinic structure of a
one-solenoid. For example, one of these algebras is an AT-algebra if
and only if they all are, and this happens if and only if the
one-solenoid is orientable. In the oriented case the heteroclinic algebra of a one-solenoid is not
so different from the heteroclinic algebra of its inverse; they are
both AT-algebra and have the same $K_1$-group, namely $\mathbb
Z$. In fact the two algebras are often, but not always, isomorphic in
this case. In contrast the heteroclinic algebra of the inverse of a
one-solenoid differs substantially from the heteroclinic algebra of
the one-solenoid itself in the non-orientable case. It
turns out that the heteroclinic algebra of the inverse of a
one-solenoid is the stabilization of a crossed product $C^*$-algebra
coming from a free and minimal action of the infinite dihedral group
on the Cantor set, and the results on its structure are obtained
through a study of such crossed products. Again the work of Lin on
algebras with tracial rank zero plays a fundamental role, one of the
crucial steps being that we can use results of N.C. Phillips on finite
group actions with the tracial Rokhlin property, \cite{Ph4},\cite{Ph5}, to
conclude that the crossed product $C^*$-algebra
coming from a free and minimal action of the infinite dihedral group
on the Cantor set has tracial rank zero. Among the results we obtain
is that
the $K_1$-group vanishes for the heteroclinic algebra of the inverse of
a non-orientable one-solenoid, but nonetheless the algebra is not an
AF-algebra because the two-torsion now pops up in the $K_0$-group.

The results we obtain on the structure of the homoclinic algebra of a one-solenoid is
obtained from the results on the heteroclinic algebras by use of a
result of Ian Putnam from \cite{Pu} which says that the homoclinic
algebra is stably isomorphic to the crossed product of the two
heteroclinic algebras.

\smallskip

\emph{Acknowledgement:}  I want to thank Marcy Barge and Ian Putnam for their
          help with this paper. In particular, I thank Marcy for
          convincing me that an unorientable one-solenoid is a Smale space,
          and Ian for many discussions and for pointing out a mistake in the first version
          which concealed the role of the infinite dihedral group.

\section{The dynamical setup} 

\subsection{The homoclinic and heteroclinic algebras of an expansive
  homeomorphism}

Let $(X,d)$ be a compact metric space and $\psi : X \to X$ a
homeomorphism. Recall that $\psi$ is \emph{expansive} when there is a
$\delta > 0$ such that
\begin{equation*}
x \neq y \ \Rightarrow \ \sup_{k \in \mathbb Z}
d\left(\psi^k(x),\psi^k(y)\right) \geq \delta .
\end{equation*}  
We say that two points $x,y \in X$ are \emph{locally conjugate} when
there are open neighborhoods $U$ and $V$ of $x$ and $y$,
respectively, and a homeomorphism $\chi : U \to V$ such that
$\chi(x) = y$ and $\lim_{k \to \pm \infty} \sup_{z \in U}
d\left(\psi^k(z),\psi^k(\chi(z))\right) = 0$. This is an equivalence
relation on $X$ and its graph 
\begin{equation*}
R_{\psi}(X) = \left\{ (x,y) \in X^2 : \ x \ \text{is locally conjugate
    to} \ y \right\}
\end{equation*}
is a locally compact Hausdorff space in the topology for which every
local conjugacy $(U,V,\chi)$ defines an element of a canonical base:
\begin{equation*}
\left\{ \left(z,\chi(z)\right) : \ z \in U \right\} .
\end{equation*}
This topology is typically different from the topology which
$R_{\psi}(X)$ inherits from $X \times X$; it has more open sets. The
crucial fact is that $R_{\psi}(X)$ is what is nowadays known as an
\'etale equivalence relation, cf. \cite{Ph1}, \cite{GPS1}, \cite{Th1},
so that the reduced groupoid $C^*$-algebra
$C^*_r\left(R_{\psi}(X)\right)$ of Renault, \cite{Re}, can
be defined. We call this \emph{the homoclinic algebra} of $(X,\psi)$
and denote it by $A_{\psi}(X)$. 
When $(X,\psi)$ is a Smale space, as
defined by Ruelle in \cite{Ru1}, the homoclinic algebra
$A_{\psi}(X)$ is the asymptotic algebra of Ruelle and Putnam,
\cite{Ru2},\cite{Pu}. When $X$ is a sub-shift the homoclinic algebra is
the AF-algebra whose dimension group was defined by Krieger in
\cite{Kr1} and Section~2 of \cite{Kr2}.

Let $\Per {\psi}$ denote the set of $\psi$-periodic
points. We assume that $\Per {\psi} \neq \emptyset$.
For $p \in \Per
\psi$, set
\begin{equation*}
W^u(p) = \bigl\{ x \in X : \ \lim_{k \to - \infty}
  d\left(\psi^k(x), \psi^k(p)\right) \bigr\} = 0.
\end{equation*}  
Since $\psi$ is expansive each $W^u(p)$ is a locally compact
Hausdorff space in a topology with base 
\begin{equation*}
\left\{ y \in X : d\left( \psi^j (y),
    \psi^j (x)\right) < \epsilon , j \leq k\right\},
\end{equation*}
where $x \in W^u(p), \ k \in \mathbb Z$, $\epsilon \in
\left]0,\epsilon_p\right[$ are arbitrary,
and $\epsilon_p > 0$ only depends on $p$. See Lemma 4.6 of \cite{Th1}. The set of \emph{post-periodic
  points} of $\left(X,\psi \right)$ is 
\begin{equation*}
W_{X, \psi} \ = \bigcup_{p \in \Per
  \psi} W^u(p), 
\end{equation*}
and it is a locally compact Hausdorff space in the topology, which we
call \emph{the Wagoner topology}, defined such that each $W^u(p)$ is
open in $W_{X, \psi}$ and has the topology we
have just described above.

Define the equivalence relation $\sim$ on $W_{X, \psi}$ such
that $x \sim y$ if and only if there are open neighborhoods, $U$ of
$x$ and $V$ of $y$ in $W_{X, \psi}$, and a
homeomorphism $\gamma : U \to V$, again called \emph{a local conjugacy}, such
that $\gamma(x) =y$ and
\begin{equation*}
\lim_{k  \to \infty} \sup_{z \in U} d\left(\psi^k(z),
  \psi^k(\gamma(z))\right) = 0 .
\end{equation*} 
This is an \'etale equivalence relation
$R_{\psi}\left(X, W_{X,\psi}\right)$, cf. \cite{Th1}, and the corresponding (reduced) groupoid
$C^*$-algebra
$C^*_r\left(R_{\psi}\left(X, W_{X,\psi}\right)\right)$ is \emph{the heteroclinic algebra} of
$\left(X, \psi\right)$. As in \cite{Th1} we denote it by
$B_{\psi}\left(X\right)$. When $(X,\psi)$ is a mixing Smale space the
heteroclinic algebra is $*$-isomorphic to the stabilization of Putnams
'stable algebra', cf. \cite{Pu} and Theorem 4.17 of \cite{Th1}. In particular, it is
a higher dimensional analogue
of the 'AF-core' in the Cuntz-Krieger construction, cf. \cite{CuK}.

\subsection{Generalized  one-dimensional solenoids}\label{onedimsol} Let
$\Gamma$ be a finite (unoriented) graph with vertexes $\mathbb V$ and
edges $\mathbb E$. Consider a continuous map $h : \Gamma \to \Gamma$ such that the following conditions are satisfied for
some metric $d$ for the topology of~$\Gamma$: 
\begin{enumerate}
\item[$\alpha$)] (Expansion) There are constants $C > 0$ and $\lambda > 1$ such that
\begin{equation*}
d\left(h^n(x),h^n(y)\right) \geq C\lambda^n d(x,y)
\end{equation*} 
for every $n \in \mathbb N$ when $x,y \in e \in \mathbb E$ and there
is an edge $e' \in \mathbb E$ with $h^n\left(\left[x,y\right]\right)
\subseteq\nobreak  e'$. ($\left[x,y\right]$ is the interval in $e$
between $x$ and $y$.)
\item[$\beta$)] (Non folding) $h^n$ is locally injective on $e$ for each $e \in \mathbb E$ and each $n \in \mathbb N$.
\item[$\gamma$)] (Markov) $h\left(\mathbb V \right) \subseteq \mathbb
  V$.
\item[$\delta$)] (Mixing) For every edge $e \in \mathbb E$ there is an
  $m \in \mathbb N$ such that $\Gamma \subseteq h^m(e)$.
\item[$\epsilon$)] (Flattening)  There is a $d \in \mathbb N$ such that for all $x \in \Gamma$ there is a neighborhood $U_x$ of $x$ with $h^d\left(U_x\right)$ homeomorphic to $]-1,1[$.
\end{enumerate}
When all conditions hold we are in a setting first introduced by
Williams in \cite{Wi2} and later studied by I.Yi in
\cite{Y1},\cite{Y2},\cite{Y3},\cite{Y4},\cite{Y5}.  
We say then that $(\Gamma, h)$ is
a \emph{pre-solenoid}. \footnote{Unlike what I
  believed when \cite{Th1} was written, the combined conditions
are strictly stronger than the conditions considered by Williams and Yi. If for
example $\Gamma$ consists of two disjoint circles and $h$ takes one
circle twice around the other in an appropriate expanding way (like $z
\mapsto z^2$ on the unit circle in the complex plane), then $\delta)$
fails, but all the Axioms 0-5 of \cite{Y1} are satisfied. The example shows that Lemma 2.14
of \cite{Y1} is wrong. The mistake occurs in the proof of 1.6 Lemma in \cite{Wi2}.}

Set
\begin{equation*}
\overline{\Gamma} = \left\{ \left(x_i\right)_{i=0}^{\infty} \in \Gamma^{\mathbb N} : \ h\left(x_{i+1}\right) = x_i, \ i = 0,1,2, \dots \right\} .
\end{equation*}
We consider $\overline{\Gamma}$ as a compact metric space with the metric
\begin{equation*}
D\bigl(\left(x_i\right)_{i=0}^{\infty} , 
\left(y_i\right)_{i=0}^{\infty}\bigr) = \sum_{i=0}^{\infty} 2^{-i} d\left(x_i,y_i\right) .
\end{equation*}
Define $\overline{h} : \overline{\Gamma} \to \overline{\Gamma}$ such
that $\overline{h}(x)_i = h\left(x_i\right)$ for all $i \in \mathbb
N$. $\overline{h}$ is a homeomorphism with inverse
\begin{equation*}
\overline{h}^{-1}\left(z_0,z_1,z_2, \dots \right) = \left(z_1,z_2,
  z_3, \dots \right).
\end{equation*}
Following Williams and Yi, \cite{Wi2}, \cite{Y1}, we call
$\left(\overline{\Gamma}, \overline{h}\right)$ a \emph{generalized
  one-dimensional solenoid} or just a
\emph{1-solenoid}.

It was shown by Williams in \cite{Wi1} that expanding attractors of
certain diffeomorphisms of compact manifolds are 1-solenoids via a
conjugacy which turns the restriction of the diffeomorphism into $\overline{h}$,
and that each 1-solenoid arises in this way from a diffeomorphism of
the 4-sphere.

As we shall see there is a dichotomy in the class of 1-solenoids which
is crucial for the structure of the $C^*$-algebras they give rise to and
which depends on whether or not the 1-solenoid is orientable. To
formalize this notion we define \emph{an orientation of $\Gamma$} to
be a collection of homeomorphisms $\psi_e : [0,1] \to e, \ e \in
\mathbb E$. We say that $h$ is \emph{positively oriented}
(resp. \emph{negatively oriented})  with respect
to the orientation $\psi_e, e \in \mathbb E$, when the function
\begin{equation}\label{function}
\psi_{e'}^{-1} \circ h \circ \psi_e : \ \psi_e^{-1}\left(e \cap
  h^{-1}\left(e'\right) \right) \to [0,1]
\end{equation}
is increasing (resp. decreasing) for every pair $e,e' \in \mathbb
E$. A pre-solenoid $(\Gamma,h)$ is \emph{positively oriented}
(resp. \emph{negatively oriented}) when there is an
orientation of the edges in $\Gamma$ such that $h$
is positively oriented (resp. negatively) oriented with respect to
that orientation. $\left(\Gamma, h\right)$ is \emph{oriented} when it
is either positively or negatively oriented.

The 1-solenoid
$\left(\overline{\Gamma},\overline{h}\right)$ is \emph{orientable}
where there is an oriented pre-solenoid $\left(\Gamma_1,h_1\right)$
such that $\left(\overline{\Gamma},\overline{h}\right)$ is conjugate
to $\left(\overline{\Gamma_1},\overline{h_1}\right)$. When
$\left(\Gamma_1,h_1\right)$ can be chosen to be positively
(resp. negatively) oriented we say that $\left(\overline{\Gamma},
  \overline{h}\right)$ is \emph{positively (resp. negatively) orientable}.

\begin{remark}\label{exampl} Our terminology concerning orientation of
  pre-solenoids is more elaborate than that used by Yi, cf. \cite{Y1}-\cite{Y5}. To illustrate the difference, let $\Gamma$ be the wedge of two
  circles $a$ and $b$, and let $h : \Gamma \to \Gamma$ be given by the
  wrapping rule 
$$
a \mapsto a^{-1}b^{-1}, \ b \mapsto a^{-1}b^{-1}
$$
in the sense of Definition 2.6 in \cite{Y1}. Then $(\Gamma, h)$ is
a pre-solenoid which is negatively oriented but not positively
oriented. In contrast, the wrapping rule 
$$
a \mapsto ba, \ b \mapsto ba
$$
is positively oriented and not negatively oriented. Despite the fact that the square of the
two maps are identical and positively oriented, the 1-solenoids they
define are not conjugate.
\end{remark}

Through the study of the heteroclinic $C^*$-algebra
$B_{\overline{h}}\left(\overline{\Gamma}\right)$ of a 1-solenoid
$\left(\overline{\Gamma},\overline{h}\right)$, in particular a careful
examination of its $K_1$-group, we shall obtain the following result.

\begin{thm}\label{orientthm} Let $(\Gamma, h)$ be a pre-solenoid. Then
  $\left(\overline{\Gamma}, \overline{h}\right)$ is positively
  orientable (resp. negatively orientable) if and
  only if $\left(\Gamma, h\right)$ is positively oriented
  (resp. negatively oriented).
\end{thm}

\section{The heteroclinic algebra of a 1-solenoid}\label{het}

In this section we study the structure of the heteroclinic algebra of a 1-solenoid. The point of departure is the inductive limit
decomposition of the heteroclinic algebra obtained in Theorem 5.17 of \cite{Th1}.

\subsection{The building blocks}\label{bblocks}

Let $\Gamma$ be a finite graph and $h : [-1,1] \to \Gamma$ a locally injective continuous map. We define an equivalence relation $\sim$ on $]-1,1[$ such that $t \sim s$ if and only if $h(t) = h(s)$ and there are open neighborhoods $U_s$ and $U_t$ of $s$ and $t$ in $]-1,1[$, respectively, such that $h\left(U_s\right) = h\left(U_t\right) \simeq ]-1,1[$. Set
\begin{equation*}
R_h = \left\{ (s,t) \in ]-1,1[^2 : \ s \sim t\right\} .
\end{equation*}
Give $R_h$ the topology inherited from $]-1,1[^2$.

\begin{lemma}\label{etaleR} $R_h$ is an \'etale equivalence relation.
\begin{proof} This is part of Lemma 5.13 in \cite{Th1}.  
\end{proof}
\end{lemma}

We are going to use the \'etale equivalence relations of Lemma
\ref{etaleR} in the special case where $h(-1), h(1) \in \mathbb V$ and
$h(]-1,1[) = \Gamma$. When this holds we say that $R_h$ is an \emph{open
  interval-graph relation}.

To give a manageable description of $C^*_r\left(R_h\right)$, let $m, n \in \mathbb N$. Consider some $a = \left(a(1),a(2), \dots,
  a(m)\right) \in \mathbb N^m, \ b = \left( b(1),b(2), \dots,
  b(n)\right) \in \mathbb N^n$ and two $n \times m$-matrices, $I,U$,
with $\{0,1\}$-entries. Assume that
\begin{equation}\label{standalg1}
\sum_{k =1}^m I_{ik}a(k) \leq b(i)
\end{equation}
and
\begin{equation}\label{standalg2}
\sum_{k =1}^m U_{ik}a(k) \leq b(i)
\end{equation}
for all $i = 1,2, \dots, n$. Set $F_{a} = M_{a(1)} \oplus M_{a(2)}
\oplus \dots \oplus M_{a(m)}$ and $F_{b} = M_{b(1)} \oplus M_{b(2)} \oplus
\dots \oplus M_{b(n)}$, where $M_k$ denotes the $C^*$-algebra of
complex $k\times k$ matrices. Let $\varphi^{I}, \varphi^{U} : F_{a} \to F_{b}$
be $*$-homomorphisms with partial maps $\varphi^I_{ij}, \varphi^U_{ij}
: M_{a(j)} \to M_{b(i)}$ of multiplicity $I_{ij}$ and $U_{ij}$,
respectively. Set
\begin{equation*}
A\left(a, b, I,U\right) = \left\{ \left(x,f\right) \in
  F_{a} \oplus C\left([0,1], F_{b}\right) : \ \varphi^I(x) =
  f(0), \ \varphi^{U}(x) = f(1) \right\} .
\end{equation*}
In the following we shall only be interested in algebras of this type when
\begin{equation}\label{equal2}
\sum_{i=1}^n I_{ik} + \sum_{i=1}^n U_{ik} = 2
\end{equation}
for all $k = 1,2, \dots, m$. We will refer to such an algebra as a
\emph{building block}. Note that $A\left(a, b, I,U\right)$ is unital if and only if
\begin{equation*}\label{unitalcond}
\sum_{k=1}^m U_{ik}a(k) = \sum_{k =1}^m I_{ik}a(k) = b(i) 
\end{equation*}
for all $i$.

For the explicit
identification of $C^*_r\left(R_h\right)$ with a building block we need
to decide for each edge which of its endpoints corresponds to $0 \in
[0,1]$ and which to $1 \in [0,1]$. For this purpose we give each edge
of $\Gamma$ an (arbitrary) orientation so that they may be considered as directed arrows instead of undirected
edges. The start vertex of the arrow $e$ is then denoted by $s(e)$ and the
terminal vertex of $e$ by $t(e)$.  For each edge $e$ of $\Gamma$ choose a homeomorphism
  $\psi_e : e \to [0,1]$ such that $\psi_e(s(e)) = 0, \ \psi_e(t(e)) =
  1$.

Let $-1=x_0 < x_1 < x_2 < \dots < x_N = 1$ be the elements of
$h^{-1}\left(\mathbb V\right)$. \emph{A passage in $\Gamma$} is a closed subset $J \subseteq \Gamma$
containing exactly one vertex $v$ such that there
is a homeomorphism $\varphi : J \to \left[-1,1\right]$ with
$\varphi^{-1}(0) = v$. We will identify two passages $J$ and $J'$ in
$\Gamma$ when $J\cap J'$ is also a passage in $\Gamma$. 
 A \emph{$h$-passage} is a passage in
$\Gamma$ which contains $h\left(\left[x_i-\epsilon,
    x_i+\epsilon\right]\right)$ for some $i \in \{1,2, \dots, N-1\}$
and all sufficiently small $\epsilon > 0$. An arrow $e$ is then an \emph{entry-arrow} in the $h$-passage
$p$ when $s(e) \in p$ and an \emph{exit-arrow} when $t(e) \in p$.

Note that for every element $x \in h^{-1}(\mathbb V)$ there is a
unique $h$-passage $p_x$ such that $p_x = h\left(\left[x-\epsilon, x+
    \epsilon\right]\right)$ for all small $\epsilon >0$.
Let $\mathbb A_h$
denote the finite-dimensional $C^*$-algebra generated by the
matrix-units $e_{x,y}, x,y \in h^{-1}(\mathbb V)$, such that $p_x =
p_y$. Similarly, let
$\mathcal I_h$ denote the set of connected components of $]0,1[ \backslash
h^{-1}\left(\mathbb V\right)$ and let $\mathbb B_h$ denote the
finite-dimensional $C^*$-algebra generated by the matrix units
$e_{I,J}$, where $I,J \in \mathcal I_h$ and
$h\left(\overline{I}\right) = h\left( \overline{J}\right)$. Define
$\pi^I : \mathbb A_h \to \mathbb B_h$ such that
$\pi^I\left(e_{x,y}\right) = \sum_{J,K} e_{J,K}$ where we sum over
the pairs $J,K \in \mathcal I_h$ with the property that $x \in
\overline{J}, y \in \overline{K}, h\left(\overline{J}\right) = h
\left(\overline{K}\right)$ and $ h\left(\overline{J}\right)$ is an 
entry-arrow in both $p_x$ and $p_y$. There are at most two such pairs
- when there is none we set $\pi^I\left(e_{x,y}\right) = 0$. Define $\pi^U : \mathbb A_h \to \mathbb B_h$ such that
$\pi^U\left(e_{x,y}\right) = \sum_{J,K} e_{J,K}$ where we sum over
the pairs $J,K \in \mathcal I_h$ with the property that $x \in
\soverline{J}, y \in \soverline{K}, h\left(\soverline{J}\right) =
h\left(\soverline{K}\right)$ and $h\left(\soverline{J}\right)$ is an 
exit-arrow in both $p_x$ and $p_y$. There are at most two such pairs
- when there is none we set $\pi^U\left(e_{x,y}\right) = 0$. Set
\begin{equation*}
\mathbb D_h = \left\{ (x,f) \in \mathbb A_h \oplus C\left([0,1],
    \mathbb B_h\right) : \ f(0) = \pi^I(x), \ f(1) =
  \pi^U(x)\right\} .
\end{equation*}
We can then define a $*$-isomorphism $\Phi_h : C^*_r\left(R_h\right) \to
\mathbb D_h$ such that 
\begin{equation*}
\Phi_h (f)  = \biggl( \sum_{x,y \in h^{-1}(\mathbb V)} f(x,y)e_{x,y} ,
  \sum_{I,J \in \mathcal I_h} f\Bigl( \left(h|_{\overline{I}}\right)^{-1} \circ
    \psi_{h\left(\overline{I}\right)}^{-1} ( \cdot), \left(h|_{\overline{J}}\right)^{-1} \circ
    \psi_{h\left(\overline{J}\right)}^{-1} ( \cdot)\Bigr)
  e_{I,J}\biggr) 
\end{equation*}
when $f \in C_c\left(R_h\right)$. This is essentially the
$*$-isomorphism from Lemma 5.16 of \cite{Th1}, revised to avoid the
assumption about the absence of loops in $\Gamma$ which was imposed there.

Let $p_i, i = 1,2, \dots, m$, be a numbering of the different $h$-passages and $e_i, i = 1,2, \dots, n$, a numbering of the edges in
$\Gamma$. Set
\begin{equation*} a(i) = \ \# \left\{ j \in \{1,2 \dots, N-1\} : \
  h\left(\left[x_{j} -\epsilon,x_{j}+\epsilon\right]\right) =  p_i  \
  \text{for all small} \ \epsilon \right\},
\end{equation*}
$ i = 1,2, \dots, m$,
and
\begin{equation*}
b(i) = \  \# \left\{ j \in \{0,1,\dots, N-1\} : \
  h\left(\left[x_j,x_{j+1}\right]\right) =  e_i\right\},
\end{equation*}
$i = 1,2, \dots, n$. Set $a  = \left(a(1),a(2),
  \dots, a(m)\right) \in \mathbb N^m, \ 
b = \left(b(1),b(2), \dots, b(n)\right) \in \mathbb N^n$. Then there are obvious $*$-isomorphisms
$\kappa : \mathbb A_h \to F_a$ and $\kappa' : \mathbb B_h \to F_b$.   
Let
\begin{equation*}
I_{ik} = \begin{cases} 1, \ & \ \text{when \ $e_i$ \ is an
    entry-edge in} \ p_k \\ 0, & \ \text{otherwise}, \end{cases}
\end{equation*}
and
\begin{equation*}
U_{ik} = \begin{cases} 1, \ & \ \text{when \ $e_i$ \ is an
    exit-edge in} \ p_k \\ 0, & \ \text{otherwise}. \end{cases}
\end{equation*}
Note
that (\ref{equal2}) holds, i.e. $A\left(a, b, I,U\right)$ is
a building block. Furthermore, $\varphi^I = \kappa' \circ \pi^I \circ
\kappa^{-1}$ and $\varphi^U = \kappa' \circ \pi^U \circ
\kappa^{-1}$, at least up to unitary equivalences which we can safely
ignore. We get therefore a $*$-isomorphism $\Psi_h : \mathbb D_h \to
A(a,b,I,U)$ defined such that
\begin{equation*}
\Psi_h (x,f) = \left( \kappa(x), \left(\id_{C[0,1]} \otimes
    \kappa'\right)(f)\right) .
\end{equation*}
In this way we get
\begin{lemma}\label{isolemma} (cf. Lemma 5.16 of \cite{Th1})  $\Psi_h
  \circ \Phi_h : C^*_r(R_h) \ \to \ A\left(a, b, I,U\right)$ is a $*$-isomorphism.

\end{lemma}

Note that the building block $A(a,b,I,U)$ of Lemma \ref{isolemma}
satisfies condition (\ref{nybeting}) of the following lemma by construction.

\begin{lemma}\label{projex} Let $A(a,b,I,U)$ be a building block
  corresponding to the vectors $a = \left(a(1),a(2), \dots, a(m)\right)
  \in \mathbb N^m$
  and $b = \left(b(1),b(2), \dots, b(n)\right) \in \mathbb N^n$. Assume
  that  
\begin{equation}\label{nybeting}
\sum_{j=1}^m U_{kj} \geq 1,  \ \ \sum_{j=1}^m
  I_{kj} \geq 1
\end{equation}
for all $k = 1,2, \dots, n$, and $a(i) \geq 2n+1$ for all $i = 1,2, \dots,
m $. Then $A(a,b,I,U)$ contains a non-zero projection.

\begin{proof} We define a labeled graph $\mathcal G$ as follows: The vertexes of
$\mathcal G$ consist of the tuples $(i,\to), i =1,2, \dots, n$, and
the tuples $(i, \leftarrow), i = 1,2, \dots, n$. The arrows
$\rightsquigarrow$ in $\mathcal G$ are labeled by the elements of
$\{1,2, \dots, m\}$ and there is a labeled arrow $(i,\to)
%\overset{j}\rightsquigarrow 
\xsquid{j}
(i',\to)$ when $\varphi^U_{ij} \neq 0$
and $\varphi^I_{i'j} \neq 0$, there is a labeled arrow $(i,\to)
%\overset{j}\rightsquigarrow 
\xsquid{j}
(i',\leftarrow)$ when $i \neq i'$, $\varphi^U_{ij}
\neq 0$ and $\varphi^U_{i'j} \neq 0$, there is a labeled arrow
$(i,\leftarrow ) 
%\overset{j}\rightsquigarrow 
\xsquid{j}
(i',\leftarrow)$ when
$\varphi^I_{ij} \neq 0$ and $\varphi^U_{i'j} \neq 0$ and there is a
labeled arrow $(i,\leftarrow ) 
%\overset{j}\rightsquigarrow 
\xsquid{j}
(i',\to)$
when $i \neq i'$, $\varphi^I_{ij} \neq 0$ and $\varphi^I_{i'j} \neq
0$. 

%If there is an arrow in $\mathcal G$ of the form 
%\begin{equation}\label{form0}
%(i,\to)
%%\overset{j}\rightsquigarrow 
%\xsquid{j}
%(i,\to)
%\end{equation}
%we choose rank 1 projections $p
%\in C\left( [0,1], M_{b(i)}\right)$ and $q \in M_{a(j)}$ such that $p(1)
%= \varphi^U_{ij}(q)$ and $\varphi^I_{ij}(q) = p(0)$. Then $(q,p)$ is
%a non-zero projection in $A(a,b,I,U)$. A similar construction gives a
%non-zero projection when there is an arrow of the form
%\begin{equation}\label{form00}
%(i,\leftarrow)
%\overset{j}\rightsquigarrow 
%\xsquid{j}
%(i,\leftarrow) .
%\end{equation}
%We will therefore assume in the following that there are no arrows in
%$\mathcal G$ of the form (\ref{form0}) or (\ref{form00}).

Given an
arrow $(i,\to) 
%\overset{j}\rightsquigarrow 
\xsquid{j}
(i',\to)$ in $\mathcal G$ we write $p
%\overset{q}\rightsquigarrow 
\xsquid{q}
p'$ when $p \in
C\left([0,1],M_{b(i)}\right)$, $p' \in
C\left([0,1],M_{b(i')}\right)$ and $ q \in M_{a(j)}$ are rank 1 projections
such that $p(1) = \varphi^U_{ij}(q)$ and $\varphi^I_{i'j}(q) =
p'(0)$. We say that $p
%\overset{q}\rightsquigarrow 
\xsquid{q}
p'$ \emph{realizes} $(i,\to )
%\overset{j}\rightsquigarrow 
\xsquid{j}
(i',\to)$. In the following all
projections, with a single obvious exception, will be rank 1
projections. We say that $p
%\overset{q}\rightsquigarrow 
\xsquid{q}
p'$ realizes $(i,\to )
%\overset{j}\rightsquigarrow 
\xsquid{j}
(i',\leftarrow)$ when  $p(1) = \varphi^U_{ij}(q)$ and $\varphi^U_{i'j}(q) =
p'(1)$, that $p
%\overset{q}\rightsquigarrow 
\xsquid{q}
p'$ realizes $(i,\leftarrow )
%\overset{j}\rightsquigarrow 
\xsquid{j}
(i',\leftarrow)$ when  $p(0) = \varphi^I_{ij}(q)$ and $\varphi^U_{i'j}(q) =
p'(1)$ and finally that $p
%\overset{q}\rightsquigarrow 
\xsquid{q}
p'$ realizes $(i,\leftarrow )
%\overset{j}\rightsquigarrow 
\xsquid{j}
(i',\to)$ when  $p(0) = \varphi^I_{ij}(q)$ and $\varphi^I_{ij}(q) =
p'(0)$.
Given a finite path
\begin{equation*}
c = \left(i_1,*_1\right)  
%\overset{j_1}\rightsquigarrow
\xsquid{j_1}
\left(i_2,*_2\right) 
%\overset{j_2}\rightsquigarrow 
\xsquid{j_2}
\cdots
%\overset{j_{d-1}}\rightsquigarrow
\xsquid{j_{d-1}}
\left(i_d,*_d\right)   
\end{equation*}
in $\mathcal G$, where $*_k \in \left\{ \leftarrow, \to\right\}$, we
say that $c$ is realized by $p_1 
%\overset{q_1}\rightsquigarrow 
\xsquid{q_1}
p_2
%\overset{q_2}\rightsquigarrow 
\xsquid{q_2}
p_3 
%\overset{q_3}\rightsquigarrow 
\xsquid{q_3}
\cdots
%\overset{q_{d-1}}\rightsquigarrow 
\xsquid{q_{d-1}}
p_d$ when $p_k
%\overset{q_k}\rightsquigarrow 
\xsquid{q_k}
p_{k+1}$ realizes $\left(i_k,*_k \right)
%\overset{j_k}\rightsquigarrow 
\xsquid{j_k}
\left(i_{k+1},*_{k+1}\right)$ for all $k
= 1,2, \dots, d-1$, and, in addition, the projections $p_1,p_2, \dots,
p_d$ are mutually orthogonal in $C\left([0,1], F_b\right)$, and
$q_1,q_2, \dots, q_{d-1}$ mutually orthogonal in $F_a$.

We claim that every path in $\mathcal G$ of length $\leq k-1$ can be
realized when $ \min_i a(i) \geq k$. To prove this by induction in
$k$, observe that a path of length $1$ has the form $(i_1,*_1)
%\overset{j}\rightsquigarrow 
\xsquid{j}
(i_2,*_2)$ which is trivial to realize when $a(j) \geq \min_i a(i)
\geq 2$. Assume then that the assertion is true for some $k\geq 2$, that $\min_i a(i) \geq k+1$ and consider a path
\begin{equation*}
c = \left(i_1,*_1\right)  
%\overset{j_1}\rightsquigarrow
\xsquid{j_1}
\left(i_2,*_2\right) 
%\overset{j_2}\rightsquigarrow 
\xsquid{j_2}
\cdots
%\overset{j_{k}}\rightsquigarrow
\xsquid{j_k}
\left(i_{k+1},*_{k+1}\right)   
\end{equation*}
of length $k$. By induction hypothesis there are projections such
that $p_1
%\overset{q_1}\rightsquigarrow 
\xsquid{q_1}
p_2 
%\overset{q_2}\rightsquigarrow 
\xsquid{q_2}
p_3
%\overset{q_3}\rightsquigarrow 
\xsquid{q_3}
\cdots 
%\overset{q_{k-1}}\rightsquigarrow
\xsquid{q_{k-1}}
p_k$ realizes $ \left(i_1,*_1\right)  
%\overset{j_1}\rightsquigarrow
\xsquid{j_1}
\left(i_2,*_2\right) 
%\overset{j_2}\rightsquigarrow 
\xsquid{j_2}
\cdots
%\overset{j_{k-1}}\rightsquigarrow
\xsquid{j_{k-1}}
\left(i_{k},*_{k}\right)$. Consider
first the case where $*_k =  \ \to$ and $*_{k+1} = \ \to$. Since
$a(j_k) \geq k+1$ there is a projection $q_{k} \in M_{a(j_k)}$
which is orthogonal in $F_a$ to each $q_{i}, i \leq k-1$. It follows that $\varphi^U_{i_k j_k}\left(q_{k}\right)$ is orthogonal in
$F_b$ to $p_i(1)$ for all $i \leq k-1$, with the possible exception of
$p_1(1)$ when $*_1 = \ \leftarrow$ and $i_1 = i_k$. When this happens we
change $p_1(t)$ for $t$ in a neighborhood of $1$ to arrange that
$\varphi^U_{i_k j_k}\left(q_{k}\right) \perp p_1(1)$, while keeping
$p_1,p_2, \dots, p_{k}$ orthogonal. This is possible because
$b\left(i_k\right) \geq \min_j a(j) \geq k+1$. So we can assume, in
all cases, that $\varphi^U_{i_k j_k}\left(q_{k}\right)$ is orthogonal
to $p_i(1)$ for all $i \leq k-1$. There is then
a rank 1 projection $p'_k \in C\left( [0,1], M_{b(i_k)}\right)$ such that
$p'_k(t) = p_k(t), t \leq \frac{1}{2}$, $p'_k \perp p_i$ in
$C\left([0,1], F_b\right)$ for all $i \leq k-1$, and $p'_k(1) =
\varphi^U_{i_k j_k}\left(q_{k}\right)$. Exchanging $p'_k$ for $p_k$
we may thus assume that $p_k(1) =
\varphi^U_{i_k j_k}\left(q_{k}\right)$.

Note that
$\varphi^I_{i_{k+1} j_k}\left(q_{k}\right)$ is orthogonal to
$p_i(0)$ in $F_b$ for $i \leq k$, with the possible exception of
$p_1(0)$ when $i_1 = i_{k+1}$ and $*_1 = \ \to$. When this happens we
change $p_1(t)$ for $t$ in a neighborhood of $0$ to arrange that
$\varphi^I_{i_{k+1} j_k}\left(q_{k}\right) \perp p_1(0)$, while keeping
$p_1,p_2, \dots, p_{k}$ orthogonal. This is possible because
$b\left(i_{k+1}\right) \geq \min_j a(j) \geq k+1$. So we can assume, in
all cases, that $\varphi^I_{i_{k+1} j_k}\left(q_{k}\right)$ is orthogonal
to $p_i(0)$ for all $i \leq k$. Using once more that $b(i_{k+1}) \geq k+1$ we can find a projection $p_{k+1}$ in
$C\left([0,1], M_{b(i_{k+1})}\right)$ which is orthogonal in
$C\left([0,1], F_b\right)$ to $p_i, i \leq k$, and
satisfies that $p_{k+1}(0) =\varphi^I_{i_{k+1}
  j_k}\left(q_{k}\right)$. Then $p_1
%\overset{q_1}\rightsquigarrow 
\xsquid{q_1}
p_2 
%\overset{q_2}\rightsquigarrow 
\xsquid{q_2}
p_3
%\overset{q_3}\rightsquigarrow 
\xsquid{q_3}
\cdots 
%\overset{q_{k-1}}\rightsquigarrow
\xsquid{q_{k-1}}
p_k 
%\overset{q_k}\rightsquigarrow 
\xsquid{q_k}
p_{k+1}$ realizes $c$. Since
exactly the same proof works in the other three case, $*_k = \ \to$ and
$*_{k+1} = \ \leftarrow$, $*_k =  \ \leftarrow$ and $*_{k+1} = \ \to$, and
finally $*_k = \ \leftarrow$ and $*_{k+1} = \ \leftarrow$, we have
completed the induction step. %The case $k =1$ follows from the
%observation that any arrow $(i,*)
%%\overset{j}\rightsquigarrow 
%\xsquid{j}
%(i',*)$ in $\mathcal G$ can be realized
%because we have abandoned the cases where $i = i'$.  

Note that it follows from conditions (\ref{equal2}) and
(\ref{nybeting}) that $\mathcal G$ has no sinks or sources. There is
therefore a loop
\begin{equation*}
c = \left(i_1,*_1\right)  
%\overset{j_1}\rightsquigarrow
\xsquid{j_1}
\left(i_2,*_2\right) 
%\overset{j_2}\rightsquigarrow 
\xsquid{j_2}
\cdots
%\overset{j_{d-1}}\rightsquigarrow
\xsquid{j_{d-1}}
\left(i_d,*_d\right)
%\overset{j_d}\rightsquigarrow
\xsquid{j_d}
\left(i_1,*_1\right) ,
\end{equation*}
in $\mathcal G$ whose length ($=d$) does not exceed $2n$. Set
$\left(i_{d+1},*_{d+1}\right) =
\left(i_1,*_1\right)$. Since we assume that
$a(i) \geq 2n+1$ for all $i$, it follows from the preceding that there
are projections such that $p_1
%\overset{q_1}\rightsquigarrow 
\xsquid{q_1}
p_2 
%\overset{q_2}\rightsquigarrow 
\xsquid{q_2}
p_3
%\overset{q_3}\rightsquigarrow 
\xsquid{q_3}
\cdots 
%\overset{q_{d}}\rightsquigarrow
\xsquid{q_d}
p_{d+1}$ realizes $c$. Assume $*_1 = \ \leftarrow$. Change $p_1$ to
$p'_1$ such that $p'_1(t) =p_1(t), t \leq \frac{1}{2}$, and $p'_1(1) =
\varphi^U_{i_1j_d}\left(q_d\right)$. Since
$\varphi^U_{i_1j_d}\left(q_d\right)$ is orthogonal to $p_j(1),
2 \leq j \leq d$, in $F_b$ and $p_1$ to $p_j, 2 \leq j \leq d$, in $C\left([0,1],F_b\right)$, we can do this such that the
elements of $\left\{p'_1,p_2, \dots, p_{d}\right\}$ are mutually
orthogonal in $C\left([0,1], F_b\right)$. Then $\bigl( \sum_{j=1}^d
  q_j, p'_1 + \sum_{j=2}^d p_j\bigr)$ is a non-zero projection in
$A(a,b,I,U)$. This completes the proof because the case $*_1 = \ \to$ is
completely analogous.
\end{proof}
\end{lemma}

Note that there is an extension
\begin{equation*}
\begin{xymatrix}{
0 \ar[r] &  C_0(0,1) \otimes F_b \ar[r] & A(a,b,I,U) \ar[r] & F_a \ar[r]  & 0 
}
\end{xymatrix}
\end{equation*}
which makes it easy to
calculate the $K$-theory groups of a building block $A\left(a, b,
  I,U\right)$:
\begin{lemma}\label{Ktheory1}
$K_0\left( A\left(a, b, I,U\right)\right)$ is isomorphic, as
a partially ordered group, to
\begin{equation*}
\left\{ z =\left(z_1,z_2, \dots, z_m\right) \in \mathbb Z^m : \
  \left(I -U\right)z = 0 \right\} ,
\end{equation*}
and 
\begin{equation*}
K_1\left(  A\left(a, b, I,U\right)\right) \ \simeq \
\text{coker} \ (I - U) \ = \ \mathbb
Z^n/\left(I - U\right)\left(\mathbb Z^m\right) .
\end{equation*}
\end{lemma}

\subsection{The inductive limit decomposition}

In the notation of Lemma 5.14 of \cite{Th1}, there is
an open interval $J \subseteq \Gamma$, contained in some edge of
$\Gamma$, with endpoints in
$h^{-1}(\mathbb V)$, such that $h^d(J) = \Gamma$ for some $d \in \mathbb N$, and such that
$B_{\overline{h}}\left(\overline{\Gamma}\right)$ is $*$-isomorphic to
$A \otimes \mathbb K$, where $A$ is the inductive limit of the
sequence
\begin{equation*}\label{limit0}
C^*_r\left(R^d\right) \ \subseteq \ C^*_r\left(R^{d+1}\right) \
\subseteq \ C^*_r\left(R^{d+2}\right) \ \subseteq \ \cdots 
\end{equation*}
Here $R^i = R_{h^i}$ is the open interval-graph relation defined by $h^i : J
\to \Gamma$ and the inclusions $C^*_r\left(R^i\right) \, \subseteq \,
C^*_r\left(R^{i+1}\right)$ are induced by considering $R^i$ as an open
sub-relation of $R^{i+1}$. To give a more detailed description
of the sequence (\ref{limit0}) we need some terminology which we now
introduce. A $*$-homomorphism $\varphi : C\left([0,1], M_{n}\right) \to C\left([0,1], M_m\right)$ is \emph{regular} when there are
continuous functions $g_1,g_2, \dots, g_k : [0,1]\to [0,1]$ such that
\begin{equation*}
\varphi\left(f\right) =  \diag \left( f \circ g_1, f\circ g_2,
  \dots, f \circ g_k, 0, \dots, 0\right) .
\end{equation*}
Of course, this requires that $kn \leq m$. $\varphi$ is unital when $kn =
m$. The functions $g_1, \dots ,g_k$ will be called \emph{the
  characteristic functions} of $\varphi$. A $*$-homomorphism $\chi :
\oplus_{i=1}^N C\left([0,1], M_{n(i)}\right)  \to \oplus_{i=1}^M C\left([0,1], M_{m(i)}\right)$ is \emph{regular} when the
$*$-homomorphisms $\chi_{ij} : C\left([0,1], M_{n(j)}\right) \to C\left([0,1],
  M_{m(i)}\right), i = 1,2, \dots, M, j = 1,2, \dots, N$, which it
defines are all regular. The union of the characteristic functions of
the $\chi_{ij}$'s will be called \emph{the characteristic functions}
of $\chi$. A $*$-homomorphism $\chi :  \oplus_{j=1}^N M_{n(j)}  \to\nobreak
\oplus_{j=1}^M M_{m(j)}$ will be said to be \emph{non-increasing in
  rank} when its multiplicity matrix $A^{\chi} =
\left(A^{\chi}_{ij}\right)$ is given by a map $\iota: \left\{1,2, \dots,
  N\right\} \to \left\{1,2, \dots, M\right\}$ such that
$$
A^{\chi}_{ij} = \begin{cases} 1 & \ \text{when $\iota(j) = i$} \\ 0 &
  \ \text{otherwise.} \end{cases}
$$

Let $A\left(a,  b, I,U\right)$  and $A(a',b',I',U')$ be building blocks. A $*$-homomorphism 
\begin{equation}\label{reghm}
\psi : A\left(  a,   b, I,U\right)  \to A\left( a',  
  b', I',U'\right)
\end{equation}
is \emph{regular} when there are $*$-homomorphisms $\mu :
C\left([0,1], F_{  b}\right) \to F_{b'}$, $\chi : F_{ 
  a} \to F_{ a'}$ and a regular $*$-homomorphism $\varphi :C\left([0,1], F_{  b}\right)
\to C\left([0,1], F_{b' }\right)$ such that 
\begin{enumerate}
\item[1)] $\mu(f)\chi(x) = \chi(x)\mu(f) = 0, \
f \in C\left([0,1], F_{ b}\right), x \in F_{  a}$,
\item[2)] $\psi(x,f) = \left( \chi(x) + \mu(f), \varphi(f)\right),
  \ (x,f) \in  A\left(  a,   b, I,U\right)$, and
\item[3)] $\chi$ is non-increasing in rank.
\end{enumerate}

We will refer to the $*$-homomorphism $\mu$ as \emph{the skew map}
of $\psi$. Note that the composition of two regular $*$-homomorphisms
is again a regular $*$-homomorphism.

Let $F$ be a finite-dimensional $C^*$-algebra. We define the \emph{spectral
  variation} $\var f$ of an element $f \in C\left([0,1],F\right)$ to
be the number
\begin{equation*}
\var f = \inf_{U} \sup_{s,t \in [0,1]} \left\|U(s)f(s)U(s)^* -
  f(t)\right\|,
\end{equation*}
where we take the infimum over all unitaries $U$ in
$C\left([0,1],F\right)$. The \emph{spectral variation} $\var z$ of an element $z = (x,f) \in
 A\left(  a,   b, I,U\right)$ of a building block is then defined to
 be the number $\var z = \var f$.

\begin{lemma}\label{limit} The heteroclinic
  algebra $B_{\overline{h}} \left(\overline{\Gamma}\right)$ is
  $*$-isomorphic to the inductive limit of a sequence of
  building blocks
\begin{equation}\label{sequence}
\begin{xymatrix}{
A'_1 \ar[r]^-{\pi_1} & A'_2 \ar[r]^-{\pi_2} & A'_3 \ar[r]^-{\pi_3} &
\dots }
\end{xymatrix} 
\end{equation}
 and injective regular $*$-homomorphisms with the property
that 
\begin{equation}\label{varlim}
\lim_{k \to \infty} \var \left[\pi_k \circ \pi_{k-1} \circ
\pi_{k-2} \circ \dots \circ \pi_i(a)\right] = 0
\end{equation}
for all $a \in A'_i$ and
all $i \in \mathbb N$.
\begin{proof} For each $i\geq d$ we let $A'_i$ be the building block
  obtained from $C^*_r\left(R_{h^i}\right)$ by use of Lemma
  \ref{isolemma}. Define $\chi : \mathbb A_{h^i} \to \mathbb A_{h^{i+1}}$ such that
$\chi\left(e_{x,y}\right) = e_{x,y}$ when $x,y \in h^{-i}(\mathbb V)$ and $p_x = p_y$. This is well-defined because $h$
maps $h^i$-passages to $h^{i+1}$-passages by the Markov condition
$\gamma)$ and the non folding condition $\beta)$. Define
$\mu : C\left([0,1], \mathbb B_{h^i}\right) \to \mathbb
A_{h^{i+1}}$ such that
\begin{equation*}
\mu\left( f \otimes e_{I,J}\right) = \sum_{(x,y) 
\in N_{I,J}} \ f\bigl(
  \psi_{h^i\left(\overline{I}\right)} \circ h^{i}(x)\bigr)
e_{x,y} 
\end{equation*}
when $f \in C[0,1]$ and $I,J \in \mathcal I_{h^i}$ are such that
$h^{i}(I) = h^i(J)$, and we sum over the set
\begin{equation*}
N_{I,J} = \left\{ (x,y) \in h^{-i-1}(\mathbb V)^2: \ x
  \in I, y \in J, \ h^i(x) = h^i(y) \right\} .
\end{equation*}
Define a $*$-homomorphism $\varphi : C\left([0,1], \mathbb B_{h^i}\right) \to
C\left([0,1], \mathbb B_{h^{i+1}}\right)$ such that
\begin{equation*}
\varphi\left( f \otimes e_{I,J}\right) = \sum_{\left(I_1,J_1\right)
  \in M_{I,J}} f \circ
  \psi_{h^i\left(\overline{I}\right)} \circ
  \bigl(h|_{h^i\left(\overline{I}\right) \cap h^{-1}\left(h^{i+1} \left(\overline{I_1}\right)\right)}\bigr)^{-1} \circ
  \psi_{h^{i+1}\left(\overline{I_1}\right)}^{-1} \otimes e_{I_1,J_1},\end{equation*} 
for $f \in C[0,1]$ and $I,J \in \mathcal I_{h^i}$ with
$h^i\left(\overline{I}\right) = h^i\left(\overline{J}\right)$, where
we sum over the set
\begin{equation*}M_{I,J} = \left\{\left(I_1,J_1\right) \in \mathcal
  I_{h^{i+1}}^2 : \ I_1 \subseteq I, \ J_1 \subseteq J, \
  h^{i+1}\left(\overline{I_1}\right) = h^{i+1}\left(\overline{J_1}\right)\right\}.
 \end{equation*}
Then $\chi,\mu, \varphi$ define a $*$-homomorphism $\psi : \mathbb D_{h^i}
\to \mathbb D_{h^{i+1}}$ such that $\psi(x,f) = \linebreak\left(\chi(x) +
  \mu(f),\varphi(f)\right)$. Note that the diagram
 \begin{equation*}\label{commdiag2}
\begin{xymatrix}{
\mathbb D_{h^i} \ar[r]^-{\psi}                                    &
\mathbb D_{h^{i+1}} \\
C^*_r\left(R_{h^{i}}\right) \ar[u]^{\Phi_{h^i}}
\ar@{}[r]|-{\displaystyle\subseteq}
&  C^*_r\left(R_{ h^{i+1}}\right) \ar[u]_{\Phi_{h^{i+1}}}  }
\end{xymatrix}
\end{equation*}
commutes. Set $\pi_i = \Psi_{h^{i+1}} \circ \psi \circ \Psi_{h^i}^{-1}$. Then
$\pi_i$ is a regular $*$-homomorphism and the diagram
\begin{equation*}\label{commdiag}
\begin{xymatrix}{
A'_i\ar[r]^-{\pi_i}                                    & A'_{i+1} \\
C^*_r\left(R_{h^{i}}\right) \ar[u]^{\Psi_{h^i} \circ \Phi_{h^i}}
\ar@{}[r]|-{\displaystyle\subseteq}
&  C^*_r\left(R_{h^{i+1}}\right) \ar[u]_{\Psi_{h^{i+1}} \circ \Phi_{h^{i+1}}}  }
\end{xymatrix}
\end{equation*}
commutes. The characteristic functions of $\pi_i$ are of the form
\begin{equation}\label{charfunct}
\psi_e \circ \bigl(h|_{\overline{I_j}}\bigr)^{-1}  \circ
\psi_{e'}^{-1} 
\end{equation}
where $e,e'$ are edges of $\Gamma$
and $I_1,I_2, \dots, I_k$ are the intervals in $e \cap
h^{-1}\left(e'\right)$ such that $h: \overline{I_j} \to e'$ are
homeomorphisms. It is then clear that the expansion condition $\alpha)$
on $h$ ensures that $\lim_{k \to \infty} \var \left[\pi_k \circ \pi_{k-1} \circ
\pi_{k-2} \circ \dots \circ \pi_i(a)\right] = 0$ for all $a \in A'_i$
and all~$i$.
 \end{proof}
\end{lemma}

 Let $H_{ij}$ be the number of times the edge
$e_i$ is covered by the edge $e_j$ under $h$, i.e. for arbitrary $x
\in \inter e_i$,
\begin{equation}\label{H}
H_{ij} = \# h^{-1}(x) \cap e_j .
\end{equation} 
It follows from the mixing condition $\epsilon)$ that $H$ is mixing,
i.e. for large enough $k$ the matrix $H^k$ has no zero entries. As seen in the proof of Lemma
  \ref{limit} the sequence (\ref{sequence}) is stationary in the sense
  that the characteristic functions of the connecting
  $*$-homomorphisms, the $\pi_i$'s, are the same. They are given by
  continuous functions $\chi^k_{ij} : [0,1] \to  [0,1], i,j = 1,2,
  \dots, n$, $k = 1,2,\dots, H_{ij}$, where $n$ is the number of edges
  in $\Gamma$, and $H_{ij}$ is the number from (\ref{H}). The
  characteristic functions have the following properties:
\begin{enumerate}
\item[n1)] The characteristic functions are injective.
\item[n2)] $ \bigcup_{j,k} \chi_{ij}^k\left([0,1]\right) = [0,1]$ for 
  all $i = 1,2, \dots, n$.
\item[n3)] When $(j,k) \neq (j',k')$ the intersection
$\chi_{ij}^k\left([0,1]\right) \cap \chi_{ij'}^{k'}\left([0,1]\right)$
consists of at most one point.
\item[n4)] For each $i,j$ the points $\chi^k_{ij}(0), k = 1,2, \dots,
  H_{ij}$, are different. 
\end{enumerate}
If there is a pair $i,j$ such that there are both increasing and
decreasing functions among $\chi_{ij}^k, k = 1,2, \dots, H_{ij}$, we
can use the fact that $H$ is mixing and telescope the sequence
(\ref{sequence}) to arrange that this is the case for all $i,j$. If
instead it holds that for all $i,j$ the functions $\chi_{ij}^k, k =
1,2, \dots, H_{ij}$, are either all increasing or all decreasing we
telescope first (\ref{sequence}) by removing every second level to
make sure that, in addition, the functions
$\chi_{11}^k, k = 1,2, \dots, H_{11}$, all are increasing. If this is
not possible we are back in the first case so assume that
$\chi_{11}^k, k = 1,2, \dots, H_{11}$, are all increasing. When $j$ is
such that the functions $\chi_{1j}^k, k = 1,2, \dots, H_{1j}$, are all
decreasing we can arrange by telescoping that
$\left\{\chi^k_{ij}\right\}$ contains both increasing and decreasing
functions \emph{unless} $\chi^k_{j1}, k = 1,2, \dots, H_{j1}$, are also
  all decreasing. In this case we change the orientation of $e_j$. After this is done for
each such $j$ we telescope the sequence again by removing every second
level. There is then, for each $i,j$, at least one $k$ such that
$\chi^k_{ij}$ is increasing. So unless we are back in the first case,
all the characteristic functions must be increasing. In short, we can arrange by telescoping that 
\begin{enumerate}
\item[n5)] the
  characteristic functions are either
  all increasing or $\left\{ \chi^k_{ij} : \ k =
    1,2, \dots, H_{ij}\right\}$ contains both
  increasing and decreasing functions for all $i,j$.
\end{enumerate}

For convenience we shall substitute $h$ with a power $h^m$ of $h$ in
order to ensure further convenient properties of the connecting maps
of (\ref{sequence}). Such a substitution is justified by the following
lemma.

\begin{lemma}\label{invlemma1} For every $m \in \mathbb N$ the map
  $h^m : \Gamma \to \Gamma$ is a pre-solenoid,
  $B_{\overline{h^m}}\left(\overline{\Gamma}\right) \simeq
  B_{\overline{h}}\left(\overline{\Gamma}\right)$ and $B_{\overline{h^m}^{-1}}\left(\overline{\Gamma}\right) \simeq
  B_{\overline{h}^{-1}}\left(\overline{\Gamma}\right)$.
\begin{proof} The first statement is straightforward to check, and the
  others follow from the observation that $\left(\overline{\Gamma},
    \overline{h^m}\right)$ is conjugate to $\bigl(\overline{\Gamma},
    \overline{h}^m\bigr)$, combined with the observation that the
  heteroclinic algebra of a homeomorphism is $*$-isomorphic to the
  heteroclinic algebra of any of its positive powers.
\end{proof}
\end{lemma}

It follows from the mixing condition $\delta)$ that there is a an $m
\in \mathbb N$ such that $h^m(e) = \Gamma$ for every edge $e$ in
$\Gamma$.  
Hence, for a study of the heteroclinic algebra, we may
assume, by Lemma \ref{invlemma1}, that
\begin{equation}\label{mixing2}
h(e) = \Gamma, \ e \in \mathbb E.
\end{equation}
In fact we are going to move to higher powers of $h$ in order to
obtain other properties. Note that this will not violate
(\ref{mixing2}). Let $\mathcal P_i$ denote the $h^i$-passages. It
follows from (\ref{mixing2}) that $\mathcal P_i \subseteq \mathcal
P_{i+1}$. Since there are only finitely many passages in $\Gamma$ it
follows that there is an $m$ so big that $\mathcal P_{i+m} = \mathcal
P_m$ for all $i \geq 0$. A similar argument shows that if $m$ is big
enough the set of passages contained in $h^i\left(\inter e\right)$ is the
same for each $i \geq m$ and each $e \in \mathbb E$. Then
(\ref{mixing2}) shows that this 'stabilized' set of passages is
independent of~$e$. Hence by substituting $h$ with $h^m$ for some
sufficiently large $m$ we can arrange that there is a set $Q$ of
passages in $\Gamma$ such that
\begin{equation}\label{stab1}
\text{$Q$ is the set of $h^i$-passages for all $i \in \mathbb N$,}
\end{equation}
and
\begin{equation}\label{stab2}
\text{$Q$ is the set of passages contained in $h\left(\inter e\right)$
  for all $e \in \mathbb E$.}
\end{equation}  
By using the properties (\ref{mixing2}), (\ref{stab1}) and
(\ref{stab2}) in the proof of Lemma \ref{limit} we can regularize the
resulting sequence of $C^*$-algebras further. To describe the property
of the connecting maps which we obtain from (\ref{stab2}), consider a
regular $*$-homomorphism as in (\ref{reghm}) with skew map
$\mu$. We say that $\psi$ is \emph{full} when the composition
\begin{equation*}
C_0(0,1) \otimes M_{b(i)}  \hookrightarrow  C\left([0,1],F_b\right) \overset{\mu}\to
F_{a'} \to M_{a'(j)}
\end{equation*}
is non-zero for all $i,j$. By using the properties (\ref{mixing2}),
(\ref{stab1}) and (\ref{stab2}) in the proof of Lemma \ref{limit}, we obtain the following

\begin{lemma}\label{stabprop} There are natural numbers $n,m \in
  \mathbb N$ and $n \times m$ $\{0,1\}$-matrices $I,U$ and a set
  $\chi^k_{ij} : [0,1] \to [0,1]$, $k = 1,2, \dots,H_{ij}, i,j = 1,2,
  \dots, n$, of continuous functions such that $H_{ij} \geq 1$ for all
  $i,j$, and sequences
  $\{a_i\} \subseteq \mathbb N^m$, $\{b_i\} \subseteq \mathbb N^n$, such
  that $B_{\overline{h}}\left(\overline{\Gamma}\right)$ is the
  inductive limit of a sequence
\begin{equation}\label{sequence2}
\begin{xymatrix}{
A\left(a_1,b_1,I,U\right)  \ar[r]^-{\pi_1} & A\left(a_2,b_2,I,U\right) \ar[r]^-{\pi_2} & A\left(a_3,b_3,I,U\right) \ar[r]^-{\pi_3} &
\dots }
\end{xymatrix} 
\end{equation}
such that \eqref{varlim} holds and each $\pi_i$ is a full regular $*$-homomorphism with
$\left\{\chi_{ij}^k\right\}$ as characteristic functions.
\end{lemma}

\begin{lemma}\label{projlemma} Let $\varphi : A \to B$ be a regular
  $*$-homomorphism between building blocks $A$ and $B$. Let $e \in A$
  be a projection. There are unital building blocks, $A'$ and $B'$, a
  unital regular $*$-homomorphism $\varphi' : A' \to B'$ and
  $*$-isomorphisms $eAe \to A'$ and $\varphi(e)B\varphi(e) \to B'$
  such that
\begin{equation}
\begin{xymatrix}{
e A e \ar[d]_{\varphi} \ar[r] & A' \ar[d]^-{\varphi'} \\
\varphi(e) B \varphi(e) \ar[r] & B' }
\end{xymatrix}
\end{equation}
commutes.
\begin{proof} Left to the reader.
\end{proof}
\end{lemma}

In the following we let $\mathbb K$ denote the $C^*$-algebra of
compact operators on an infinite-dimensional separable Hilbert space.

\begin{prop}\label{limit3} There are natural numbers $n,m \in
  \mathbb N$ and $n \times m$ $\{0,1\}$-matrices $I,U$ and a set
  $\chi^k_{ij} : [0,1] \to [0,1]$, $k = 1,2, \dots,H_{ij}, i,j = 1,2,
  \dots, n$, of continuous functions such that $H_{ij} \geq 1$ for all
  $i,j$, and sequences
  $\{a_i\} \subseteq \mathbb N^m$, $\{b_i\} \subseteq \mathbb N^n$, such
  that $B_{\overline{h}}\left(\overline{\Gamma}\right)$ is
  $*$-isomorphic to $A \otimes \mathbb K$ where $A$ is the
  inductive limit of a sequence
\begin{equation}\label{sequence3}
\begin{xymatrix}{
A\left(a_1,b_1,I,U\right)  \ar[r]^-{\pi_1} & A\left(a_2,b_2,I,U\right) \ar[r]^-{\pi_2} & A\left(a_3,b_3,I,U\right) \ar[r]^-{\pi_3} &
\dots }
\end{xymatrix} 
\end{equation}
of unital building blocks such that
(\ref{varlim}) holds and each $\pi_i$ is a unital full regular $*$-homomorphism with
$\bigl\{\chi_{ij}^k\bigr\}$ as the characteristic functions.
\begin{proof} Let $n,m, U,I,H$ and $\left\{\chi^k_{ij}\right\}$ be as
  in Lemma \ref{stabprop}. The mixing condition combined with the
  fullness of the connecting maps in (\ref{sequence2}) guarantees that
\begin{equation*}  
\lim_{i \to \infty} \min_k a_i(k) = \infty
\end{equation*} 
and it follows therefore from Lemma \ref{projex} that the building
blocks in (\ref{sequence2}) contain projections, at least from a
certain stage. The fullness of the connecting maps ensure that these
projections become full projections in the sense of \cite{Br}. $B_{\overline{h}}
\left(\soverline{\Gamma}\right)$ is simple and stable by Lemma 5.10
of \cite{Th1} so it follows from \cite{Br} that $B_{\overline{h}}
\left(\soverline{\Gamma}\right) \simeq pB_{\overline{h}}
\left(\soverline{\Gamma}\right)p \otimes \mathbb K$ for any non-zero
projection $p \in B_{\overline{h}} \left(\soverline{\Gamma}\right)$.
As we have just argued we can assume that $p$ is a full projection in
the first building block occurring in (\ref{sequence2}). The proof is
then completed by use of Lemma \ref{projlemma}.
\end{proof}
\end{prop}

\subsection{Real rank zero and the consequences}

\begin{lemma}\label{uniquetrace} The $C^*$-algebra $A$ of Proposition \ref{limit3} 
has a unique trace state.
\begin{proof} Consider the sequence (\ref{sequence3}), and set
$A_i = A\left(a_i,b_i,
  I,U\right)$. There is a $*$-homomorphism $\iota_i : A_i \to
C\left([0,1], F_{b_i}\right) =B_i$ such that $\iota_i(x,f) = f$. Since the
connecting maps are regular we get an
infinite commuting diagram
\begin{equation}\label{sequence44}
\begin{xymatrix}{
A_1 \ar[r]^-{\pi_1} \ar[d]^-{\iota_1} & A_2 \ar[r]^-{\pi_2}
\ar[d]^-{\iota_2} & A_3 \ar[r]^-{\pi_3} \ar[d]^-{\iota_3} & \hdots \\
B_1 \ar[r]^-{\varphi_1} & B_2 \ar[r]^-{\varphi_2} &  B_3
\ar[r]^-{\varphi_3} & \hdots }
\end{xymatrix}
\end{equation}
where the $\varphi_i$'s are regular $*$-homomorphisms between
interval algebras. It follows then from condition (\ref{varlim}) that the inductive limit
$C^*$-algebra $B = \smash{\varinjlim} \left( B_i, \varphi_i\right)$ is
AF. Furthermore, the sequence $K_0\left(B_i\right)
%\overset{{\varphi_i}_*} \to 
\smash{\xrightarrow{{\varphi_i}_*}}
K_0\left(B_{i+1}\right)$ is stationary:
$K_0\left(B_i\right) \simeq \mathbb Z^n$ where $n$ is the number of
edges in $\Gamma$ and ${\varphi_i}_*$ is given by the matrix $H$, cf.
(\ref{H}). Since $H$ is mixing it follows that
$B$ is a unital simple AF-algebra with a unique
trace state, cf. Theorem~6.1 of \cite{Ef}. 

To complete the proof we need some notation. When $D$ is a unital $C^*$-algebra
  we denote by $T(D)$ its tracial state space and by $\Aff T(D)$ the
  order unit space of real-valued affine continuous functions on $T(D)$. When $\nu : D \to C$ is
  a $*$-homomorphism between unital $C^*$-algebras we let
  $\widehat{\nu} :  D_{sa} \to \Aff T(C)$ be the bounded linear map
  from the self adjoint part of $D$ to $\Aff T(C)$ given by
\begin{equation*}
\widehat{\nu}(d)(\omega) = \omega(\nu(d)), \ \omega \in T(C).
\end{equation*}
Note that
\begin{equation}\label{traceeq1}
\bigl\|\widehat{\nu' \circ \nu}\bigr\| \leq
\bigl\|\widehat{\nu'}\bigr\|\left\|\widehat{\nu}\right\|
% \left\|\widehat{\nu' \circ \nu}\right\| \leq
% \left\|\widehat{\nu'}\right\|\left\|\widehat{\nu}\right\|
\end{equation}
when $\nu' : C \to E$ is another $*$-homomorphism between unital $C^*$-algebras. 
Let
  $\pi_{k,i}(x,f) = \left(\chi_{k,i}(x)+ \mu_{k,i}(f),\varphi_{k,i}(f)\right)$ be the
  decomposition of $\pi_{k,i} = \pi_{k-1} \circ \pi_{k-2} \circ
  \dots \circ \pi_i : A_i \to A_k$ arising from the fact that the
  $\pi_j$'s, and hence also $\pi_{k,i}$, are regular
  $*$-homomorphisms. Let $z = (x,f) \in A_i$ be a
  self-adjoint element in the unit ball of $A_i$ and let $\epsilon >
  0$. We will show that there is a constant $\lambda \in
  [-1,1]$ such that 
\begin{equation}\label{traceeq2}
\left\|\widehat{\pi_{\infty,i}}(z) - \lambda 1
  \right\| \leq 2\epsilon
\end{equation} 
in $\Aff T(A)$, where $\pi_{\infty,i} : A_i \to A$ is the canonical
embedding going with the inductive limit construction. The desired
  conclusion follows easily from (\ref{traceeq2}) since $i, z$ and
  $\epsilon > 0$ are arbitrary. Observe first that since $B$ has a unique trace state there is a
constant $\lambda \in [-1,1]$ such that $\lim_{k \to \infty}
\left\|\widehat{\varphi_{k,i}}(f) - \lambda 1 \right\| = 0$. Hence
$\left\|\widehat{\varphi_{j,i}}(f) - \lambda 1 \right\| \leq \epsilon$
for all large $j$. Write $\pi_{j,i} (z) = (y,g)$ where $g =
\varphi_{j,i}(f)$. Since $\varphi_{j,i}$ is unital it follows from
(\ref{traceeq1}) that 
\begin{equation}\label{traceeq3}
\left\|\widehat{\mu_{k,j}}(g) - \lambda
  \widehat{\mu_{k,j}}(1)\right\|  = \left\| \widehat{\mu_{k,j} \circ
    \varphi_{j,i}}(f - \lambda 1)\right\|  \leq
\left\|\widehat{\varphi_{j,i}}(f) - \lambda 1 \right\| \leq \epsilon
\end{equation}
 for all $k \geq j$. Since $\chi_{k,i}$ does not increase rank by condition 3) on a
regular $*$-homomorphism, we find that 
\begin{equation}\label{traceeq4}
\lim_{k \to \infty} \left\|\widehat{\chi_{k,j}}(y - \lambda 1)\right\|
  = 0.
\end{equation}
By combining (\ref{traceeq3}) and (\ref{traceeq4}) we find that
\begin{equation*}
\left\|
  \widehat{\chi_{k,j}}(y) + \widehat{\mu_{k,j}}(g) - \lambda 1\right\|
= \left\|\widehat{\chi_{k,j}}(y -\lambda 1) + \widehat{\mu_{k,j}}(g) - \lambda
\widehat{\mu_{k,j}}(1)\right\|
  \leq 2\epsilon
\end{equation*} 
for all large $k$. It follows that (\ref{traceeq2}) holds.
\end{proof}
\end{lemma}

\begin{lemma}\label{RR0}  $B_{\overline{h}}
  \left(\overline{\Gamma}\right)$ has real rank zero.
  \begin{proof} By \cite{BP} we must show that the $C^*$-algebra $A$
    of Proposition~\ref{limit3} has real rank zero. Since $A$ is a
    simple unital recursively subhomogenous $C^*$-algebra in the sense
    of \cite{Ph3} and since $A$ only has one trace state by
    Lemma~\ref{uniquetrace} it suffices by Theorem~4.2 of \cite{Ph3}
    to show that there are non-zero projections in $A$ of arbitrary
    small trace. For this purpose observe that the projection $p$
    constructed in the proof of Lemma~\ref{projex} satisfies the
    inequality
\begin{equation*}
\tau(p) \leq \frac{2n+1}{{\min_k a(k)}}  
\end{equation*}
for all $\tau \in T\left(A(a,b,I,U)\right)$. As we have already used above,
the mixing condition and the fullness of the connecting
$*$-homomorphisms in (\ref{sequence3}) imply that the number $\min_k
a(k)$ becomes arbitrarily large in the building blocks of (\ref{sequence3}). It follows that
$A$ has non-zero projections of arbitrarily small trace.
  \end{proof}
\end{lemma}

\begin{remark} 
My original proof of Lemma \ref{RR0} was more direct, but also more
complicated. I am grateful to N.C. Phillips for his suggestion to use the
result from \cite{Ph3} to shorten the proof. \qed
\end{remark}

We know now that $B_{\overline{h}}
  \left(\overline{\Gamma}\right) \simeq A \otimes \mathbb K$ where $A$
 
\begin{enumerate}
\item[$\bullet$] is the inductive limit of a sequence of unital
  building blocks with unital connecting maps,
\item[$\bullet$] has real rank zero, and
\item[$\bullet$] has a unique trace state.
\end{enumerate}

This allows us to use the work of H. Lin on simple $C^*$-algebras of tracial rank 0.

\begin{lemma}\label{TTR0} $A$ has tracial rank zero in the sense of
  H. Lin, cf. \cite{Lin2}.
\begin{proof} By Proposition 5.4 of \cite{Lin3} $T(A)$ is approximately
  AC. It is easily shown that the invertibles are dense in a unital building
  block and it follows that the same is true in $A$, i.e. $A$ has
  stable rank one. Furthermore, it follows from Lemma \ref{Ktheory1} that
  $K_0(A)$ is unperforated. Hence Theorem 4.15 of \cite{Lin3} shows that
  $A$ has tracial rank zero. 
\end{proof}
\end{lemma}

By an \emph{AH-algebra with no dimension-growth} we mean in the
following the inductive limit of a sequence of $C^*$-algebras 
of the
form
\begin{equation}\label{pseq}
p_1C\left(X_1,M_{k_1}\right)p_1 \to p_2C\left(X_2,M_{k_2}\right)p_2
 \to p_3C\left(X_3,M_{k_3}\right)p_3 \to \cdots,
\end{equation}
where each $X_i$ is a compact metric space, $p_i$ is a projection in
$C\left(X_i,M_{k_i}\right)$ and $\sup_i \dim X_i < \infty$. In the
case where each of the $X_i$'s is a disjoint union of circles the
resulting algebra is called an \emph{AT-algebra.}

\begin{thm}\label{main} Let $B_{\overline{h}}
  \left(\overline{\Gamma}\right)$ be the heteroclinic algebra of the 1-solenoid $\left(\overline{\Gamma},
    \overline{h}\right)$. Then $B_{\overline{h}}
  \left(\overline{\Gamma}\right)  \simeq A \otimes \mathbb K$ where
  $A$ is a unital simple AH-algebra with real rank zero, no dimension
  growth and a unique trace. 
\begin{proof} As pointed out by N.C. Phillips in Proposition 4.7 of
  \cite{Ph2} this follows from the work of Elliott, Gong and Lin, cf.
  \cite{Lin3} and \cite{EG}.  
\end{proof}
\end{thm}

In particular, it follows now, thanks to the work of Elliott and Gong,
\cite{EG}, or of Elliott, Gong and Li, \cite{EGL}, that
$\smash{B_{\overline{h}} \left(\overline{\Gamma}\right)}$ is classified by
K-theory. In fact, since
$K_0\left(B_{\overline{h}}\left(\overline{\Gamma}\right)\right)$ is torsion
  free it follows from the concluding remark in \cite{EGL}
that $B_{\overline{h}} \left(\overline{\Gamma}\right)$ can also be
realized by the inductive limit of a sequence of direct sums of circle
algebras and matrix algebras over so-called dimension drop algebras,
i.e. algebras of the form
\begin{equation*}
\left\{ f \in C\left([0,1], M_k\right) : \ f(0) \in \mathbb C, f(1)
  \in \mathbb C \right\} .
\end{equation*}
Thus $B_{\overline{h}}
  \left(\overline{\Gamma}\right)$ is actually in the class of real
  rank zero algebras classified in Elliotts breakthrough paper,
  \cite{Ell2}. Furthermore, since the dimension drop algebras are only
  needed in the presence of torsion in $K_1$, it will follow from results
  below that the dimension drop algebras can be omitted if and only if
  $\left(\overline{\Gamma}, \overline{h}\right)$ is orientable.

\subsection{On the K-theory of $B_{\overline{h}}\left(\overline{\Gamma}\right)$}

With the detailed information on the connecting maps at hand it is
easy to obtain the following description of the $K$-theory of
$B_{\overline{h}}\left(\overline{\Gamma}\right)$ from the
'stationary' sequence of Lemma \ref{limit3}.

\begin{prop}\label{Kth}  The partially ordered
  group $K_0\left(B_{\overline{h}}\left(\overline{\Gamma}\right)\right)$ is the inductive limit of the sequence
\begin{equation*} 
H \xrightarrow{\chi_* + \mu_* \circ I}
H \xrightarrow{\chi_* + \mu_* \circ I}
H \xrightarrow{\chi_* + \mu_* \circ I} \cdots,
% H \overset{\chi_* + \mu_* \circ I}\to H \overset{\chi_* + \mu_* \circ I}\to
% H \overset{\chi_* + \mu_* \circ I}\to \dots,
\end{equation*}  
where $H = \left\{x \in \mathbb Z^m : \ Ix = Ux \right\}$,
and
\begin{equation*}
  K_1\left(B_{\overline{h}}\left(\overline{\Gamma}\right) \right) 
  \simeq \coker (I-U)  =  \mathbb Z^n /
  (I-U)\left(\mathbb Z^m\right) .
\end{equation*}
\begin{proof} Set $C_k = A\left(a_k,b_k, I,U\right)$ and identify
  $K_0\left(F_{a_k}\right)$ with $\mathbb Z^m$ and
  $K_0\left(F_{b_k}\right)= K_1\left(S F_{b_k}\right)$ with $\mathbb Z^n$ . The conclusion
  regarding $K_0$ follows from Lemma \ref{limit3} and Lemma \ref{Ktheory1} by
  observing that
\begin{equation*}
\begin{xymatrix}{
K_0\left(C_k\right) \ar[d]_-{{\pi_k}_*} \ar[r] & \mathbb Z^m
\ar[d]_-{\chi_* + \mu_* \circ I}
\ar[r]^-{I - U}  & \mathbb Z^n \\
K_0\left(C_{k+1}\right) \ar[r]  & \mathbb Z^m
\ar[r]^-{I - U}  & \mathbb Z^n
}
\end{xymatrix}
\end{equation*}
commutes. Similarly, by using the properties n1)-n5) of the
characteristic functions it follows that also
\begin{equation*}
\begin{xymatrix}{
\mathbb Z^m   \ar[r]^-{I-U}  & \mathbb Z^n \ar@{=}[d]
\ar[r]  & K_1\left(C_k\right) \ar[d]^-{{\pi_k}_*} \\
\mathbb Z^m   \ar[r]_-{I-U}  & \mathbb Z^n 
\ar[r]  & K_1\left(C_{k+1}\right)
}
\end{xymatrix}
\end{equation*}
commutes. This implies the statement concerning $K_1$.
\end{proof}
\end{prop}

Note that
$K_0\left(B_{\overline{h}}\left(\overline{\Gamma}\right)\right)$ is a
stationary partially ordered group, but in a slightly more general sense
 than the usual because $H$ is not (obviously) a simplicial group. For
 this reason we must resort to Theorem \ref{main} in order to conclude
 that $K_0\left(B_{\overline{h}}\left(\overline{\Gamma}\right)\right)$
 is a
dimension group. In the following we focus on $K_1\left(B_{\overline{h}}\left(\overline{\Gamma}\right)\right)$.

Let $A = \left(a_{ij}\right)$ be an $n \times m$ matrix. Let $G_A$ be the unoriented graph
with vertex set $\{1,2, \dots, n\}$ and an edge between $i$ and $j$ if
and only if $i \neq j$ and there is a $k \in \{1,2,\dots, m\}$ such
that $a_{ik} \neq 0$ and $a_{jk} \neq 0$. We say that $A$ is
\emph{irreducible} when $G_A$ is connected.

\begin{lemma}\label{matrix2} Let $n,m \in \mathbb N, n \geq 2$. Let $A$ be an irreducible $n \times m$ integer
  matrix with the property that each
  non-zero column
  has exactly two non-zero entries from the set $\{1,-1\}$ or exactly one
  non-zero entry from the set $\left\{2,-2\right\}$. 

  a) If $A$ contains
  an entry from $\{2,-2\}$ it follows
  that the cokernel of $A : \mathbb Z^m \to \mathbb Z^n$ is isomorphic
  to $\mathbb Z_2$.
  
  b) If each
  non-zero column of $A$
  has exactly two non-zero entries, $1$ and $-1$, it follows
  that the cokernel of $A : \mathbb Z^m \to \mathbb Z^n$ is isomorphic
  to $\mathbb Z$. 
  
\begin{proof} a) We will prove by induction in $n$ that there is a matrix
  $H \in \Gl_n(\mathbb Z)$ such that $HA\left(\mathbb Z^m\right) = 2
  \mathbb Z \oplus \mathbb Z^{n-1}$. It is easy to see that this holds when $n = 2$. Consider then the case $n > 2$. By permuting the
  rows and columns we can arrange that $a_{11} \neq 0$ and $a_{21}
  \neq 0$. (The possibility that each columb only has one non-zero
  entry is ruled out by the irreducibility of $A$.) If $a_{11}$ and $a_{21}$
  have the same sign we subtract the first row from the second;
  otherwise we add the first row to the second. Both
  operations correspond to the multiplication by an element of
  $\Gl_n(\mathbb Z)$. The resulting $n \times m$-matrix
  $A'$ has the form
\begin{equation}\label{rowop}
A' = \left( \begin{matrix} \pm 1 & * \\ 0 & B \end{matrix} \right) ,
\end{equation}
where $*$ is an $1 \times (m-1)$ matrix with entries from
$\{-2,-1,0,1,2\}$, the $0$ is $(n-1) \times 1$ and $B$ is an irreducible
$(n-1) \times (m-1)$ integer matrix for which we can use the induction
hypothesis. In particular $B$ contains an entry from $\{2,-2\}$. Let $H_0 \in \Gl_{n-1}(\mathbb Z)$ be such that $H_0B\left(
  \mathbb Z^{m-1}\right) = 2
  \mathbb Z \oplus \mathbb Z^{n-2}$. Then  
 \begin{equation*}
H = \left( \begin{matrix} 1 & 0 \\ 0 & H_0 \end{matrix} \right) \in
\Gl_n(\mathbb Z)
\end{equation*} 
has the property that $HA'\left(\mathbb Z^{m}\right) = \mathbb Z
\oplus 2\mathbb Z \oplus \mathbb Z^{n-2}$. This completes the induction step and hence
the proof.

b) is proved in the same way: Use induction to show that there is a matrix
  $H \in \Gl_n(\mathbb Z)$ such that $HA\left(\mathbb Z^m\right) = 0
  \oplus \mathbb Z^{n-1}$.
\end{proof}
\end{lemma}

\begin{lemma}\label{matrix} Let $n,m \in \mathbb N, n \geq 2$. Let $A$ be an irreducible $n \times m$ integer
  matrix with the property that each
  non-zero column
  has exactly two non-zero entries from the set $\{1,-1\}$. It follows
  that the cokernel of $A : \mathbb Z^m \to \mathbb Z^n$ is isomorphic
  to $\mathbb Z$ or to $\mathbb Z_2$.
\begin{proof} After appropriate permutations of rows and columns and
  a single permissible row operation we obtain a matrix $A'$ as in
  (\ref{rowop}) such that the matrix $B$ is covered by a) or b) in Lemma \ref{matrix2}. 
\end{proof}
\end{lemma}

\begin{thm}\label{dichotomy} Let $B_{\overline{h}}
  \left(\overline{\Gamma}\right)$ be the heteroclinic algebra of the 1-solenoid $\left(\overline{\Gamma},
    \overline{h}\right)$. Then $K_1\left(B_{\overline{h}}
  \left(\overline{\Gamma}\right)\right)  \simeq \mathbb Z$ when $h
: \Gamma \to \Gamma$ is oriented and $K_1\left(B_{\overline{h}}
  \left(\overline{\Gamma}\right)\right)  \simeq \mathbb Z_2$ when it
is not.
\begin{proof} When $h$ is oriented there is an orientation of the
  edges of $\Gamma$ such that $h^2$ is 
 orientation preserving and then the matrix $U-I$ will have the
property described in b) of Lemma \ref{matrix2} and hence $K_1\left(B_{\overline{h}}
  \left(\overline{\Gamma}\right)\right)  \simeq \mathbb Z$. 

Assume
that $K_1\left(B_{\overline{h}}
  \left(\overline{\Gamma}\right)\right)  \simeq \mathbb Z$. Let $Q$ be
the passages of (\ref{stab2}). Let $e,f$ be two edges of $\Gamma$. The
passages from $Q$ which involve both $e$ and $f$ are of the form
\begin{alignat*}{4}
+\quad & : &\quad& \overset{e}\to \overset{f}\to  &\quad&\text{or} &\quad&\overset{f}\to \overset{e}\to ,
\\
 - \quad & : &\quad&  \overset{e}\to \overset{f}\leftarrow && \text{or} && \overset{e}\leftarrow \overset{f}\to .
\end{alignat*}
% \begin{equation*}
% + \ \ \ : \ \ \ \ \ \overset{e}\to \overset{f}\to  \ \ \ \ \text{or} \ \ \ \
% \overset{f}\to \overset{e}\to ,
% \end{equation*}
% \begin{equation*}
%  - \ \ \ : \ \  \  \ \  \overset{e}\to \overset{f}\leftarrow \ \ \ \ \text{or} \ \
% \ \ \overset{e}\leftarrow \overset{f}\to .
% \end{equation*}
We use these parities to define a labeled unoriented graph $H$ whose
vertexes are the edges of $\Gamma$. Let $e,f \in \mathbb E$. If all passages involving $e$ and $f$ are of parity $+$ draw an edge
between $e$ and $f$ with label $+$. If not, and there is in fact at
least one passage involving both $e$ and $f$, draw an edge between $e$
and $f$ with label $-$. Note that $H$ is connected. A change in the orientation of the edges in
$\Gamma$ will not change the graph $H$ - only its labels. Let
$\mathbb G$ denote the collection of connected subgraphs $H'$ of $H$ with the
property that there is an orientation of the edges in $\Gamma$ such
that all edges of $H'$ are labeled $+$. Then $\mathbb G$ is not empty
and is partially ordered by inclusion. Let $H_{max}$ be a maximal
element. If $H_{max} \neq H$ there must be a loop in $H$ of the form 
\begin{equation*}
\begin{xymatrix}{
e_1 \ar@{-}[r]^-{-}  & e_2 \ar@{-}[r]^-{+} & e_3\ar@{-}[r]^-{+} & e_4
\ar@{-}[r]^-{+} & \dots \dots \ar@{-}[r]^-{+} & e_1 
}
\end{xymatrix}
\end{equation*}
Writing down what this means for $U-I$ it is easy to see that we can obtain a matrix $B$
from $U-I$ by permutations and addition and/or subtraction of rows such that 
\begin{equation*}
B = \left( \begin{matrix}  1 & * & * & * & \dots & * \\ 0 & 1 & * & *
    & \dots & * \\ 
    \vdots & \vdots & \vdots & \vdots & \ddots & \vdots
    \\ 0 & 0 & 0 & 0 & \hdots & B_0  \end{matrix} \right) 
\end{equation*}
where $B_0$ is an irreducible matrix for which a) Lemma \ref{matrix2}
applies. It follows that $K_1\left(B_{\overline{h}}
  \left(\overline{\Gamma}\right)\right)  \simeq \mathbb Z_2$ which
contradicts the assumption, proving that $H_{max} = H$. 

Give the edges of $\Gamma$ an orientation such that all the passages,
$Q$, are of type $+$. It
follows that for all large enough $i$, $h^i : J \to \Gamma$ will
either be orientation preserving or orientation reversing. If $h^i: J
\to \Gamma $
and $h^{i+1} : J \to \Gamma$ are either both orientation preserving or
both orientation reversing, $h : \Gamma \to \Gamma$ will be
orientation preserving. Otherwise $h$ must be orientation reversing. It follows that $h$ is either positively or
negatively oriented.

We have shown that $(\Gamma, h)$ is oriented if and only if
$K_1\left(B_{\overline{h}}\left(\overline{\Gamma}\right)\right) \simeq
\mathbb Z$. It follows then from Lemma \ref{matrix} that 
$K_1\left(B_{\overline{h}}\left(\overline{\Gamma}\right)\right) \simeq
\mathbb Z_2$ when $\left(\Gamma ,h\right)$ is not oriented. 
\end{proof}
\end{thm}

\begin{thm}\label{orient2} Assume that the pre-solenoid $h : \Gamma \to
  \Gamma$ is oriented. It follows that the heteroclinic
  algebra $B_{\overline{h}}\left(\overline{\Gamma}\right)$ is an AT-algebra.
\begin{proof} By Theorem \ref{dichotomy} and Proposition \ref{Kth} the
  K-theory of $B_{\overline{h}}\left(\overline{\Gamma}\right)$ is
  torsion-free. Therefore the conclusion follows from Theorem
  \ref{main} and Theorem \cite{EG}.
 \end{proof}
\end{thm}

We can now give the proof of Theorem \ref{orientthm}:
\begin{proof} Let $\left(\Gamma, h\right)$ be a negatively orientable
  pre-solenoid. There is then a negatively oriented pre-solenoid 
  $\left(\Gamma',h'\right)$ such that
  $\left(\overline{\Gamma},\overline{h}\right)$ is conjugate to
  $\left(\overline{\Gamma'},\overline{h'}\right)$. Then
  $B_{\overline{h}}\left(\overline{\Gamma}\right)$ is $*$-isomorphic
  to $B_{\overline{h'}}\left(\overline{\Gamma'}\right)$ and hence the
  two algebras have the same $K_1$-group. It follows therefore from
  Theorem \ref{dichotomy} that $\left(\Gamma,h\right)$ is
  oriented. Assume to reach a contradiction that
  $\left(\Gamma,h\right)$ is positively oriented. Give $\Gamma$ an
  orientation such that $h$ is positively oriented with respect to it,
  and give
  $\Gamma'$ an orientation such that $h'$ is negatively oriented with
  respect to it. Then both $h^2: \Gamma \to \Gamma$ and $h'^2 :
  \Gamma' \to \Gamma'$ preserve the orientation and hence
  $\overline{\Gamma} $ and $\overline{\Gamma'}$ 
  are both oriented matchbox manifolds in the sense of 
  \cite{AO}. Since $\left(\Gamma, h\right)$ and
  $\left(\Gamma',h'\right)$ are conjugate there is a homeomorphism
  $\mu : \overline{\Gamma} \to \overline{\Gamma'}$ such that 
\begin{equation}\label{orinv}
\overline{h'}\circ
  \mu = \mu \circ \overline{h}.
\end{equation} 
The 'matchboxes' in $\overline{\Gamma}$ and $\overline{\Gamma'}$, in
the sense of \cite{AO}, are all the product of $[-1,1]$ and the Cantor
set. Since there are no isolated points in the Cantor set it is easy
to see that $\mu$ must either preserve or
  invert the orientation. None of the two possibilities are consistent
  with (\ref{orinv}) since $\overline{h'}$ inverts the orientation of
  $\Gamma'$ while $\overline{h}$ preserves the orientation of
  $\overline{\Gamma}$. This contradiction shows that $\left(\Gamma,
    h\right)$ must be negatively oriented.

The case where $\left(\Gamma, h\right)$ is a positively orientable
  pre-solenoid is handled in exactly the same way.
\end{proof}

\section{The heteroclinic algebra of the inverse of a 1-solenoid}\label{invhet}

 We
seek here to clarify the structure of the heteroclinic algebra of the
inverse of $\overline{h}$, i.e. we seek to describe
$B_{\overline{h}^{-1}}\left(\overline{\Gamma}\right)$. Since $\left(\overline{\Gamma}, \overline{h}\right)$ is a Smale space,
$B_{\overline{h}^{-1}}\left(\overline{\Gamma}\right)$ is the same
algebra as the stabilization of the 'unstable algebra' of Putnam,
cf. \cite{Pu} and Section 4.4 of \cite{Th1}. However, the
technicalities required to show that a 1-solenoid has a Smale space
structure are quite involved and tedious, and do not seem to exist in
the litterature. For this reason we continue to ignore this additional
structure since it is not really needed. Only in Section \ref{homo}
below will it be necessary to appeal to the Smale space structure.

\subsection{Simplicity and stability}\label{secsimstab}

In this section we show that
$B_{\overline{h}^{-1}}\left(\overline{\Gamma}\right)$ is simple and stable.

\begin{lemma}\label{wagtop} Let $x \in W_{\overline{\Gamma},
    \overline{h}^{-1}}$. The sets
$$
U(x,n) = \left\{ y \in \overline{\Gamma} : \ y_n = x_n \right\}, \ n
\in \mathbb N, \ 
$$
form a neigborhood base of $x$ for the Wagoner
topology of $W_{\overline{\Gamma}, \overline{h}^{-1}}$.
\begin{proof} Consider the $\kappa > 0$ and the $l \in \mathbb N$ from
  Lemma 2.9 of \cite{Y1}. When $y \in U(x,n)$, $z \in
  \overline{\Gamma}$ and
  $D\left(\overline{h}^m(z),\overline{h}^m(y)\right) < \kappa, \ m \geq
  - l - n$, it follows that
  $d\left(h^k\left(z_{l+n}\right),h^k\left(y_{l+n}\right)\right) <
  \kappa$ for all $k \in \mathbb N$. By Lemma 2.9 of \cite{Y1} this
  implies that $y_n = h^l\left(y_{l+n}\right) =h^l\left(z_{l+n}\right)
  = z_n$, proving that
$$
\left\{ z \in \overline{\Gamma} : \
  D\left(\overline{h}^m(z),\overline{h}^m(y)\right) < \kappa, \ m \geq
  - l - n \right\} \subseteq U(x,n) .
$$
By the definition of the Wagoner topology of $W_{\overline{\Gamma},
  \overline{h}^{-1}}$ this implies that $U(x,n)$ is open in $W_{\overline{\Gamma},
  \overline{h}^{-1}}$, cf. Lemma 4.6 \cite{Th1}. To prove that the $U(x,n)$'s form a base consider an element
$$
\Omega = \left\{ y \in \overline{\Gamma} : \ D\left(\overline{h}^m(y),
    \overline{h}^m(x)\right) < \epsilon , \ m \geq k \right\}
$$ 
of the base described in Lemma 4.6 of \cite{Th1}. Choose $L \in
\mathbb N$ such that 
$$
\sum_{j=L+1}^{\infty} 2^{-j} \sup_{a,b \in \Gamma} d(a,b) <
\epsilon.
$$ 
Then $U\left(x,|k|+L\right) \subseteq \Omega$. 
 \end{proof} 
\end{lemma}

\begin{lemma}\label{opendense} Let $x,y \in \overline{\Gamma}$ be
  post-periodic for $\overline{h}^{-1}$. Assume that there is an $N
  \in \mathbb N$ such that 
$$
\forall i \geq N \ \exists e_i \in \mathbb E : \ \left[x_i,y_i\right]
\subseteq \inter e_i 
.
$$
It follows that $x$ is locally conjugate to $y$.
\begin{proof} Let $x' \in U(x,N)$. Note that there is an edge $f_{N+1}$ such
  that $x'_{N+1} \in \inter f_{N+1} $ and that there is a unique interval
  $\left[x'_{N+1},y'_{N+1}\right] \subseteq \inter f_{N+1}$ such that
  $h\left(\left[x'_{N+1},y'_{N+1}\right]\right) =
  \left[x_N,y_N\right]$. Continuing recursively we obtain edges $f_i,
  i \geq N+1$, and intervals $\left[x'_{i},y'_{i}\right] \subseteq \inter f_{i}$ such that
  $h\left(\left[x'_{i},y'_{i}\right]\right) =
  \left[x'_{i-1},y'_{i-1}\right]$ for all $i \geq N+2$. We define
  $\chi : U(x,N) \to U(y,N)$ such that
  $$
\chi(x')_i = \begin{cases} y_i, \ i \leq N \\ y'_{i}, \ i \geq N+1 .
\end{cases} 
$$
Note that $\chi(x) =y$. By symmetry we can define $\mu : U(y,N) \to
U(x,N)$ in a similar way and
it is then clear that $\chi \circ \mu = \id_{U(y,N)}$ while $\mu \circ
\chi = \id_{U(x,N)}$. Furthermore, it follows easily from Lemma
\ref{wagtop} that both $\chi$ and $\mu$ are continuous for the Wagoner
topology. Finally, it follows from the expansion condition that there
are constants $C' > 0$ and $\lambda >1$ such that $d\left(x'_i, \chi(x')_i\right) \leq C'\lambda^{-i}$
for all $x' \in U(x,N)$ and all $i \geq N$. In particular, 
$$
\lim_{ n \to \infty} \sup_{x' \in U(x,N)} D
\left(\overline{h}^{-n}(x'),\overline{h}^{-n}\left(\chi(x')\right)\right)
= 0,
$$
completing the proof that $\left(U(x,N),U(y,N),\chi\right)$ is a local
conjugacy from $x$ to $y$.
\end{proof}
\end{lemma}

\begin{lemma}\label{DSB} Let $x,y \in \overline{\Gamma}$ be
  post-periodic for $\overline{h}^{-1}$. Assume that $x_N \notin
  \mathbb V$ for some $N$. For every $n \in \mathbb N$ there is an
  element $z \in U(y,n)$ such that $x$ is locally conjugate to $z$.
\begin{proof} It follows from the mixing condition that there is a $k
  \in \mathbb N$ such that $h^k\left(\inter e\right) = \Gamma$ for all
  $e \in \mathbb E$. Choose $e\in \mathbb E$ such that $x_m \in
  \inter e$, where $m = \max \{n,N\} + k$. There is then an element
  $z_{m}\in \inter e$ such that $h^k\left(z_m\right) = y_{m-k}$. We
  can then define $z_i, i \geq m+1$, recursively by the requirement
  that $\left[z_i,x_i\right] \cap \mathbb V = \emptyset$ and
  $h\left(z_i\right) = z_{i-1}$. Set $z_j = h^{m-j}(z_m), j \leq
  m$. Then $z = \left(z_j\right)_{j=0}^{\infty} \in U(y,n)$ and it
  follows from Lemma \ref{opendense} that $z$ is locally conjugate to $x$.
\end{proof}
\end{lemma}

\begin{remark}\label{notflatremark} With future investigations in mind
  I want to point out that the last three lemmas, unlike the following, do not depend on
  the flattening condition $\epsilon$). 
\end{remark}

Thanks to Lemma \ref{invlemma1}, in order to study the structure of
$B_{\overline{h}^{-1}}\left(\overline{\Gamma}\right)$ we may
substitute $h$ with any positive power of $h$. As in \cite{Y4} we can
use this together with a result of Williams, 5.2 of \cite{Wi2}, to
restrict our attention to the case where $\Gamma$ is the wedge of $n$ circles, each of which is an edge in
$\Gamma$, such that the common point of the circles is the only vertex
in $\Gamma$. To emphasize this assumption we denote the graph by
$\Gamma_c$.
Let $e_i, i=1,2,\dots, n$, be the edges
of $\Gamma_c$ and let $v$ be the unique vertex in $\Gamma_c$. By the
flattening condition
$\epsilon$) and Lemma \ref{invlemma1} we can also assume that
there is an open neighborhood $U$ of $v$ in $\Gamma_c$ such that $h(U)
\simeq \ ]-1,1[$ and that $h(e_i) = \Gamma_c$ for all $i$. There are numbers $a_{ij} \in
\mathbb N$ and disjoint intervals $I_{ij}^l, l = 1,2, \dots, a_{ij}$, such that
\begin{equation}\label{union}
\inter \left(e_i\right) \cap h^{ -1}\left(\inter \left(e_j\right)\right) =
\bigcup_{l=1}^{a_{ij}} I^l_{ij} .
\end{equation}
Set $\mathcal J = \left\{ I^l_{ij}: \ l = 1,2, \dots, a_{ij}, \ i,j =
  1,2, \dots, n\right\}$. Since we assume that $h(U) \simeq ]-1,1[$ for all small neighborhoods
$U$ of $v$ we can renumber the edges of $\Gamma_c$ and the intervals in
$\mathcal J$ to arrange that for
some $j \in \{1,n\}$ it holds that
\begin{equation}\label{firstcond}
h(U) \subseteq \overline{I_{11}^1} \cup \overline{I_{jj}^{a_{jj}}}
\end{equation}
or
\begin{equation*}
h(U) \subseteq \overline{I_{1j}^1} \cup \overline{I_{j1}^{a_{j1}}}
\end{equation*}
If the latter occurs and $j \neq 1$, or if $1=j$ and
$h\left({I^1_{11}}\right) \supseteq {I^{a_{jj}}_{jj}}$, we substitute $h$ with $h^2$, as we
may by Lemma \ref{invlemma1}, to ensure that (\ref{firstcond}) holds
and that
\begin{equation}\label{secondcond2}
h\left(\overline{I^1_{11}}\right) \supseteq \overline{I^1_{11}}, \ \
h\left(\overline{I^{a_{jj}}_{jj}}\right) \supseteq \overline{I^{a_{jj}}_{jj}}.
\end{equation}

\begin{lemma}\label{thevertex} Let $x, y \in \overline{\Gamma_c}$ be
  such that $x_i = v$, $y_i \neq v$ for all $i$, and $\lim_{i \to
    \infty} d\left(x_i, y_i\right) = 0$. It follows that $x$ and $y$
  are locally conjugate in $W_{\overline{\Gamma_c}, \overline{h}^{-1}}$.
\begin{proof} Since (\ref{firstcond}) and (\ref{secondcond2}) hold there is an $N \in \mathbb N$
  such that $y_i \in I^1_{11}$ for all $i \geq N $ or $y_i \in I^{a_{jj}}_{jj}$
  for all $i \geq N$. Without loss of generality we assume that $y_i
  \in I^1_{11}$ for all $i \geq N $. We define $\chi : U(y,N) \to
  U(x,N)$ in the following way. Set $\chi(y) = x$. When $z \in
  U(y,N) \backslash \{y\}$, set $i_z = \max \left\{ j \in \mathbb N: \
    z_i = y_i, \ i \leq j \right\}$. Then $z_{i_z+1} \in \inter J_1$ for
  some $J_1 \in \mathcal J$ and there is a unique point $x'_{i_z+1} \in \overline{J_1}$ such that
  $h\left(x'_{i_z+1}\right) = x_{i_z} = v$ and $h\left(\left[z_{i_z+1},
      x'_{i_z+1}\right]\right) = \left[z_{i_z} ,x_{i_z}\right]
  \subseteq \overline{I^1_{11}}$. We can then continue
  recursively: Choose $x'_{i_z+2} \in J_2 \in \mathcal J$ such that
  $h\left(x'_{i_z+2}\right) = x'_{i_z+1}$ and $h\left(\left[z_{i_z+2},
      x'_{i_z+2}\right]\right) = \left[z_{i_z+1} ,x'_{i_z+1}\right]
  \subseteq \overline{J_1}$ and so
  on. We define 
$$
\chi(z)_j = \begin{cases} x_j = v, & j \leq i_z, \\  x'_j, \ & j \geq
  i_z+1. \end{cases}
$$
Since $\chi(z)_j$ only depends on $z_j$ it follows from Lemma
\ref{wagtop} that $\chi$ is continuous for the Wagoner topology. Note
that the expansion condition ensures that
\begin{equation}\label{firstexp}
\lim_{n \to \infty} \sup_{z \in U(x,N)}
D\left(\overline{h}^{-n}\left(\lambda(z)\right),
  \overline{h}^{-n}(z)\right) = 0.
\end{equation}
We define a map $\mu : U(x,N) \to U(y,N-1)$ in a similar way: Set
$\mu(x) = y$. For $z \in U(x,N) \backslash \{x\}$, set $i_z = \max
\left\{ j \in \mathbb N: \ z_i = x_i = v, \ i \leq j\right\}$. Since
(\ref{firstcond}) holds there is a unique pair
$\left(J_1',y'_{i_z+1}\right) \in \mathcal J \times \Gamma_c$ such that
$z_{i_z+1} \in \overline{J'_1}, \ y'_{i_z+1} \in J'_1$, $h^2\left(y'_{i_z +1}\right) =
y_{i_z-1}$ and $h^2\left(\left[y'_{i_z+1},
    z_{i_z+1}\right]\right) = \left[z_{i_z-1},y_{i_z-1}\right]
\subseteq \overline{I^1_{11}}$. There
is then a unique $J'_2 \in \mathcal J$ such that $z_{i_z+2} \in J'_2$ and a unique
element $y'_{i_z+2} \in J'_2$ such that $h\left(y'_{i_z +2}\right) =
y'_{i_z+1}$ and $h\left(\left[y'_{i_z+2},
    z_{i_z+2}\right]\right) =
\left[z_{i_z+1},y'_{i_z+1}\right] \subseteq \overline{J'_1}$. Again we can continue and we define      
$$
\mu(z)_j = \begin{cases} y_j , & j \leq i_z -1 \\
  h\left(y'_{i_z+1}\right) , & j = i_z, \\  y'_j, \ & j \geq i_z+1 . \end{cases}
$$
Since $\mu(z)_j$ only depends on $z_{j+1}$ it follows from Lemma
\ref{wagtop} that $\mu$ is continuous for the Wagoner topology, and as
above the expansion condition ensures that
\begin{equation}\label{secondexp}
\lim_{n \to \infty} \sup_{z \in U(y,N)}
D\left(\overline{h}^{-n}\left(\mu(z)\right),
  \overline{h}^{-n}(z)\right) = 0.
\end{equation}
By combining (\ref{firstexp}) and (\ref{secondexp}) we find that
$$
\lim_{n \to \infty} \sup_{z \in U(y,N)}
D\left(\overline{h}^{-n}\left(\mu \circ \lambda(z)\right),
    \overline{h}^{-n}(z)\right) = 0.
$$
The expansiveness of $\overline{h}$ on $\overline{\Gamma_c}$ and the
definition of the Wagoner topology therefore imply that $\mu \circ
\lambda(z) = z$ for all $z$ in a Wagoner-neighborhood of $y$. Since we
also have that
$$
\lim_{n \to \infty} \sup_{z \in U(x,N) \cap \mu^{-1}\left(U(y,N)\right)}
D\left(\overline{h}^{-n}\left( \lambda \circ \mu(z)\right),
    \overline{h}^{-n}(z)\right) = 0,
$$
we conclude in the same way that $\lambda \circ \mu(z) = z$ for all
$z$ in a Wagoner-neighborhood of $x$. It follows that $\lambda$, suitably
restricted, is a local conjugacy from $y$ to $x$.
\end{proof}
\end{lemma}

\begin{lemma}\label{simstab} The heteroclinic algebra
  $B_{\overline{h}^{-1}}\left(\overline{\Gamma_c}\right)$ of
  $\left(\overline{\Gamma_c}, \overline{h}^{-1}\right)$ is stable and
  simple.
\begin{proof} It follows from Lemma \ref{DSB} and Lemma
  \ref{thevertex} that every element of
  $W_{\overline{\Gamma_c},\overline{h}^{-1}}$ is locally conjugate to an
  element of any given open set. It follows
  then from Proposition 4.7 of \cite{Re} that
  $B_{\overline{h}^{-1}}\left(\overline{\Gamma_c}\right)$ is
  simple. That $B_{\overline{h}^{-1}}\left(\overline{\Gamma_c}\right)$
  is stable follows also from Lemma \ref{DSB} and Lemma
  \ref{thevertex} by the method used to
  prove Lemma 5.10 of \cite{Th1}; specifically the two lemmas show that
  for any compact subset $K$ of $W_{\overline{\Gamma_c},
    \overline{h}^{-1}}$ there are periodic points and local
  conjugacies such that a),b),c) and d) from the proof of Lemma 5.10
  of \cite{Th1} are all satisfied, and then the results of Hjelmborg
  and Rørdam, \cite{HR}, imply that
  $B_{\overline{h}^{-1}}\left(\overline{\Gamma_c}\right)$ is stable.
\end{proof}
\end{lemma}

\subsection{A free minimal action by the infinite dihedral group}\label{dihedral}

We retain in this section the assumptions about $\left(\Gamma_c,
  h\right)$ made in the previous section. In particular, $\Gamma_c$ is a
wedge of circles. 
For each $e \in
\mathbb E$ choose homeomorphisms $\psi^e_{\pm} :
[0,1] \to e$ such that $\psi^e_+(t) = \psi^e_-(1-t), t \in [0,1]$. We will consider $e$ as
an oriented edge with an orientation defined such that $\psi^e_+$ is
orientation-preserving and hence such that $\psi^e_-$ is
orientation-reversing. The orientations we choose are basically completely
arbitrary, but for notational convenience we will ensure that
\begin{equation}\label{eventually}
\overline{I^1_{11}} =
\psi^{e_1}_-\left(\left[0,\frac{1}{4}\right]\right) \ \text{and} \
\overline{I^{a_{jj}}_{jj}} =  \psi^{e_j}_+\left(\left[0,\frac{1}{4}\right]\right). 
\end{equation}
 Set
$$
J = \psi^{e_1}_-\left(\left[0,\frac{1}{4}\right]\right) \cup
\psi^{e_j}_+\left(\left[0,\frac{1}{4}\right]\right) .
$$
Then
\begin{equation}\label{conda} 
h \ \text{is injective on} \  J \ \text{and} \ h(J) \supseteq J.
\end{equation}

For each $n = 0,1,2, \dots$, let $\mathcal I^n$ denote the connected components
of $\Gamma_c \backslash h^{-n}(v)$. In particular, $\mathcal I^0$ consists of the
interiors, $\inter e$, of the edges $e$ in $\Gamma_c$, and $\mathcal I^1
= \mathcal J$. Note that with
$C$ and $\lambda$ the constants from the expansion condition we have
the estimate
\begin{equation}\label{diam} 
\diam \sigma \leq K C \lambda^{-n}
\end{equation}
for all $\sigma \in \mathcal I^n$, where $K = \diam \Gamma_c$ is the diameter of $\Gamma_c$.

Let $\Sigma$ denote the collection of sequences
$$
\sigma = \left( \pm, \sigma_0,\sigma_1,\sigma_2,\sigma_3, \dots \right),    
$$
where $\sigma_j$ is an interval in $\Gamma_c$ with $\inter \sigma_j \in
\mathcal I^j$ and
$h\left(\sigma_{j+1}\right) = \sigma_j$ for all $j$. The
first entry, $+$ or $-$, in $\sigma$ will be called \emph{the sign of
  $\sigma$} and will be denoted by $sign(\sigma)$. We consider $\Sigma$ as a
compact Hausdorff space in the product topology. 
%{\begin{lemma}\label{bgn1} There are exactly 4 elements $\sigma$ of
%$\Sigma$  such that $v \in \sigma_j$ for all $j = 0,1,2, \dots$.\end{lemma}}

%\begin{lemma}\label{bgn2} There are exactly 2 elements $\sigma$ of $\Sigma$
 % such that $v \in \sigma_j$ for all $j = 0,1,2, \dots$, and such that
 % $\lim_{t \uparrow 1} \left(h|_{\sigma_1}\right)^{-1} \circ \psi^{\sigma_0}_{sign(\sigma)}(t) = v$.

Note that for all $i \geq 1$, the set
$$
\left\{ I \in \mathcal I^i : \ v \in \overline{I} \  \subseteq J \right\} 
$$
consists of exactly two elements, say $J^+_i$ and $J^-_i$. We choose
the superscripts such that $J^+_i \subseteq
\psi^{e_j}_+\left(\left[0,\frac{1}{4}\right]\right)$ and $J^-_i
\subseteq \psi^{e_1}_-\left(\left[0,\frac{1}{4}\right]\right)$. It
follows then from (\ref{secondcond2}) that
$h\left(J^{\pm}_{i+1}\right) = J^{\pm}_i$. Set $a^-_0 =
h\left(I^1_{11}\right) = \inter e_1$ and
$$
a^-_i = J^-_i, \ i \geq 1 .
$$
Alternatively $a^-_i \in \mathcal I^i$ is determined, for all $i$, by
the requirement that $a^-_i \cap
\psi^{e_1}_-\left(\left[0,\frac{1}{4}\right]\right)$ is a non-degenerate
interval which contains $v$. Define an element $a^- \in \Sigma$ such that
$$
a^- = \left(sign\left(a^-\right),a^-_0,a^-_1,a^-_2, \dots \right) \in
\Sigma,
$$
where $sign\left(a^-\right)$ is determined by the condition that
$$
\lim_{t \uparrow 1} \left( h|_{a^-_1} \right)^{-1}\left(\psi^{a^-_0}_{sign\left(a^-\right)}(t)\right) = v .
$$
Similarly we can define an element 
$$
b^- = \left(sign\left(b^-\right),b^-_0,b^-_1,b^-_2, \dots \right)
$$
in $\Sigma$ by the
requirements that $b^-_0 = h\left(I^{a_{jj}}_{jj}\right) = \inter e_j$,
$$
b^-_i = J^+_i, \ i \geq 1,
$$
and
$$
\lim_{t \uparrow 1} \left(h|_{b^-_1}\right)^{-1} \left(\psi^{b^-_1}_{sign\left(b^-\right)}(t)\right) =
v .
$$

It is easy to check that $a^-,b^- \in \Sigma$ have the following properties:
\begin{lemma}\label{q1} 
\begin{enumerate}
\item[i)] $a^-_0 = e_1, \ b^-_0 = e_j, \ sign(a^-) = +$ and $sign(b^-) =
  - \ $.
\item[ii)] $\lim_{t \uparrow 1}
  \left(h^k|_{a^-_k}\right)^{-1}\left(\psi^{e_1}_+(t)\right) = v$
  and $\lim_{t \uparrow 1}
  \left(h^k|_{b^-_k}\right)^{-1}\left(\psi^{e_j}_-(t)\right) = v$
  for all $k \in \mathbb N$.
\item[iii)] If $\sigma \in \Sigma$ is an element such that $v \in
  \sigma_i$ for all $i$ and 
$$
\lim_{t \uparrow 1}
  \left(h|_{\sigma_1}\right)^{-1}\left(\psi^{\sigma_0}_{sign(\sigma)}(t)\right) = v,
$$ 
then $\sigma \in \left\{a^-,b^-\right\}$.
\end{enumerate}
%{\begin{proof} i) It follows from (\ref{condb}) and (\ref{condc}) that for each $k \in
%  \mathbb N$, there is an $\epsilon_k > 0$ such that $h^k :
%  \psi^{e}_-\left(\left[0,\epsilon_k\right]\right) \to 
%  \psi^e_-\left(\left[0,\frac{1}{4}\right]\right)$ is a
%  homeomorphism. Since $v \in J^-_k \subseteq
%  \psi^e_-\left(\left[0,\frac{1}{4}\right]\right)$ (for $k$ large
%  enough) this implies that $a^-_0 = e$ and that $\lim_{t \uparrow 1}
%  \left(h^k|_{a^-_k}\right)^{-1}\left(\psi^{e}_+(t)\right) = \lim_{t \downarrow 0}
%  \left(h^k|_{a^-_k}\right)^{-1}\left(\psi^{e}_-(t)\right) = v$; first for
%  all large $k$ and then for all $k$ because $v$ is fixed by $h$. In
%  particular, $sign(a^-) = +$ since $\lim_{t \uparrow 1}
%  \left(h|_{a^-_1}\right)^{-1}\left(\psi^{e}_+(t)\right) = v$. This proves i) and ii) as far as $a^-$ is concerned. The arguments concerning
%b^-$ are the same.
%To prove iii) observe that if $v \in
 % \sigma_i$ for all $i$, then $v \in \sigma_i \subseteq J$ for all
 % large $i$ and hence $\sigma \in \left\{a^-,b^-\right\}$ by construction.
%\end{proof}}
\end{lemma}

Let $S : \Sigma \to
\Sigma$ be the continuous involution defined such that
$$  
S\left( \pm , \sigma_0,\sigma_1,\sigma_2,\sigma_3, \dots \right) = \left( \mp, \sigma_0,\sigma_1,\sigma_2,\sigma_3, \dots \right).
$$

Set
 $a^+ = S(a^-), b^+ = S(b^-)$. It follows from Lemma \ref{q1} that
$a^+,b^+$ are the only elements $\sigma$ in $\Sigma$
  such that $v \in \sigma_j$ for all $j = 0,1,2, \dots$, and 
  $\lim_{s \downarrow 0} \left(h|_{\sigma_1}\right)^{-1} \circ
  \psi^{\sigma_0}_{sign(\sigma)}(s) = v$. Let $\sigma  \in
\Sigma \backslash \left\{a^-,b^-\right\}$. By Lemma \ref{q1} 
\begin{equation}\label{eqq1}
j_{\sigma} = \min \left\{ k \in
  \mathbb N : \ \lim_{t \uparrow 1}
  \left(h^k|_{\sigma_k}
  \right)^{-1}\left(\psi^{\sigma_0}_{sign(\sigma)}(t)\right) \neq v
\right\}
\end{equation}
is finite.
Then $\lim_{t \uparrow 1}
  \left(h^k|_{\sigma_k}
  \right)^{-1}\left(\psi^{\sigma_0}_{sign(\sigma)}(t)\right) \neq v$
  when $k \geq j_{\sigma}$ and there is a unique element $\sigma'_k \in
  \mathcal I^k$ such that 
\begin{equation}\label{eqq2}
  \lim_{t \uparrow 1} \left(h^k|_{\sigma_k}
\right)^{-1}\left(\psi^{\sigma_0}_{sign(\sigma)}(t)\right) = \sigma_k \cap \sigma'_k.
\end{equation}
Set $\sigma' =
\left(sign(\sigma'),\sigma'_0,\sigma'_1,\sigma'_2, \dots\right)$, where $\sigma'_i, i <
j_{\sigma}$, is determined by the condition that $\sigma'_i =
h^{j_{\sigma}-i}\left(\sigma'_{j_{\sigma}}\right)$ and $sign(\sigma')$ is determined by the condition that
\begin{equation}\label{eq103}
\lim_{t \uparrow 1} \left(h^i|_{\sigma_i}\right)^{-1} \left(
\psi^{\sigma_0}_{sign(\sigma)}(t)\right) = \lim_{s \downarrow 0} \left(h^i|_{\sigma'_i}\right)^{-1} \left(
\psi^{\sigma'_0}_{sign(\sigma')}(s)\right) 
\end{equation}
for all $i \geq 1$. We define a map $\varphi :
\Sigma \to \Sigma$ such that
$$
\varphi(\sigma) = \begin{cases} \sigma', & \ \text{when} \ \sigma \notin
  \left\{a^-,b^-\right\} \\ b^+, & \ \text{when} \ \sigma = a^- \\ a^+,
& \ \text{when} \ \sigma = b^- \end{cases}
$$

\begin{lemma}\label{bgn3} $\varphi$ is a homeomorphism of $\Sigma$ and
  $S \circ \varphi = \varphi^{-1} \circ S$.
\begin{proof} The continuity of $\varphi$ at points different from
  $a^-$ and $b^-$ is obvious. Continuity at the two exceptional points
  requires a little care; if $\sigma = \left(sign(\sigma),
    \sigma_0,\sigma_1,\sigma_2, \dots \right) \in \Sigma$ agrees with
  $b^-$ on more than the first $2$ coordinates the $j_{\sigma}$ of
  (\ref{eqq1}) exceeds $2$. Note that $\sigma_0 = b^-_0 =
  e_j$ and $sign(\sigma) = sign\left(b^-\right) = -$ by Lemma
  \ref{q1}. Therefore (\ref{firstcond}) implies that
$$
\lim_{t \uparrow 1} \left(h^{j_{\sigma}}|_{\sigma_{j_{\sigma}}}\right)^{-1} \circ
\psi^{e_j}_-(t) = \lim_{s \downarrow 0} \left(h^{j_{\sigma}}|_{\sigma'_{j_{\sigma}}}\right)^{-1} \circ
\psi^{e_1}_-(s),
$$
when $\sigma'_{j_{\sigma}} \in \mathcal I^{j_{\sigma}}$ is determined by (\ref{eqq2}). We
conclude from this that $sign(\sigma') = - =
sign\left(a^+\right)$, cf. (\ref{eq103}). Furthermore, $\sigma'_i =
h^{j_{\sigma}-i}\left(\sigma'_{j_{\sigma}}\right) = a^+_i$ when
 $i
\leq j_{\sigma} - 2$. Thus we see that $\sigma'_i = a^+_i$ for all $i \leq N$ when
$sign(\sigma) = sign\left(b^-\right)$ and $\sigma_i = b^-_i$ for all $i \leq N + 2$. This proves the
continuity of $\varphi$ at $b^-$. The continuity at $a^-$ follows in a
similar way. To prove
  that $\varphi$ is a bijection, we construct the inverse: Let $\sigma
  = \left( sign(\sigma), \sigma_0,\sigma_1,\sigma_2, \dots \right) \in \Sigma
  \backslash \left\{a^+,b^+\right\}$. We can then consider
$$
k_{\sigma} = \min \left\{ k \in
  \mathbb N : \ \lim_{s \downarrow 0}
  \left(h^k|_{\sigma_k}
  \right)^{-1}\left(\psi^{\sigma_0}_{sign(\sigma)}(s)\right) \neq v
\right\}.
$$
There is a unique element $\sigma'_{k_{\sigma}} \in \mathcal I^{k_{\sigma}}$ such that 
$$
\sigma_{k_{\sigma}} \cap \sigma'_{k_{\sigma}}
=  \lim_{s \downarrow 0} \left(h^{k_{\sigma}}|_{\sigma_{k_{\sigma}}}
\right)^{-1}\left(\psi^{\sigma_0}_{sign(\sigma)}(s)\right).
$$ 
Set $\sigma' =
\left(sign(\sigma'),\sigma'_0,\sigma'_1,\sigma'_2, \dots\right)$, where $\sigma'_i, i <
k_{\sigma}$, is determined by the condition that $\sigma'_i =
h^{k_{\sigma}-i}\left(\sigma'_{k_{\sigma}}\right)$ and $\sigma'_i, i \geq k_{\sigma}$, by the condition that
$$
\sigma_i \cap \sigma'_i
= \lim_{s \downarrow 0}  \left(h^i|_{\sigma_i}
\right)^{-1}\left(\psi^{\sigma_0}_{sign(\sigma)}(s)\right).
$$ 
Finally,
$sign(\sigma')$ is determined by the condition that
$$
\lim_{t \uparrow 1} \left(h^{k_{\sigma}}|_{\sigma'_{k_{\sigma}}}\right)^{-1} \circ
\psi^{\sigma'_0}_{sign(\sigma')}(t) = \lim_{s \downarrow 0} \left(h^{k_{\sigma}}|_{\sigma_{k_{\sigma}}}\right)^{-1} \circ
\psi^{\sigma_0}_{sign(\sigma)}(s) .
$$ 
 We can then define a map $\phi :
\Sigma \to \Sigma$ such that
$$
\phi(\sigma) = \begin{cases} \sigma', & \ \text{when} \ \sigma \notin
  \left\{a^+,b^+\right\} \\ b^-, & \ \text{when} \ \sigma = a^+, \\ a^-,
& \ \text{when} \ \sigma = b^+.  \end{cases}
$$  
It is straightforward to check that $\varphi \circ \phi = \phi \circ
\varphi = \id_{\Sigma}$, which is all it takes to conclude that
$\varphi$ is a homeomorphism with inverse $\phi$. To show that
$\varphi \circ S = S \circ \phi$ it suffices to observe that
\begin{equation*}
\begin{split}
& \lim_{t \uparrow 1}
\left(h^i|_{\sigma_i}\right)^{-1}\left(\psi^{\sigma_0}_{-sign(\sigma)}(t)\right)
  = \lim_{s \downarrow 0}
  \left(h^i|_{\sigma_i}\right)^{-1}\left(\psi^{\sigma_0}_{sign(\sigma)}(s)\right) \\
& = \lim_{t \uparrow
  1}\left(h^i|_{\phi(\sigma)_i}\right)^{-1}\left(\psi^{\phi(\sigma)_0}_{sign(\phi(\sigma))}(t)\right) = \lim_{s \downarrow 0} \left(h^i|_{\phi(\sigma)_i}\right)^{-1}\left(\psi^{\phi(\sigma)_0}_{-sign(\phi(\sigma))}(s)\right).
\end{split}
\end{equation*} 
\end{proof}
\end{lemma}

Let $D_{\infty}$ denote the infinite dihedral group, considered here as
the semidirect product of $\mathbb Z$ and $\mathbb Z_2$. Since
$\varphi S = S \varphi^{-1}$ it follows that $S$ and $\varphi$ induce an action $\alpha$ of $D_{\infty}$ by homeomorphisms of $\Sigma$.

\begin{lemma}\label{bgn4} $\alpha$ is a free action, i.e. $\alpha_g(\sigma) = \sigma \ \Rightarrow \ g = e$.
\begin{proof} Assume that $\sigma \in \Sigma$, $j,k \in \mathbb Z$ and
  that 
\begin{equation}\label{freeeq}
S^j\varphi^k(\sigma) = \sigma.
\end{equation} 
Choose $i_0 \in \mathbb N$ such that
$\# \left\{ I \in \mathcal I^i : \ I \subseteq e \right\} > |k|$
for all $i \geq i_0$ and all $e \in \mathbb E$. By construction of
$\varphi$, the intervals $\varphi^l(\sigma)_{i_0}, \ -|k| \leq l \leq
|k|$, are mutually distinct. Since $\left(S^j\varphi^k(\sigma)\right)_{i_0} =
\varphi^k(\sigma)_{i_0}$, it follows then from (\ref{freeeq}) that $k
=0$. Once this is established it it clear that $j = 0$ mod 2.
\end{proof}
\end{lemma}

When $(\sigma,t) \in \Sigma \times [0,1]$,
we define $\Phi_0(\sigma,t) \in \overline{\Gamma_c}$ such that 
$$
\Phi_0(\sigma,t)_j = \begin{cases} \psi^{\sigma_0}_{sign(\sigma)}(t), & \ j
  = 0, \\ \left(h^j|_{\sigma_j}\right)^{-1}\left(
    \psi^{\sigma_0}_{sign(\sigma)}(t)\right), & \ j \geq 1, \ 0 < t < 1, \\
\lim_{s \uparrow 1} \left(h^j|_{\sigma_j}\right)^{-1}\left(
    \psi^{\sigma_0}_{sign(\sigma)}(s)\right), & \ j \geq 1, \ t =1, \\
\lim_{s \downarrow 0} \left(h^j|_{\sigma_j}\right)^{-1}\left(
    \psi^{\sigma_0}_{sign(\sigma)}(s)\right), & \ j \geq 1, \ t =0. 
  \end{cases}
$$
 It is straightforward to check that $\Phi_0$ is a continuous surjection. Recall that the
 \emph{suspension} $susp(\varphi)$ of $\varphi$ is the quotient of
 $\Sigma \times [0,1]$ obtained by identifying $(\sigma,1)$ with
 $\left(\varphi(\sigma),0\right)$. Given a pair $(\sigma,t)\in \Sigma
 \times [0,1]$ we denote its image in $susp(\varphi)$ by
 $\left[\sigma,t\right]$. Since $\Phi_0\left(\sigma,1\right) =
\Phi_0\left(\varphi(\sigma),0\right)$ it follows that $\Phi_0$ falls to a continuous surjection $\Phi : susp(\varphi) \to
\overline{\Gamma_c}$ defined such that $\Phi[\sigma,t] =
\Phi_0(\sigma,t)$.

Since $S \varphi = \varphi^{-1}S$ we can define a
homeomorphism $S : susp(\varphi) \to susp(\varphi)$ such that
$$ 
S[\sigma,t] = \left[S\sigma, 1-t\right].
$$
$S$ is an involution and we can consider $S$ as an action by
$\mathbb Z_2$ on $susp(\varphi)$.

\begin{lemma}\label{bgn17} $\Phi$ induces a homeomorphism $\Phi : \ susp(\varphi)/\mathbb Z_2 \ \to \ \overline{\Gamma_c}$.
\begin{proof} It suffices to show that 
\begin{equation}\label{vee}
\Phi(\xi) = \Phi(\xi') \ \Leftrightarrow \ \xi = \xi'  \ \vee \
  S\xi = \xi'.
\end{equation} 
The implication $\Leftarrow$ is obvious. To prove the implication in
the other direction, let $\xi = (\sigma,t), \ \xi' = (\sigma',t')$. If
$t \notin \{0,1\}$ it is clear that $\Phi(\xi) = \Phi(\xi')$
implies that $(\sigma,t) = (\sigma',t')$ or $(\sigma,t) =
(S\sigma',1-t')$. Assume then that $t =0$ and that $\Phi(\xi) =
\Phi(\xi')$. Then $t' \in \{0,1\}$. If $\Phi(\xi)_j = v$ for all $j$, we have that $v \in
\sigma_j$ and $\lim_{s \downarrow 0} \left(h^i|_{\sigma_j}\right)^{-1}
\left(\psi^{\sigma_0}_{sign(\sigma)}(s)\right) = v$ for all $j$. It
follows then from Lemma \ref{q1} that $\sigma \in \left\{a^+,b^+\right\}$. If also $t' =
0$ we find in the same way that $\sigma' \in
\left\{a^+,b^+\right\}$. Thus $\xi \in \left\{\left[a^+,0\right],
  \left[b^+,0\right]\right\}$ and $\xi' \in \left\{\left[a^+,0\right],
  \left[b^+,0\right]\right\}$. Since $\left[a^+,0\right] =
\left[\varphi^{-1}\left(a^+\right),1\right] = \left[b^-,1\right] =
S\left[b^+,0\right]$ it follows that $\xi = \xi'$ or $S\xi = \xi'$ in
this case. If instead $t' = 1$ we find that $\xi' \in \left\{\left[a^-,1\right],
  \left[b^-,1\right]\right\}$. Since $\left[a^-,1\right] =
\left[\varphi\left(a^-\right),0\right] = \left[b^+,0\right]$ and $\left[b^-,1\right] =
\left[\varphi\left(b^-\right),0\right] = \left[a^+,0\right]$ we find
again that $\xi = \xi'$ or $S\xi = \xi'$. Assume then that
$\Phi(\xi)_j \neq v$ for some $j$. Then $\lim_{s \downarrow 0} \left(h^i|_{\sigma_i}\right)^{-1}
\left(\psi^{\sigma_0}_{sign(\sigma)}(s)\right) = \lim_{s \downarrow 0} \left(h^i|_{\sigma'_i}\right)^{-1}
\left(\psi^{\sigma'_0}_{sign(\sigma')}(s)\right)$ for all $i \geq j$ (if $t' = 0$) or $\lim_{s \downarrow 0} \left(h^i|_{\sigma_i}\right)^{-1}
\left(\psi^{\sigma_0}_{sign(\sigma)}(s)\right) = \lim_{s \uparrow 1} \left(h^i|_{\sigma'_i}\right)^{-1}
\left(\psi^{\sigma'_0}_{sign(\sigma')}(s)\right)$ for all $i \geq j$ (if $t'=1$). In the
first case, where $t' = 0$, we find that $\varphi S(\sigma) = \sigma'$
and hence that $S\xi =
S\left[\sigma,0\right] = S \left[\varphi^{-1}(\sigma),1\right] =
\left[S \varphi^{-1}(\sigma),0\right] = \left[\sigma',0\right] = \xi'$
and in the second case, where $t' = 1$, we find that $\sigma =
\varphi\left(\sigma'\right)$ and hence $\xi = \left[\sigma,0\right] =
  \left[\varphi(\sigma'),0\right] = \left[\sigma',1\right] = \xi'$. 

All in all we have now established (\ref{vee}) when $t \in
[0,1[$. The case $t = 1$ is handled in way analogues to the case $t=0$.

\end{proof}
\end{lemma}

Set $Z = \left\{(x_i)_{i=0}^{\infty} \in \overline{\Gamma_c} : \
    x_0 = v \right\}$ and note that $Z$ is a compact subset of
  $\overline{\Gamma_c}$, as well as of
  $W_{\overline{\Gamma_c},\overline{h}^{-1}}$ with respect to the
  Wagoner topology. Define $\Lambda : \Sigma \to Z$ such that
$$
\Lambda(\sigma) = \Phi[\sigma,0] .
$$
It follows from (\ref{vee}) that $\Lambda$ is a continuous surjection such
  that 
\begin{equation}\label{bgn5}
\Lambda(\sigma) = \Lambda(\sigma') \ \Leftrightarrow \ \sigma'
  \in \left\{\sigma, \ \varphi
  \circ S(\sigma)\right\}. 
\end{equation}
In particular, $\Lambda$ is a local homeomorphism
  which is 2 to 1.
%\begin{proof} Let $x = \left(x_i\right)_{i=0}^{\infty} \in Z$. If $x_i
%  = v$ for all $i$, $x = \Lambda\left[a^+,0\right]$. Otherwise
%  consider the least $j$ for which $x_j \neq v$. Choose $\sigma_j \in
%  \mathcal I^j$ such that $x_j \in \sigma_j$, set $\sigma_0 = h^j\left(\sigma_j\right)$ and let%  $sign\left(\sigma\right) \in \{\pm \}$ be determined by the
%  condition that $x_j = \lim_{t \downarrow 0}
%  \left(h^j|_{\sigma_j}\right)^{-1}\left(\psi^{\sigma_0}_{sign(\sigma)}(t)\right)$. Set $\sigma_i = h^{j-i}\left(\sigma_j\right), \ i \leq j$, and let $\sigma_i, i \geq j$, be determined by the condition that $x_i \in \sigma_i$ and $h^{i-j}\left(\sigma_i\right) = \sigma_j$. Then $\sigma = \left(sign(\sigma), \sigma_0,\sigma_1, \dots \right) \in \Sigma$ and $\Lambda(\sigma) =x$. This prove the surjectivity of $\Lambda$. The remaining statements about $\Lambda$ follows from Lemma \ref{bgn3} and Lemma \ref{bgn4} .
%\end{proof}
%\end{lemma}

\begin{lemma}\label{locconj?} Let $j,k \in \mathbb Z$ and let $\sigma \in
  \Sigma$. There is an open neighborhood $U$ of $\Lambda(\sigma)$ and a
  local conjugacy $\left(U,V,\chi\right)$ in $Z$ from $\Lambda(\sigma)$ to
    $\Lambda\left(S^j\varphi^k(\sigma)\right)$ such that $\chi(z) =
    \Lambda\left(S^j\varphi^k\left(\Lambda^{-1}(z)\right)\right)$ for
    all $z \in U$.
\begin{proof} Note that for every $i \in \mathbb N$ the elements
  $\Lambda\left(\varphi(\mu)\right)_i$ and $\Lambda(\mu)_i$ of
  $\Gamma_c$ are in
neighboring intervals from $\mathcal I^i$ while $\Lambda\left(S(\mu)\right)_i$
and $\Lambda(\mu)_i$ are in the same of these intervals. It follows
therefore from (\ref{diam}) that $\lim_{i \to \infty} \sup_{\mu \in \Sigma} d\left(
    \Lambda\left(\varphi(\mu)\right)_i , \Lambda(\mu)_i\right) = 0$
  and $\lim_{i \to \infty} \sup_{\mu \in \Sigma} d\left(
    \Lambda\left(S(\mu)\right)_i, \Lambda(\mu)_i\right) = 0$. Hence 
$$
\lim_{i \to \infty} \sup_{\mu \in \Sigma} d\left(
    \Lambda\left(S^j\varphi^k(\mu)\right)_i, \Lambda(\mu)_i\right) = 0 .
$$
Since $\Lambda$ is a local homeomorphism there are open neighborhoods
$U_0$ and $V_0$ of $\sigma$ and $S^j\varphi^k(\sigma)$ in $\Sigma$
such that $\Lambda(U_0)$ and $\Lambda(V_0)$ are open in $Z$ and
$\Lambda : U_0 \to \Lambda(U_0)$, $\Lambda : V_0 \to \Lambda(V_0)$ are
homeomorphisms. Set $U = \Lambda\left(U_0 \cap \left(S^j\varphi^k\right)^{-1}\left(V_0\right)\right)$. 
\end{proof}
\end{lemma}

\begin{lemma}\label{locconj??} There is a $\delta > 0$ with the
  following property: When $x,y \in \overline{\Gamma_c}$ are such that $d\left(x_k,y_k\right) < \delta$ for all $k \in
  \mathbb N$, it follows that either
\begin{enumerate}
\item[a)] there is an $N_1 \in \mathbb N$ such that $x_k$ and $y_k$
  are contained in the closure of the same interval from $\mathcal I^1$ for each $k \geq N_1$,
  or
\item[b)] there is an $N_2 \in \mathbb N$ such that $x_k,y_k \in
  J$ for all $k \geq
  N_2$.
\end{enumerate}
\begin{proof} In the proof we assume, as we may, that $d$ is an
  arclength metric. It follows from (\ref{secondcond2}) that there is a $\delta_0 > 0$ such that 
$$
h\left(\left\{ x
      \in J^+_1 : \ d(x,v) \leq \delta_0 \right\} \right) \subseteq
  J^+_1
$$ and 
$$
h\left(\left\{ x
      \in J^-_1 : \ d(x,v) \leq \delta_0 \right\} \right) \subseteq
  J^-_1. 
$$
Set 
$$
J_0 = \left\{ x
      \in J^+_1 : \ d(x,v) \leq \delta_0 \right\} \cup \left\{ x
      \in J^-_1 : \ d(x,v) \leq \delta_0 \right\} .
$$
Then, if $a \in \left\{ x
      \in J^+_1 : \ d(x,v) \leq \delta_0 \right\}$, $b \in \left\{ x
      \in J^-_1 : \ d(x,v) \leq \delta_0 \right\}$ and
    $\left\{h^k(a),h^k(b)\right\} \nsubseteq J$ for some $k \in
    \mathbb N$, there must be an $l \in \{1,2, \dots, k\}$ such that
    $d\left(h^l(a),h^l(b)\right) \geq \delta_0$.

There is a $\delta_1 > 0$ such that $d(z,v) < \delta_1 \ \Rightarrow \
h(z) \in J_0$ and there is a $\delta_2 > 0$ such that $d(x,y) <
\delta_2 \ \Rightarrow \ d\left(h(x),h(y)\right) <
\frac{\delta_1}{2}$. Let $\delta_3 = \min \left\{d(x,y) : \ x,y \in
h^{-1}(v), \ x \neq y\right\}$. Set
$\delta = \min \left\{\delta_0, \delta_1, \delta_2, \delta_3\right\}$ and
consider elements $x , y \in \overline{\Gamma_c}$ such that
$d\left(x_j,y_j\right) < \delta$ for all $j$. Assume that there is no
$N_1$ such that a) holds and also no $N_2$ such that b) holds. There
is then an $M_0$ such that $\left\{x_{M_0},y_{M_0}\right\} \nsubseteq
J$, and there is an $M_1 > M_0+2$ such that $x_{M_1}$ and $y_{M_1}$
are not both in the closure of the same element from $\mathcal
I^1$. There is then a path in $\Gamma_c$ of length $\leq \delta$ which
contains an element $w$ from $h^{-1}(v)$ and connects $x_{M_1}$ to
$y_{M_1}$. In particular, it follows that $d\left(x_{M_1-1},v\right) < \delta_1$,
$d\left(y_{M_1-1},v\right) < \delta_1$ and hence that $x_{M_1-2}
,y_{M_1-2} \in J_0$. If $x_{M_1 -2},y_{M_1-2}$ are both in $J^+_1$ or
both in $J^-_1$, let $l \geq M_1$ be the least number such that
$d\left(x_l,v\right) > \delta_1$ or $d\left(y_l,v\right) > \delta_1$
which exists because b) fails. There is then an
  element $z$ from $h^{M_1-l}(v) \backslash \{v\}$ between $x_l$ and $y_l$
  such that $h^{l-M_1+2}$ fails to be injective on any open
  neighborhood of $z$, contradicting the nonfolding condition. Hence
  $x_{M_1-2}$ and $y_{M_1-2}$ can not both be in $\left\{ x
      \in J^+_1 : \ d(x,v) \leq \delta_0 \right\}$ or both in $\left\{ x
      \in J^-_1 : \ d(x,v) \leq \delta_0 \right\}$. Since
    $\left\{x_{M_0},y_{M_0}\right\} \nsubseteq J$, it follows that
    $d\left(x_i,y_i\right) \geq \delta_0$ for some $i$ between $M_0$
    and $M_1-2$. This violates the
    condition that $d\left(x_j,y_j\right) < \delta \leq \delta_0$ for all $j$.
\end{proof} 
\end{lemma}

\begin{lemma}\label{locconj???}  Let $\sigma, \sigma' \in \Sigma$. The
  following are equivalent:
\begin{enumerate}
\item[i)] $\Lambda(\sigma)$ and $\Lambda(\sigma')$ are locally conjugate in
  $Z$.
\item[ii)] $\lim_{k \to \infty}
  d\left(\Lambda(\sigma)_k,\Lambda(\sigma')_k\right) = 0$.
\item[iii)] $\sigma$ and $\sigma'$ are in the same $D_{\infty}$-orbit
  under the action $\alpha$.
\end{enumerate}
\begin{proof} iii) $\Rightarrow$ i) follows from Lemma \ref{locconj?}
  and i) $\Rightarrow$ ii) is obvious. To prove that ii) $\Rightarrow$
  iii) assume that $\lim_{i \to
    \infty} d\left(\Lambda(\sigma)_i, \Lambda(\sigma')_i\right) =
  0$. Lemma \ref{locconj??} gives us then two alternatives, valid for
  any $i \geq N$, where $N$ is either the $N_1$ from a) of Lemma
  \ref{locconj??} or the $N_2$ from b) of the same lemma: There is an
  interval, which is either the closure of an interval
  from $\mathcal I^1$ or the interval $J =\psi_-^{e_1}\left(\left[0,\frac{1}{4}\right]\right) \cup
  \psi^{e_j}_+\left(\left[0,\frac{1}{4}\right]\right)$, such that the
  position of $\Lambda(\sigma)_i$ and $\Lambda(\sigma')_i$ in that
  interval is
\begin{equation*}
\begin{xymatrix}{
 &  & \Lambda(\sigma)_i \ar@{=}[d] & 
    &   &   &  &  &  & \Lambda(\sigma')_i \ar@{=}[d]   \\
 & \ar@{-}[r]^{l_1} & c_1 \ar@{-}[r]^{l_2} & c_2
    \ar@{-}[r]^{l_3} & c_3  \ar@{-}[r]^{l_4} &  \hdots & \hdots & \ar@{-}[r]^{l_{k-1}} & c_{k-1} \ar@{-}[r]^{l_k} & c_k  
}
\end{xymatrix}
\end{equation*}
Here the linepieces, which have been labeled $l_1$ to $l_k$, are intervals
from $\mathcal I^i$ and the points labeled $c_1$ to $c_k$, with
$\Lambda(\sigma)_i$ equal to the first and $\Lambda(\sigma')_i$ equal
to the last, are the points from $h^{-i}(v)$ which divide them. 

Note that $l_1 = \sigma_i$ or $l_2 = \sigma_i$. Assume first that
$\sigma_i = l_1$. Then 
$$
\Lambda(\sigma)_i = \lim_{s \downarrow 0}
\left(h^i|_{\sigma_i}\right)^{-1}\left(\psi^{\sigma_0}_{sign(\sigma)}(s)\right)
= \lim_{t \uparrow 1}
\left(h^i|_{\mu_i}\right)^{-1}\left(\psi^{\mu_0}_{sign(\mu)}(t)\right)
$$
where $\mu = \varphi^{-1}(\sigma)$. Hence
$$
\Lambda\left(\varphi^{-1}(\sigma)\right)_i = \lim_{s \downarrow 0}
\left(h^i|_{\mu_i}\right)^{-1}\left(\psi^{\mu_0}_{sign(\mu)}(s)\right)
= c_2.
$$ 
Similarly,
$\Lambda\left(\varphi^{-2}(\sigma)\right)_i =
\Lambda\left(\varphi^{-1}\left(\varphi^{-1}(\sigma)\right)\right)_i =
c_3$ etc. In particular, $\Lambda(\sigma')_i =
\Lambda\left(\varphi^{-k+1}(\sigma)\right)_i$. Since $i \geq N$, $h$
is injective (and expansive) from an interval containing
$\Lambda(\sigma)_{i+1}$ and $\Lambda(\sigma')_{i+1}$ onto an interval
containing $l_1,l_2, \dots, l_k$, regardless of which of the two cases
from Lemma \ref{locconj??} we consider. It follows therefore that
$\Lambda\left(\varphi^{-k+1}(\sigma)\right)_{i+1} =
\Lambda\left(\sigma'\right)_{i+1}$, and by a repetition of
the argument $\Lambda\left(\varphi^{-k+1}(\sigma)\right)_{i'} =
\Lambda\left(\sigma'\right)_{i'}$ for all $i' \geq i$. Since
$\Lambda\left(\varphi^{-k+1}(\sigma)\right)_{i'} =
h^{i-i'}\left(\Lambda\left(\varphi^{-k+1}(\sigma)\right)_i\right) = 
h^{i-i'}\left(\Lambda\left(\sigma'\right)_i\right) =
\Lambda(\sigma')_{i'}$ when $i' \leq i$ we conclude that
$\Lambda\left(\varphi^{-k+1}(\sigma)\right) = \Lambda(\sigma')$. It follows
then from (\ref{bgn5}) that $\sigma$ and $\sigma'$ are in the
same $D_{\infty}$-orbit in this case. 

Assume next that $\sigma_i = l_2$. Set $\mu = S 
\varphi^{-1}(\sigma)$ and note that $\Lambda(\mu)_i =
l_1$. Furthermore, with $\mu' = \varphi^{-1}(\sigma)$ we find that 
\begin{equation*}
\begin{split}
&\Lambda\left(S  \varphi^{-1}(\sigma)\right)_i = \lim_{s
  \downarrow 0}
\left(h^i|_{\mu_i}\right)^{-1}\left(\psi^{\mu_0}_{sign(\mu)}(s)\right) =\lim_{t
  \uparrow 1}
\left(h^i|_{\mu'_i}\right)^{-1}\left(\psi^{\mu'_0}_{sign(\mu')}(t)\right)
\\
&= \lim_{s
  \downarrow 0}
\left(h^i|_{\sigma_i}\right)^{-1}\left(\psi^{\sigma_0}_{sign(\sigma)}(s)\right)
= \Lambda(\sigma)_i .
\end{split}
\end{equation*}
We can then argue as above, with $\sigma$ replaced by $S \varphi^{-1}(\sigma)$, to deduce that $\Lambda\left( S\varphi^{k-2}(\sigma)\right) = \Lambda \left(\varphi^{-k+1} S
  \varphi^{-1}(\sigma)\right) = \Lambda(\sigma')$. As above (\ref{bgn5}) tells us then that $\sigma$ and $\sigma'$ are in the
same $D_{\infty}$-orbit.
\end{proof}
\end{lemma}

\begin{lemma}\label{locconj!} Let $(U,V,\chi)$ be a local conjugacy
  from $x$ to $y$ in $Z$. Choose $\sigma, \sigma' \in \Sigma$ such
  that $\Lambda(\sigma) =x$ and $\Lambda(\sigma') =y$. It follows that
  there is an open and closed neighborhood $T$ of $\sigma$ and an
  element $g \in D_{\infty}$ such that 
\begin{enumerate}
\item[i)] $\Lambda(T)$ is open and closed in $Z$,
\item[ii)] $\Lambda(T) \subseteq U$, 
\item[iii)] $\Lambda : T \to \Lambda(T)$ is a homeomorphism
\item[iv)] $\alpha_g(\sigma) = \sigma'$, and
\item[v)] $\chi(z) = \Lambda\left(\alpha_g\left(
      \left(\Lambda|_{T}\right)^{-1}(z)\right)\right), \ z \in \Lambda(T)$.
\end{enumerate}
\begin{proof} Combine i) $\Rightarrow$ iii) of Lemma \ref{locconj???}
  with Lemma \ref{locconj?} and Lemma 1.4 of \cite{Th1}.
\end{proof}
\end{lemma}

\begin{thm}\label{1mainrev} The heteroclinic algebra
  $B_{\overline{h}^{-1}}\left(\overline{\Gamma_c}\right)$ of
  $\left(\overline{\Gamma_c}, \overline{h}^{-1}\right)$ is
  $*$-isomorphic to $\mathbb K \otimes \left( C(\Sigma)
    \rtimes_{\alpha} D_{\infty}\right)$.
\begin{proof} Since $\alpha$ is a free action the crossed product $C(\Sigma)
    \rtimes_{\alpha} D_{\infty}$ is $*$-isomorphic to the reduced
    groupoid $C^*$-algebra $C^*_r(R)$ where $R$ is the \'etale
    equivalence relation given by orbit equivalence under $\alpha$,
    with the topology obtained from the product topology of $\Sigma
    \times D_{\infty}$ via the bijection $(x,g) \to \left(x,
      \alpha_g(x)\right) \in R$. This follows from \cite{Re} and/or \cite{Ph1}. Since $\Lambda : \Sigma \to Z$ is a
    local homeomorphism and $\Sigma$ is totally disconnected there is
    an open and closed subset $U \subseteq \Sigma$ such that
    $\Lambda(U) \subseteq Z$ is open and $\Lambda : U \to \Lambda(U)$
    is a homeomorphism. It follows from Lemma \ref{locconj?} and
    Lemma \ref{locconj!} that $\Lambda : U \to \Lambda(U)$
    induces an isomorphism of \'etale equivalence relations $R|_U \to
    R_{\overline{h}^{-1}}\left(\overline{\Gamma_c},\Lambda(U)\right)$,
    where the latter is the \'etale equivalence relation, for the case
    at hand, from Section
    1.2 of
    \cite{Th1}. Hence $C^*_r\left(R|_U\right)$ is $*$-isomorphic to
    $A_{\overline{h}^{-1}}\left(\overline{\Gamma_c}, \Lambda(U)\right)$, in the
    notation of \cite{Th1}. Since
    $B_{\overline{h}^{-1}}\left(\overline{\Gamma_c}\right) =
    A_{\overline{h}^{-1}}\left(\overline{\Gamma_c},
      W_{\overline{\Gamma_c}, \overline{h}^{-1}}\right)$ is simple by Lemma
      \ref{simstab} it follows from Proposition 4.7 of \cite{Re} that every element of
      $W_{\overline{\Gamma_c}, \overline{h}^{-1}}$, and in particular
      every element of $Z$, is locally conjugate to an element of
      $\Lambda(U)$. It follows then from Lemma \ref{locconj???} that
      the $D_{\infty}$-orbit of every element $\sigma$ of $\Sigma$
      contains an element from $U$. This implies that $U$ is an
    abstract transversal in $R$, in the sense of \cite{MRW}, so that
    $C^*_r\left(R\right)$ is stably isomorphic to
    $C^*_r\left(R|_U\right)$ by the main result of \cite{MRW}. It follows that
    $B_{\overline{h}^{-1}}\left(\overline{\Gamma_c}\right)$ is stably
    isomorphic to $C(\Sigma)
    \rtimes_{\alpha} D_{\infty}$, which gives the theorem since
    $B_{\overline{h}^{-1}}\left(\overline{\Gamma_c}\right)$ is simple by Lemma
    \ref{simstab}.
\end{proof}
\end{thm}

\begin{cor}\label{minimal}
The action $\alpha$ of $D_{\infty}$ on $\Sigma$ is minimal.
\begin{proof} This was essentially established in the proof of Theorem
  \ref{1mainrev}, but it follows also from the statement of the
  theorem by using the wellknown characterisation of ideals in crossed
  product $C^*$-algebras of free actions, combined with the simplicity of
  $B_{\overline{h}^{-1}}\left(\overline{\Gamma_c}\right)$, cf. Lemma \ref{simstab}.  
\end{proof}
\end{cor}

\begin{lemma}\label{minimalityofphi} $\varphi$ is minimal if and only
  if $\left(\Gamma_c,h\right)$ is not oriented.
\begin{proof}

Let $\psi_e : [0,1] \to e, \ e \in \mathbb E$, be
  homeomorphisms determining an orientation of $\Gamma_c$. For each
  pair $e,e' \in \mathbb E$, $e \cap h^{-1}\left(e'\right)$ is the
  union of closed intervals, $I^e_1, I^e_2, \dots, I^e_{n(e,e')}$, determined by the condition that $h :
  I^e_j \to e'$ is a homeomorphism for each $j$. We will refer to these
  functions as the characteristic functions of $h$ since they are basically
  the same as the characteristic functions of Section \ref{het}. As
  observed in the proof of n5) of Section \ref{het} there is a $k \in
  \mathbb N$ and a choice of orientation of $\Gamma_c$ such that either all the characteristic functions of $h^k$ are
  increasing or else there are both decreasing and increasing
  functions among the functions $h^k : I^e_j \to e', \ j = 1,2,
  \dots, n(e,e')$, for every $e,e' \in \mathbb E$. If the latter occurs
  $h$ can clearly not be oriented. Conversely, if $h$ is not oriented,
  the first possibility is ruled out by Theorem \ref{dichotomy} and
  Lemma \ref{invlemma1}. On the other hand the last condition will
  hold for all orientations of $\Gamma_c$ if and only if it holds for
  one. It follows that what we must show is that $\varphi$ is minimal
  if and only if the last condition holds. In fact, thanks to the
  mixing condition it suffices to find a $k \in \mathbb N$ and a
  single pair of edges $e,e'$ such that there are both decreasing and increasing
  functions among the functions $h^k : I^e_j \to e', \ j = 1,2,
  \dots, n(e,e')$.

Assume first that $\varphi$ is minimal. Consider the element $a^+ =
\left(-, \inter e_1, J^-_1,J^-_2, \dots\right) \in \Sigma$. Since $\varphi$ is minimal
there is a $M \in \mathbb N$ such that
$sign\left(\varphi^M(a^+)\right) = +$ and $\varphi^M(a^+)_0 =
a^+_0 = \inter e_1$. Now take $k_0 \in \mathbb N$ so large that $e_1$
contains more than $M$ intervals from $\mathcal I^{k_0}$. It follows then
from the definition of $\varphi$ that the set of functions $h^{k_0} :
I^{e_1}_j \to e_1, \ j =1,2, \dots, n\left(e_1,e_1\right)$, can't all
have the same orientation; if they did we would have that
$sign\left(\varphi^M(a^+)\right) = sign(a^+) = -$. As pointed out
above this gives us the desired conclusion. 

\smallskip

(A less direct argument goes as follows: If $\left(\Gamma_c, h\right)$
is oriented, $\left(\Gamma_c, h^2\right)$ is orientable in the sense
of \cite{Y4} and it follows then from Corollary 2.7 of \cite{Y4} that
$K_1\left(B_{\overline{h}^{-1}}\left(\overline{\Gamma_c}\right)\right)
\simeq \mathbb Z$. This contradicts Theorem \ref{1mainrev} and
Corollary \ref{trivK1}. Hence $\left(\Gamma_c, h\right)$
is not oriented.)

\smallskip

Conversely, assume that $\left(\Gamma_c, h\right)$
is not oriented. Thanks to Corollary \ref{minimal} it will be enough
to consider an element $\sigma \in \Sigma$ and show that 
$$
S\sigma \in \overline{\left\{ \varphi^i(\sigma) : \ i \in \mathbb
    Z\right\}} .
$$
To this end we must consider an $N \in \mathbb N$ and show that there is a $j
\in \mathbb Z$ such that 
\begin{equation}\label{firstmil}
sign \left(\varphi^j(\sigma)\right) = sign
\left(S \sigma\right) = - sign(\sigma)
\end{equation} 
and 
\begin{equation}\label{secondmil}
\varphi^j(\sigma)_i = (S\sigma)_i = \sigma_i, \  i = 0,1,2, \dots, N.
\end{equation}
As argued above there is a $k \in \mathbb N$ such that there are both decreasing and increasing
  functions among the functions $h^k : I^e_j \to e', \ j = 1,2,
  \dots, n(e,e')$, for every $e,e' \in \mathbb E$. Let $e, e' \in
  \mathbb E$ be the edges such that $\sigma_N
  \subseteq e$ and $\sigma_{N+k} \subseteq e'$. By assumption there is a
  $\sigma' \subseteq e'$ such that $h^k\left(\sigma'\right) =
  \sigma_N$ and such that $h^k : \sigma' \to \sigma_N$ has the
  opposite orientation of $h^k : \sigma_{N+k} \to \sigma_N$. Let $l \in
  \mathbb N$ be the number of intervals from $\mathcal I^{N+k}$ which
  sit between $\sigma_{N+k}$ and $\sigma'$ in $e'$. Then $j = l$ or $j
  = -l$ will satisfy (\ref{firstmil}) and (\ref{secondmil}).   
\end{proof}
\end{lemma}

Theorem \ref{1mainrev} and Lemma \ref{minimalityofphi} will be our main tool for the study of the
structure of $B_{\overline{h}^{-1}}\left(\overline{\Gamma_c}\right)$,
but first we pause to show, in the following subsection, how
the above double cover and the action of $D_{\infty}$ is related to
the double cover of Fokkink and Yi, cf. \cite{F} and \cite{Y3}.

\subsection{The relation to the oriented double cover of Yi}

This section is not needed for the remaining part of the paper. It is
included to put the dihedral action constructed above in
perspective, and to prove Proposition \ref{doublecoveraction} which
may have some independent interest.

Starting from an arbitary pre-solenoid (not just one which is a wedge of
circles) and building on work of R. Fokkink, \cite{F}, I. Yi gave in
\cite{Y3} an explicit description of an oriented graph $\Gamma'$ and
an orientation preserving continuous map $h' : \Gamma' \to \Gamma'$
which is a double cover of $h: \Gamma \to
\Gamma$ in the sense that there is 2 to 1 continuous surjection $p :
\Gamma' \to \Gamma$ such that
\begin{equation}\label{commute}
\begin{xymatrix}{
\Gamma' \ar[r]^-{h'} \ar[d]_-p &  \Gamma' \ar[d]^-p \\
\Gamma \ar[r]_-h  & \Gamma }
\end{xymatrix}
\end{equation}
commutes. See \cite{Y3}. $p$ has the following additional properties:
\begin{enumerate}
\item[1)] $\mathbb V_{\Gamma'} = p^{-1}\left(\mathbb
    V_{\Gamma}\right)$ where $\mathbb V_{\Gamma}$ and $\mathbb
  V_{\Gamma'}$ are the vertex sets of $\Gamma$ and $\Gamma'$, respectively.
\item[2)] For every edge $e \subseteq \Gamma$ there are two edges,
  $e^+$ and $e^-$, such that $p^{-1}(e) = e^+ \cup e^-$ and
\begin{enumerate}
\item[{}] $p : e^+ \to e$ is an orientation preserving
  homeomorphism, 
\item[{}] $p : e^- \to e$ is an orientation reversing homeomorphism.
\end{enumerate}
\end{enumerate}

We need a few observations regarding Yis construction.

%\begin{lemma}\label{sep} Let $d$ be a metric for the topology of
%  $\Gamma'$. There is then a $\delta > 0$ such that $d(a,b) \geq
%  \delta$ when $a \neq b$ and $p(a) = p(b)$.
%\begin{proof} Since all metrics are equivalent we can assume that $d$
%  is the arc-length metric corresponding to a realisation of $\Gamma'$
%  where each edge has length $1$. Then $\delta
%  = \frac{1}{2}$ works because $p$ is 2 to 1..
%\end{proof}
%\end{lemma}

\begin{lemma}\label{leftresol} Let $y \in \Gamma, z_0 \in \Gamma'$
  such that $p\left(z_0\right) = h(y)$. Then
$$
\# \{ z \in p^{-1}(y) : h'(z) = z_0   \} = 1 .
$$
\begin{proof} Assume first that $x = p(z_0) \notin \mathbb
  V_{\Gamma}$. There are then edges $e,f \subseteq \Gamma$ and an
  open non-degenerate interval $I \subseteq e$ such that $y \in I$
  and $h(I) = \inter f \ni x$. It follows from 1) and 2) that
  $p^{-1}(y) = \left\{z_+,z_-\right\}$ where $z_{\pm} \in \inter
  e_{\pm}$. Similarly, $z_0 \in \inter f_+$ or $z_0 \in \inter
  f_-$. Assume that $z_0 \in \inter f_+$. If $h : I \to f$ is
  orientation preserving the commutative of (\ref{commute}) implies
  that $\{ z \in p^{-1}(y) : h'(z) = z_0   \} = \{z_+\}$ since $h'$ is
  orientation preserving. If instead $h : I \to f$ is
  orientation reversing we conclude in the same way that $\{ z \in
  p^{-1}(y) : h'(z) = z_0   \} = \{z_-\}$. Similarly, when $z_0 \in
  \inter f_-$ we find that $\{ z \in p^{-1}(y) : h'(z) = z_0   \} =
  \{z_+\}$ when $h : I \to f$ is orientation reversing and $\{ z \in
  p^{-1}(y) : h'(z) = z_0   \} = \{z_-\}$ when it is orientation
  preserving. 

 To handle the case where $x \in \mathbb V_{\Gamma}$, let $p^{-1}(y) = \left\{z_1,z_2\right\}$. Assume to reach a contradiction that $h'(z_1) =
 h'(z_2) = z_0$. Use then that $p$ is
 2 to 1 and hence a local homeomorphism to find open neighborhoods $U$
 and $V$ of $y$ and $x$, respectively, together with open sets $Z_i
 \subseteq \Gamma', W_i \subseteq \Gamma', i = 1,2$, such that $Z_1
 \cap Z_2 = \emptyset, \ z_0 \in W_1, \ z_i \in Z_i, i = 1,2, \ W_1 \cap W_2 = \emptyset$, $p\left(Z_1\right)
 = p\left(Z_2\right) = U, \ h(U) \subseteq V$ and $p\left(W_1\right) =
 p\left(W_2\right) = V$, and such that $p$ is injective on
 $W_1$. Since $h'(z_1) = h'(z_2) = z_0$ we can shrink $U$ (and hence $Z_i, i = 1,2$)
to arrange that $h'\left(Z_i\right) \subseteq W_1, i = 1,2$. For every element $\xi$ of
$U$ there is an element $a_i \in Z_i$ such that $p(a_i) =
\xi, i =1,2$. Note that $h'(a_i) \in W_1$ and that $ph'(a_i) =
hp(a_i) = h(\xi), i =1,2$. Since $p$ is injective on $W_1$ it follows
that $h'(a_1) = h'(a_2)$. This contradicts the first part of the proof
provided $h(\xi) \notin \mathbb V_{\Gamma}$. Since $h$ does not
collapse $U$ to a finite set there \emph{is} a $\xi \in U$ such that
$h(\xi) \notin \mathbb V_{\Gamma}$. This shows that $\# \{ z \in
p^{-1}(y) : h'(z) = z_0   \} \leq 1$. Finally, the conclusion that $\# \{ z \in
p^{-1}(y) : h'(z) = z_0   \} > 0$ follows by approximating $y$ with
elements $y'$ for which $h(y') \notin \mathbb V_{\Gamma}$.     
\end{proof}
\end{lemma}

Define $\overline{p} : \overline{\Gamma'} \to \overline{\Gamma}$ such
that $\overline{p}\left(\left(x_i\right)_{i=0}^{\infty}\right) =
\left( p(x_i)\right)_{i = 0}^{\infty}$. It follows from
(\ref{commute}) that $\overline{p} \circ \overline{h} = \overline{h'}
\circ \overline{p}$.

\begin{lemma}\label{yi?} Let $x,y \in \overline{\Gamma'}$ such that
  $\overline{p}(x) = \overline{p}(y)$. Then $ x =y$ if and only if
  $x_0 = y_0$. 
\begin{proof} This follows immediately from Lemma \ref{leftresol}.
\end{proof}
\end{lemma}

\begin{lemma}\label{contlemma} Let $\{x^n\}_{n=1}^{\infty}$ be a
  sequence in $\overline{\Gamma'}$ and $x \in \overline{\Gamma'}$ an
  element in $\overline{\Gamma'}$. Then $\lim_{n \to \infty} x^n = x$
  if and only if $\lim_{n \to \infty} \overline{p}\left(x^n\right) =
  \overline{p}(x)$ and $\lim_{n \to \infty} x^n_0 = x_0$.
\begin{proof} One implication follows from the continuity of
  $\overline{p}$. To prove the other assume that $\lim_{n \to \infty} \overline{p}\left(x^n\right) =
  \overline{p}(x)$ and $\lim_{n \to \infty} x^n_0 = x_0$. If $x^n$
  does not converge to $x$ there is a subsequence $\left\{
    x^{n_i}\right\}$ and an element $y \neq x$ such that $\lim_{i \to
    \infty} x^{n_i} = y$. But then the assumptions imply that
  $\overline{p}(y) = \overline{p}(x)$ and $x_0 = y_0$. This
  contradicts Lemma \ref{yi?}.
\end{proof}
\end{lemma}

To relate Yi's double cover to the constructions of the preceding
section we have to specialize, as we did above, to the case where
$\Gamma_c$ is a wedge of circles. Let $(\sigma,t) \in \Sigma \times [0,1]$. It follows from Lemma \ref{leftresol} and
Lemma \ref{yi?} that there is unique element $\Psi_0(\sigma,t) \in
\overline{\Gamma_c'}$ such that
$\overline{p}\left(\Psi_0(\sigma,t)\right) =
\Phi_0(\sigma,t)$ and 
$$
\Psi_0(\sigma,t)_0 \in \sigma_0^{sign(\sigma)}.
$$ 
Since
  $\Phi_0\left(\varphi(\sigma),0\right) =
  \Phi_0\left(\sigma,1\right)$, it
  follows from Lemma \ref{yi?} that $\Psi_0(\sigma,1) =
  \Psi_0(\varphi(\sigma),0)$. We get therefore a well-defined map
$$
\Psi : susp(\varphi) \to \overline{\Gamma_c'}
$$
defined such that $\Psi[\sigma,t] = \Psi_0(\sigma,t)$. Then the diagram
\begin{equation}\label{uio}
\begin{xymatrix}{
susp(\varphi) \ar[r]^-{\Psi} \ar[dr]_-{\Phi}& \overline{\Gamma_c'}
\ar[d]^-{\overline{p}} \\
& \overline{\Gamma_c}
}
\end{xymatrix}
\end{equation}
 commutes by construction.

\begin{lemma}\label{psi0homeo} $\Psi : susp(\varphi) \to
  \overline{\Gamma_c'}$ is a homeomorphism.
\begin{proof} It follows from Lemma \ref{contlemma} that $\Psi_0$
  is continuous and the same is therefore true for $\Psi$. Let $x
  \in \overline{\Gamma_c'}$. Since $\Phi$ is surjective there is an element
  $(\sigma,t) \in \Sigma \times [0,1]$ such that $\overline{p}(x) =
  \Phi(\sigma,t)$. Then $x_0 \in \sigma_0^{sign(\sigma)}$ or $x_0 \in
  \sigma_0^{-sign(\sigma)}$. In the first case, $x =
  \Psi_0(\sigma,t)$ and in the second $\Psi_0\left(S\sigma, 1-t\right)
    = x$. This shows that $\Psi$ is surjective. Thanks to this it
    follows from (\ref{uio}) and the fact that both $\overline{p}$ and
    $\Phi$ are 2 to 1 that $\Psi$ is also injective.   
\end{proof}
\end{lemma}

To see how the $D_{\infty}$-action fits into this picture, note that we can also consider $susp(\varphi)$ as the
  quotient of $\Sigma \times \mathbb R$ under the equivalence relation
  $\sim$ where $(\sigma,t) \sim (\sigma',t')$ means that $t-t' \in \mathbb Z$
  and $\varphi^{t-t'}(\sigma) = \sigma'$. We obtain then a flow
  $\varphi_t, t\in \mathbb R$, on $susp(\varphi)$ defined such that
$$
\varphi_s\left[\sigma,t\right] = \left[\sigma, t+s\right] .
$$
Note that $S \circ \varphi_t = \varphi_{-t} \circ S$ so that we have
in fact an action $\beta$ of the generalized dihedral group
$D_{\mathbb R} = \mathbb R \rtimes \mathbb Z_2$ on $susp(\varphi)$. Formally $D_{\mathbb R}$ is the
semidirect product $\mathbb R \rtimes \mathbb Z_2$ in the category of
locally compact groups coming from the automorphism $t \mapsto -t$ of
$\mathbb R$. If we write the elements of $D_{\mathbb R}$ as tuples
$(i,s)$ where $s\in \mathbb R$ and $i \in \mathbb Z_2 =\{0,1\}$,
$$
\beta_{(i,s)}[\sigma,t] = S^i\varphi_s[\sigma,t] = \left[\left(\varphi
    S\right)^i\sigma,
  (-1)^i(s+t)\right] .
$$

\begin{lemma}\label{freeRaction} $\beta$ is a free action,
  i.e. $\beta_x(\xi) = \xi \ \Rightarrow \ x = 0$.
\begin{proof} If $x = (i,s)$ and $\beta_{(i,s)}[\sigma,t] = \left[(\varphi S)^i\sigma,
  (-1)^i(s+t)\right] = \left[\sigma,t\right]$, where $t \in [0,1]$, we have that
\begin{equation}\label{refeq1}
\left[ \varphi^{\left[(-1)^i(s+t)\right]+i}S^{i}\sigma, r\right] =
[\sigma,t],
\end{equation}
 where $\left[(-1)^i(s+t)\right]$ denotes the integer part
of $(-1)^i(s+t)$ and $r = (-1)^i(s+t) - \left[(-1)^i(s+t)\right]$. It
follows from (\ref{refeq1}) that $r= t$ and
$\varphi^{\left[(-1)^i(s+t)\right]+i}S^{i}\sigma = \sigma$ or $r  = 0$ and $t = 1$. In the first case we deduce from Lemma \ref{bgn4}
$$
\left[(-1)^i(s+t)\right]+i = 0, \ (-1)^i = 1, \ r = t.
$$
It follows that $x = 0$ in this case. In the second case we conclude from
Lemma \ref{bgn4} that
$$
\left[(-1)^i(s+t)\right] +i - 1 = 0, \ (-1)^i = 1
$$
which shows that $x = 0$, also in this case.
\end{proof}
\end{lemma}

\begin{lemma}\label{something} Let $x,y \in \overline{\Gamma_c'}$. Then
  $\lim_{n \to \infty} d\left(\overline{h}^{-n}\left(\Phi(x)\right),
    \overline{h}^{-n}\left(\Phi(y)\right)\right) = 0$ if and only if
   $\beta_{a}(x) = y$ for some $a \in D_{\mathbb R}$.
\begin{proof} Assume first that $\beta_{a}(x) = y$ for some $a \in
  D_{\mathbb R}$. Write $a = (i,t)$ where $t \in \mathbb R$ and $i \in
  \mathbb Z_2$. Since $\Phi \circ S = \Phi$ we may assume that $i
  =0$ so that $\varphi_t(x) = y$. Write $x = (\sigma,s)$ where $s \in
  [0,1]$. Then $\varphi_t(x) = \left[\varphi^{k}(\sigma),r\right]$ where $k
  \in \mathbb Z$ and $r \in [0,1]$. Since
  $\Phi\left(\left[\varphi^k(\sigma),r\right]\right)_n$ and
  $\Phi\left(\left[\varphi^{k}(\sigma),0\right]\right)_n$ are in the
  same interval from $\mathcal I^n$ for each $n$ it follows from
  (\ref{diam}) that $\lim_{n \to \infty}
  d\left(\overline{h}^{-n}\left(\Phi\left(\left[\varphi^k(\sigma),r\right]\right)\right)
    ,\overline{h}^{-n}\left(\Phi\left(\left[\varphi^{k}(\sigma),0\right]\right)\right)\right) =
  0$. By Lemma \ref{locconj?} $\lim_{n \to \infty}
  d\left(\overline{h}^{-n}\left(\Phi\left(\left[\varphi^k(\sigma),0\right]\right)\right)
    ,\overline{h}^{-n}\left(\Phi\left(\left[\sigma,0\right]\right)\right)\right) =
  0$, and by a repetition of the preceding argument $\lim_{n \to \infty}
  d\left(\overline{h}^{-n}\left(\Phi\left(\left[\sigma,0\right]\right)\right)
    ,\overline{h}^{-n}\left(\Phi\left(\left[\sigma,s\right]\right)\right)\right) =
  0$, proving that $\lim_{n \to \infty}
  d\left(\overline{h}^{-n}\left(\Phi\left(\varphi_t(x)\right)\right)
    ,\overline{h}^{-n}\left(\Phi\left(x\right)\right) \right) = 0$,
  as dsired.

Assume next that $\lim_{n \to \infty} d\left(\overline{h}^{-n}\left(\Phi(x)\right),
    \overline{h}^{-n}\left(\Phi(y)\right)\right) = 0$. Write $x =
  [\sigma,t]$ and $y = \left[\sigma',t'\right]$ where $\sigma,\sigma'
  \in \Sigma, t,t' \in [0,1]$. Since $[\sigma,t] =
  \varphi_t[\sigma,0]$ and $[\sigma',t'] = \varphi_{t'}[\sigma',0]$
  and $\lim_{n \to \infty} d\left(\overline{h}^{-n}\left(\Phi(\varphi_{-t}(x))\right),
    \overline{h}^{-n}\left(\Phi(\varphi_{-t'}(y))\right)\right) = 0$
  by the first part of the proof we may assume that $t=t' = 0$, and
  then the desired conclusion follows from Lemma \ref{locconj???}.
\end{proof}
\end{lemma}

\begin{prop}\label{doublecoveraction} Let $\Gamma$ be a pre-solenoid and $\Gamma'$ the corresponding oriented double cover of
  $\Gamma$. Let $\overline{p} : \overline{\Gamma'} \to
  \overline{\Gamma}$ be the corresponding canonical surjection. There is then a free action $\beta$ of $D_{\mathbb R}$ on
  $\overline{\Gamma'}$ with the following property:

\smallskip

Two elements $x,y \in \overline{\Gamma'}$ satisfy that 
$$
\lim_{n \to \infty} d\left(\overline{h}^{-n}(\overline{p}(x)), \overline{h}^{-n}(\overline{p}(y))\right)
= 0
$$
if and only if $\beta_a(x) = y$ for some $a \in D_{\mathbb
  R}$.

\begin{proof} When $\Gamma$ is a wedge of circles the existence follows
  from the preceding. For the general case we use a result of
  R. Williams which shows that for
  an arbitrary pre-solenoid $(\Gamma, h)$ there is an $m \in
\mathbb N$ and a pre-solenoid $\left(\Gamma_1,h_1\right)$ such that
$\Gamma_1$ is a wedge of circles and $\left(\overline{\Gamma},
  \overline{h}^m\right)$ is conjugate to
$\left(\overline{\Gamma_1},\overline{h_1}\right)$, cf. Proposition 2.2
of \cite{Y4}. Since this conjugacy lifts to a homeomorphism between the
double covers $\overline{\Gamma'}$ and $\overline{\Gamma'_1}$ by
Proposition 3.7 of \cite{Y3} this gives us the general case. 
\end{proof}
\end{prop}

\begin{remark}\label{actioncor} It follows from Corollary \ref{minimal}
  that the action $\beta$ of Proposition \ref{doublecoveraction} is
  minimal. It is a consequence of Proposition \ref{trace} which we
  prove below that the action is also uniquely ergodic.
\end{remark}

\subsection{Affability and a unique trace}

In this section we obtain the conclusion that
$B_{\overline{h}^{-1}}\left(\overline{\Gamma}\right)$ has a densely
defined lower-semicontinous trace which is unique up to multiplication
by a scalar. For this purpose we return to the setting of Section
\ref{dihedral}. In particular, we consider the intervals $\mathcal J$
given by (\ref{union}) and use them to define an ordered Bratteli diagram $\mathbb B$ as
follows: The first level has $2n$ vertices, labeled $(+,i)$  and $(-,i),
i = 1,2, \dots, n$. Beyond this first level the diagram is stationary
and has $n$ vertices
at each level and there are $a_{ik}$ arrows, labeled by $I^l_{ik}, l = 1,2,
\dots, a_{ik}$, from
vertex $k$ at each level to vertex $i$ at the next level. From the
first to the second level there are $a_{ik}$ arrows, labeled by $I^l_{ik+}, l = 1,2,
\dots, a_{ik}$, from
vertex $(+,k)$ at the first level to vertex $i$ at the next level, and there are $a_{ik}$ arrows, labeled by $I^l_{ik-}, l = 1,2,
\dots, a_{ik}$, from
vertex $(-,k)$ at the first level to vertex $i$ at the next
level. Then every element $\sigma \in \Sigma$ defines an infinite path $\xi(\sigma)  = \left(I^{l_0}_{i_0j_0sign(\sigma)},
  I^{l_1}_{i_1j_1}, I^{l_2}_{i_2j_2}, \dots \right)$ in $\mathbb B$ by
the requirement that $\sigma_{k+1} \subseteq I^{l_k}_{i_kj_k}, \ k
\geq 0$.
Let $X_{\mathbb B}$ be the space of infinite paths in $\mathbb B$ equipped
with the relative topology inherited from the infinite product space $\left\{ I^l_{ik} : \ l =
  1,2, \dots , a_{ik},  \ i,k = 1,2, \dots, n\right\}^{\mathbb N}$. It
is straightforward to check that $\xi :  \Sigma \to X_{\mathbb B}$ is a homeomorphism.

Consider a vertex $w$ in $\mathbb B$, say the vertex $i$ at some level
beyond the second level
of the Bratteli diagram. Then the arrows in $\mathbb B$ terminating
at $w$ are labeled by $I^k_{il}, k = 1,2, \dots, a^k_{il}, l = 1,2,
\dots, n$. We equip these arrows with a total ordering such that
$I^1_{i1}$ is the minimal element and $I^{a_{ij}}_{ij}$ the
maximal. If $w =i$ is at the second level the arrows in $\mathbb B$ terminating
at $w$ are labeled by $I^k_{il\pm}, k = 1,2, \dots, a^k_{il}, l = 1,2,
\dots, n$. We equip these arrows with a total ordering such that
$I^1_{i1-}$ is the minimal element and $I^{a_{ij}}_{ij+}$ the
maximal. Add a
zeroth level with one vertex with $2n$ emitting arrows, one to each of
the vertices at the first level. With these additional features $\mathbb
B$ becomes a simple ordered Bratteli diagram in the sense of \cite{GPS1} and by
construction there is a unique minimal and a unique maximal path in
$X_{\mathbb B}$. Hence the Bratteli-Vershik map $\psi : X_{\mathbb B}
\to X_{\mathbb B}$ is defined and minimal. For $x \in X_{\mathbb B}$
and $\sigma \in \Sigma$ we let $orb_{\psi}(x)$ denote the orbit of $x$
under $\psi$ and $orb_{D_{\infty}}(\sigma)$ the orbit of $\sigma$
under $\alpha$.

\begin{lemma}\label{orbiteq} $\xi\left(orb_{D_{\infty}}(\sigma)\right) =
  orb_{\psi}\left(\xi(\sigma)\right)$ for all $\sigma \in \Sigma$.
\begin{proof} Assume $\sigma' \in orb_{D_{\infty}}(\sigma)$. It
  follows from the implication iii) $\Rightarrow$ ii) of Lemma
  \ref{locconj???} that $\lim_{j \to \infty}
  \dist\left(\sigma_j,\sigma'_j\right) = 0$. By Lemma \ref{locconj??}
  this implies that there is an $N \in \mathbb N$ such that either
  $\xi(\sigma)_k = \xi(\sigma')_k$ for all $k \geq N$ or else
  $\xi(\sigma)_k  \neq \xi(\sigma')_k$ and $\left\{\xi(\sigma)_k,
    \xi(\sigma')_k\right\} = \left\{I^1_{11},I^{a_{jj }}_{jj}\right\}$
  for all $k \geq N$. In both cases we conclude that $\xi(\sigma') \in
  orb_{\psi}\left(\xi(\sigma)\right)$, cf. \cite{GPS1}, in particular
  Remark (ii) on page 72 of \cite{GPS1}.

Conversely, if $\xi(\sigma') \in orb_{\psi}\left(\xi(\sigma)\right)$,
it follows from \cite{GPS1} that there is an $N \in \mathbb N$ such that either
  $\xi(\sigma)_k = \xi(\sigma')_k$ for all $k \geq N$ or else
  $\xi(\sigma)_k  \neq \xi(\sigma')_k$ and $\left\{\xi(\sigma)_k,
    \xi(\sigma')_k\right\} = \left\{I^1_{11},I^{a_{jj}}_{jj}\right\}$
  for all $k \geq N$. In both cases we conclude that $\lim_{k \to
    \infty} d\left(\Lambda(\sigma)_k, \Lambda(\sigma')_k\right) = 0$
  and then Lemma \ref{locconj???} implies that $\sigma' \in orb_{D_{\infty}}(\sigma)$.   
\end{proof}
\end{lemma}

It follows from Lemma \ref{orbiteq} that the \'etale equivalence
relation given by orbit equivalence under $\alpha$ is affable in the
sense of \cite{GPS2}. For the purposes of the present paper it is the
following consequence of Proposition \ref{orbiteq} we shall need:

\begin{prop}\label{trace} The action $\alpha$ of $D_{\infty}$ on
  $\Sigma$ has only one invariant Borel probability measure, and hence
  $C(\Sigma) \rtimes_{\alpha} D_{\infty}$ has a unique trace state.
\begin{proof} The simplex of $\psi$-invariant Borel probability
  measures on $X_{\mathbb B}$ agrees with the trace-simplex of the
  unital AF-algebras whose Bratteli-diagram is $\mathbb B$,
  cf. Theorem 1.13 and Theorem 3.7 of \cite{GPS1}. But
  $\mathbb B$ is eventually stationary with an inclusion pattern given
  by a primitive matrix and hence this AF-algebra has a unique trace
  state by \cite{Ef}. In this way the first statement becomes a well-known consequence of
  Lemma \ref{orbiteq} and the second is then a well-known consequence
  of the first, cf. e.g. Theorem 4.5 of \cite{Th2}.
\end{proof}
\end{prop}

It follows from Proposition \ref{trace} and Theorem \ref{1mainrev}
that the heteroclinic algebra $B_{\overline{h}^{-1}}\left(\overline{\Gamma}\right) \simeq B_{\overline{h}^{-1}}\left(\overline{\Gamma_c}\right)$ has an
essentially unique densely defined lower-semicontinuous trace. The
remaining information about the structure which we shall obtain in this
paper comes
from the fact that it is (stably) given as the crossed product of a free minimal action of
$D_{\infty}$ on the Cantor set. The next two sections are therefore
devoted to a study of the crossed products by such actions.

\subsection{On crossed products $C(\Omega) \rtimes D_{\infty}$}\label{Ktheorysection}

An action $\alpha$ of the infinite dihedral group $D_{\infty}$ on the
Cantor set $\Omega$ is generated by two homeomorphism $\varphi$ and
$\sigma$ of $\Omega$ such that
\begin{equation}\label{order2}
\sigma^2 = \id
\end{equation}
and
\begin{equation}\label{conjcond}
\sigma \varphi \sigma = \varphi^{-1} .
\end{equation}
Crossed product $C^*$-algebras arising from actions of $D_{\infty}$ on
$\Omega$ were studied by Bratteli, Evans and Kishimoto in \cite{BEK}
under the assumptions that $\varphi$ is minimal and that $\sigma$ or
$\varphi \sigma$ has a fixed point. The main result of \cite{BEK} is that
the crossed product is then AF. As we have seen, actions of
$D_{\infty}$ on $\Omega$ arising from 1-solenoids are always free so
the last assumption in \cite{BEK} is never met for such
actions. Minimality of $\varphi$, on the other hand, is equivalent to
orientability of the 1-solenoid, as it follows by combining Lemma
\ref{minimalityofphi} with Theorem \ref{orientthm}.

We will assume throughout, in this and the following section, that $\alpha$ is a free and minimal action; i.e. for
every $x \in \Omega$ the points $\varphi^j(x), \sigma\varphi^k(x), j,k
\in \mathbb Z$, are all different and constitute a dense subset of
$\Omega$.

%To study the structure of the crossed product we will use an
%$\alpha$-invariant Borel probability measure $\mu$ on $\Omega$ to
%represent the crossed product $C(\Omega) \rtimes_{\alpha} D_{\infty}$
%on the Hilbert space $L^2(\Omega) = L^2(\Omega, \mu)$ in the standard
%way: The functions $f \in C(\Omega)$ act as multiplication operators:
%$$
%(f\xi)(x) = f(x)\xi(x),
%$$
%and there are two unitaries $u$ and $u_{\sigma}$ given by the formulas
%$$
%$(u\xi)(x) = \xi\left(\varphi^{-1}(x)\right)
%$$
%and
%$$
%\left(u_{\sigma}\xi\right)(x) = \xi\left(\sigma(x)\right),
%$$
%respectively. These operators combined generate a copy of
%$C(\Omega)\rtimes_{\alpha} D_{\infty}$. $C(\Omega)$ and $u$ generate a copy of the
%crossed product $C(\Omega)\rtimes_{\varphi} \mathbb Z$ which
%$u_{\sigma}$ leaves globally invariant. Conjugation by $u_{\sigma}$ gives rise to an order two
%automorphism of $C(\Omega)\rtimes_{\varphi} \mathbb Z$ which we denote
%by $\sigma$. In this way we get an identification
%$$
%C(\Omega) \rtimes_{\alpha} D_{\infty} =
%\left(C(\Omega)\rtimes_{\varphi} \mathbb Z\right) \rtimes_{\sigma}
%\mathbb Z_2 .
%$$

\begin{lemma}\label{dichotomi} One of the following two alternatives
  occur:
\begin{enumerate}
\item[i)] $\varphi$ is minimal on $\Omega$ or
\item[ii)] there is a clopen $\varphi$-invariant subset $F \subseteq
  \Omega$ such that $\varphi|_F$ is minimal (on $F$), $\sigma(F) \cap F =
  \emptyset$ and $\sigma(F) \cup F = \Omega$.
\end{enumerate} 
If $F'$ is another subset of $\Omega$ with the properties specified in
ii), $F' = F$ or $F' =\sigma(F)$. 
\begin{proof} Assume that $\varphi$ is not minimal. There is then a non-trivial
  closed subset $F  \subseteq \Omega$ such that $\varphi(F) =
  F$. Since $F \cap \sigma(F)$ is closed and invariant under both
  $\varphi$ and  $\sigma$, the minimality of $\alpha$ implies that $F
  \cap \sigma(F) = \emptyset$ or $F \cap \sigma(F) =\Omega$. The
  latter possibility can not occur since $F$ is non-trivial, i.e. $F
  \cap \sigma(F) = \emptyset$. A similar reasoning, with $F \cup
  \sigma(F)$ replacing $F \cap \sigma(F)$, shows that $F \cup
  \sigma(F) = \Omega$. In particular, $F$ is clopen. If $F_0$ is any non-empty
  closed $\varphi$-invariant subset of $F$ the same arguments show that
  $F_0 \cap \sigma(F_0) = \emptyset$ and $F_0 \cup \sigma(F_0) =
  \Omega$. It follows that $F_0$ must be all of $F$, and hence
  $\varphi|_F$ is minimal. This proves the existence of $F$. The
  essential uniqueness of $F$, as described in the statement, is established along similar lines.
\end{proof}
\end{lemma}

Let $A$ and $B$ be isomorphic $C^*$-algebras, and let $\gamma : A \to
B$ be a $*$-isomorphism. The automorphism of $A \oplus B$ sending
$(a,b)$ to $\left(\gamma^{-1}(b),\gamma(a)\right)$ will be called \emph{a swop of
  $A$}. 

\begin{lemma}\label{swop} Let $\gamma : A \oplus B \to A \oplus B$ be
  a swop of $A$. Then
$$
(A \oplus B) \rtimes_{\gamma} \mathbb Z_2 \ \simeq \ M_2(A).
$$
\begin{proof} Left to the reader.
\end{proof} 
\end{lemma}

\begin{thm}\label{easycase} Assume that $\varphi$ is not minimal, and
  let $F$ be the clopen subset from ii) of Lemma \ref{dichotomi}. Then
\begin{equation}\label{isoeq}
C(\Omega) \rtimes_{\alpha} D_{\infty} \ \simeq \ M_2\left(C(F)
  \rtimes_{\varphi} \mathbb Z\right) ,
\end{equation}
and hence $C(\Omega) \rtimes_{\alpha} D_{\infty}$ is a simple unital
AT-algebra of real rank zero with $K_1$-group $\mathbb Z$.
\begin{proof} Since both $F$ and $\sigma(F)$ are $\varphi$-invariant
  we have a decomposition
$$
C(\Omega)\rtimes_{\varphi} \mathbb Z \ \simeq \ \left(C(F)\rtimes_{\varphi}
 \mathbb Z\right) \oplus \left(C(\sigma(F))\rtimes_{\varphi} \mathbb Z\right) .
$$
This isomorphism turns $\sigma$ into the swop of $C(F)\rtimes_{\varphi}
 \mathbb Z$ given by the $*$-isomorphism $C(F)\rtimes_{\varphi}
 \mathbb Z \to  C(\sigma(F))\rtimes_{\varphi} \mathbb Z$ which sends
 $f \in C(F)$ to $f \circ \sigma \in C\left(\sigma(F)\right)$ and $u$
 to $u^*$ when $u$ denotes the canonical unitary going with a crossed
 product by $\mathbb Z$. Hence Lemma \ref{swop} applies. The last statement follows from (\ref{isoeq}) and \cite{HPS}. 
\end{proof}
\end{thm}

We turn to the case where $\varphi$ is minimal.

\begin{prop}\label{tower3} Assume that
$\varphi$ is minimal. Let $R \in \mathbb N$ be
  given. There is then a partition 
$$
C(k,i), \ i=0,1,2, \dots, J_k-1, \ k = 1,2, \dots, N,
$$
of $\Omega$ into clopen sets such that
\begin{enumerate}
\item[i)] $J_k \geq R, k =1,2, \dots, N$,
\item[ii)] $\sigma\left(C(k,i)\right) = C\left(k, J_k-i-1\right)$ for
  all $i,k$, 
\item[iii)] $\varphi\left(C(k,i)\right) = C(k,i+1), i = 0,1,\dots,
  J_k-2$, for all $k$,
\item[iv)] $\varphi\left(\bigcup_{k=1}^N C(k,J_k-1)\right) =
  \bigcup_{k=1}^N C(k,0)  $.
\end{enumerate}  
\begin{proof} Since $D_{\infty}$ acts freely, any point in $\Omega$
  has a clopen neighborhood $Z$ such that 
$$
\left[\varphi^k\left(Z\right)
  \cup \varphi^k\left(\varphi \sigma(Z)\right)\right] \cap \left[\varphi^l\left(Z\right)
  \cup \varphi^l\left(\varphi \sigma(Z)\right)\right] = \emptyset
$$
when $k \neq l$ and $0 \leq k,l \leq R$. The desired partition can
then be obtained by applying Lemma 1.4 of \cite{BEK} with $Y = Z \cup \varphi\sigma(Z)$.
\end{proof}
\end{prop}

\begin{thm}\label{tracRokhlin} Assume that $\varphi$ is minimal. The action of $\mathbb Z_2$ on $C(\Omega)
  \rtimes_{\varphi} \mathbb Z$ given by $\sigma$ has the tracial
  Rokhlin property of Phillips \cite{Ph4}.
\begin{proof} By \cite{HPS} $C(\Omega)\rtimes_{\varphi} \mathbb Z$ is a
  unital simple AT-algebra of real rank zero. It follows therefore
  from (2.21) of \cite{EG} that $C(\Omega)
  \rtimes_{\varphi} \mathbb Z$ has tracial rank zero. By Theorem 2.9 of \cite{Ph5} the conclusion will
  follow if we for a given $\epsilon > 0$ exhibit functions $f_1,f_2 \in
  C(\Omega)$ such that $0 \leq f_i \leq 1, \ \left\|f_i - f_i \circ
    \varphi\right\| \leq \epsilon, i =1,2$, $f_1f_2 = 0$, $f_1 \circ \sigma
  = f_2$ and $\omega\left(1 - f_1 -f_2\right) \leq \epsilon$ for every
  tracial state $\omega$ of $C(\Omega)
  \rtimes_{\varphi} \mathbb Z$. To this end, choose first an $M \in
  \mathbb N$ such that $\frac{1}{M} \leq \epsilon$ and then an $R \geq
  M$ such that $\frac{4M}{R} \leq \epsilon$. With this choice of $R$,
  consider a partition $C(k,i)$ with the properties spelled out in
  Proposition \ref{tower3}. For each $k = 1,2, \dots, N$ we define a
  function $g_k : \Omega  \to [0,1]$ in the
  following way: Let $l_k$ be the integer part of
  $\frac{J_k-1}{2}$. Choose numbers $\alpha_i \in [0,1], i = 0,1,2,
  \dots, l_k$, such that $\alpha_0 = \alpha_{l_k} = 0$, $\left|\alpha_{j+1} -
    \alpha_j\right|\leq \epsilon$ for all $j = 0,1,2, \dots, l_k-1$,
  and $\alpha_i = 1, \ M \leq i \leq l_k-M-1$. This is possible since
  $\frac{1}{M} \leq \epsilon$. Now define ${g}_k$ to be constant
  equal to $\alpha_i$ on $C(k,i)$ and to be zero on $\Omega \backslash
  \left(\bigcup_{i=0}^{l_k} C(k,i)\right)$. %Let $x \in \Omega \backslash
  %\left(\bigcup_{i=0}^{l_k} C(k,i)\right)$. If $\varphi(x) \in
  %\bigcup_{i=0}^{l_k} C(k,i)$ it must be in $C(k,0)$ and then
  %$g_k\left(\varphi(x)\right) = g_k(x) = 0$. The same conclusion holds
  %when $\varphi(x) \notin \bigcup_{i=0}^{l_k} C(k,i)$. If $\varphi^{-1}(x) \in
  %\bigcup_{i=0}^{l_k} C(k,i)$, it must be in $C(k,l_k)$ and then
  %$g_k\left(\varphi^{-1}(x)\right) = g_k(x) = 0$ again. The same conclusion holds
  %when $\varphi^{-1}(x) \notin \bigcup_{i=0}^{l_k} C(k,i)$.  

%Let then $x \in
 % \bigcup_{i=0}^{l_k} C(k,i)$. If $\varphi(x) \notin
 % \bigcup_{i=0}^{l_k} C(k,i)$, $x \in C(k,l_k)$ and $g_k(\varphi(x)) =
  %g_k(x) = 0$. If $\varphi(x) \in \bigcup_{i=0}^{l_k} C(k,i)$ we have
  %instead that $\left|g_k\left(\varphi(x)\right) - g_k(x)\right| \leq
  %\epsilon$. Similarly, if $\varphi^{-1}(x) \notin
  %\bigcup_{i=0}^{l_k} C(k,i)$, $x \in C(k,0)$ and $g_k(\varphi^{-1}(x)) =
  %g_k(x) = 0$. If $\varphi^{-1}(x) \in \bigcup_{i=0}^{l_k} C(k,i)$ we have
  %instead that $\left|g_k\left(\varphi^{-1}(x)\right) - g_k(x)\right| \leq
  %\epsilon$. 
It is then easy to check that $\left\|g_k\circ \varphi - g_k\right\|
\leq \epsilon$ and $\left\|g_k\circ \varphi^{-1} - g_k\right\|
\leq \epsilon$, and that $g_k \circ \varphi, \ g_k \circ
\varphi^{-1}$ are both supported in $\bigcup_{i=0}^{l_k} C(k,i)$.  
Set
$$
f_1 = \sum_{k=1}^N g_k  
$$ 
and $f_2 = f_1 \circ \sigma$.
It is then easy to check that $0 \leq f_i \leq 1, \ \left\|f_i - f_i \circ
    \varphi\right\| \leq \epsilon, i =1,2$, $f_1f_2 = 0$, $f_1 \circ \sigma
  = f_2$. To check that $\omega\left(1 - f_1 -f_2\right) \leq \epsilon$ for every
  tracial state $\omega$ of $C(\Omega)
  \rtimes_{\varphi} \mathbb Z$, observe first that every
  $D_{\infty}$-invariant trace state is given by integration with
  respect to a $D_{\infty}$-invariant Borel probability measure on
  $\Omega$; this follows e.g. from a more general result which can be found
  as Theorem 4.5 of \cite{Th2}. Hence it suffices to show that
\begin{equation}\label{esttrce}
\int_{\Omega}  1 - f_1 -f_2 \ d\mu \leq \epsilon   
\end{equation}
when $\mu$ is a $D_{\infty}$-invariant Borel probability measure on
$X$. To this end note that 
$$
\int_{\Omega} 1 - f_1 -f_2 \ d\mu \leq  \sum_{k=1}^N 4M \beta_k
$$   
where $\beta_k = \mu\left( C(k,0)\right)$. Since $\frac{4M}{J_k} \leq
\frac{4M}{R} \leq \epsilon$ and $\sum_{k=1}^N J_k\beta_k = 1$ this
yields (\ref{esttrce}).
\end{proof}
\end{thm}

\begin{cor}\label{trrank0} Assume that $\varphi$ is minimal. Then $C(\Omega) \rtimes_{\alpha} D_{\infty}$ has
  tracial rank zero in the sense of Lin, \cite{Lin2}.
\begin{proof} This follows from Theorem \ref{tracRokhlin} and Theorem 2.6 of \cite{Ph4}.
\end{proof}
\end{cor}

\begin{remark}\label{rokhlin}
The $\mathbb Z_2$ action of Theorem \ref{tracRokhlin}
does not have the (strict) Rokhlin property. Indeed, if it did the
crossed product $C(\Omega) \rtimes_{\alpha} D_{\infty}$ would be a simple
AT-algebra by Corollary 3.4 of \cite{OP} which is not possible since
we shall show below that the $K_0$-group of $C(\Omega) \rtimes_{\alpha}
D_{\infty}$ contains torsion.
\end{remark}

It follows \emph{almost} from the preceding theorem that $C(\Omega)
\rtimes_{\alpha} D_{\infty}$ is classified by K-theory, thanks to the
work of H. Lin. What is missing is that we need to show that $C(\Omega)
\rtimes_{\alpha} D_{\infty}$ satisfies the UCT. This technicality
is dealt with in the following proposition. 

\begin{defn}\label{UCTdef} A separable $C^*$-algebra $A$
  \emph{satisfies the UCT} when the canonical map
$$
\gamma_{A,D} : KK(A,D) \to \ \Hom \left(K_*(A),K_*(D)\right)
$$
is an isomorphism for every separable $C^*$-algebra $D$ such that
$K_0(D)$ and $K_1(D)$ are both divisible groups.
\end{defn}

That this formulation is equivalent to the more conventional, which
appears in the litterature on the classification of $C^*$-algebras, follows from
Proposition 23.8.1 of \cite{Bl}.

\begin{prop}\label{UCTthm} Let $A_1$ and $A_2$ be a separable
  $C^*$-algebras with a common $C^*$-subalgebra $B$. Assume that
  $A_1$ and $A_2$ satisfy the UCT, and assume that $B$
  is AF. Then the
  amalgamated free product $A_1 *_B A_2$ satisfies the UCT.
\begin{proof} Note that $A_1 *_B A_2 = \varinjlim_n A_1 *_{B_n} A_2$ when
  $B_1 \subseteq B_2 \subseteq B_3 \subseteq \dots $ is an increasing
  sequence of finite-dimensional $C^*$-subalgebras with dense union in
  $B$. As shown in \cite{Th3} $A_1 *_{B_n} A_2$ is KK-equivalent to
  the mapping cone for an embedding of $B_n$ into $A_1 \oplus
  A_2$. Therefore $A_1 *_{B_n} A_2$ satisfies the UCT since $A_1$,
  $A_2$ and $B_n$ all do. Since the class of separable $C^*$-algebras which satisfy the
  UCT is closed under inductive limits, cf. \cite{Bl}, it follows that
  $A_1 *_B A_2$ is in the class.
\end{proof}
\end{prop}

\begin{cor}\label{UCTforthecase} $C(\Omega) \rtimes_{\alpha} D_{\infty}$
  satisfies the UCT.
\begin{proof} As observed and used in Section 4 of \cite{BEK},
  $C(\Omega) \rtimes_{\alpha} D_{\infty}$ can be realized as the
  amalgamated free
  product
$$
\left(C(\Omega) \rtimes_{\sigma} \mathbb Z_2\right) *_{C(\Omega)}
\left(C(\Omega) \rtimes_{\varphi\sigma} \mathbb Z_2\right) .
$$ 
Both $C(\Omega) \rtimes_{\sigma} \mathbb Z_2$ and $C(\Omega)
\rtimes_{\varphi \sigma} \mathbb Z_2$ are sub-homogeneous and type I
(see Section 4 of \cite{BEK} or (\ref{isocrossed}) below), and $C(\Omega)$ is
AF, so Proposition \ref{UCTthm} applies.  
\end{proof}
\end{cor}

We now turn to the calculation of $K_*\left(C(\Omega)
    \rtimes_{\alpha} D_{\infty}\right)$. When $\varphi$ is not
minimal Theorem \ref{easycase} shows that $C(\Omega) \rtimes_{\alpha}
D_{\infty}$ is, up to an $M_2$-tensor factor, also the
crossed product coming from a minimal homeomorphism of the Cantor set,
and an obvious construction shows that all such algebras arise this
way. In particular, the $K$-theory can be found from \cite{HPS}; the
$K_1$-group is $\mathbb Z$ and the $K_0$-group is a simple dimension
group. Also from \cite{HPS} we conclude that all simple dimension
groups occur here.

It remains to consider the case where $\varphi$ is minimal. It follows then from \cite{N} and Lemma 4.2 of \cite{BEK} that
$K_0\left(C(\Omega) \rtimes_{\varphi} \mathbb Z\rtimes_{\sigma} \mathbb
  Z_2\right)$ is the cokernel of the map
\begin{equation}\label{themap}
\begin{xymatrix}{
K_0\left(C(\Omega)\right) \ar[r]^-{\left({i_1}_*,  {i_2}_*\right)} & K_0\left(C(\Omega)
  \rtimes_{\sigma} \mathbb Z_2\right) \oplus K_0\left(C(\Omega)
  \rtimes_{\varphi \sigma} \mathbb Z_2\right)}
\end{xymatrix}
\end{equation}
where $i_1 : C(\Omega) \to C(\Omega) \rtimes_{\sigma} \mathbb Z_2$ and
$i_2 : C(\Omega) \to C(\Omega) \rtimes_{\varphi \sigma} \mathbb Z_2$
are the canonical maps. Note that there is an isomorphism
\begin{equation}\label{isocrossed}
C(\Omega) \rtimes_{\sigma} \mathbb Z_2 \ \simeq \ \left\{ f \in
  C\left(\Omega, M_2\right) : \ f \circ \sigma = \Ad u (f)\right\}
\end{equation}
induced by the covariant representation $\left(\pi, u\right)$ where $u
=\left( \begin{smallmatrix}
  0 & 1 \\ 1 & 0 \end{smallmatrix} \right)$ and $\pi :
C(\Omega) \to C(\Omega, M_2)$ is given by
$$
 \pi(f) = \left( \begin{smallmatrix}
  f & 0 \\ 0 & f \circ \sigma \end{smallmatrix} \right) .
$$
For simplicity of notation we identify $C(\Omega) \rtimes_{\sigma}
\mathbb Z_2$ and $\left\{ f \in
  C\left(\Omega, M_2\right) : \ f \circ \sigma = \Ad u (f)\right\}$ in
the following. There is a map $Tr : C(\Omega)
  \rtimes_{\sigma} \mathbb Z_2 \to C(\Omega)$
defined such that 
$$
Tr(f)(x) = f_{11}(x) + f_{22}(x)
$$
when $f \in C(\Omega) \rtimes_{\sigma}
\mathbb Z_2$. This gives rise to a group homomorphism
$$
Tr_* : \ K_0\left(C(\Omega) \rtimes_{\sigma} \mathbb Z_2\right) \to
\left\{g \in C(\Omega,\mathbb Z) : \ g \circ \sigma = g\right\} .
$$

\begin{lemma}\label{freeorder2} There is a compact and open set $K
  \subseteq \Omega$ such that $K \cap \sigma(K) = \emptyset$ and $K \cup
  \sigma(K) = \Omega$.
\begin{proof} Let $x \in \Omega$. Since $\sigma(x) \neq x$ there is a
  clopen subset $U \subseteq \Omega$ such that $U \cap \sigma(U) =
  \emptyset$. By compactness of $\Omega$ we get a finite collection
  $U_1,U_2, \dots, U_N$ of clopen sets in $\Omega$ such that $U_i \cap
  \sigma(U_i) = \emptyset$ for all $i$ and $\bigcup_{i=1}^N U_i =
  \Omega$. Set $B_1 = U_1$ and define $B_i, i
  \geq 2$, recursively, such that 
$$
B_i = U_i \backslash \left[ \bigcup_{j=1}^{i-1} B_j \cup
  \sigma(B_j)\right].
$$
Set $K = B_1 \cup B_2 \cup \dots \cup B_N$.
\end{proof}
\end{lemma}

\begin{lemma}\label{iso} Let $K$
  be the clopen set of Lemma \ref{freeorder2}, and define $\beta
 : C(\Omega) \rtimes_{\sigma} \mathbb Z_2 \ \to \
C(K,M_2)$ such that $\beta(f)  = f|_K$. Then $\beta$ is a $*$-isomorphism.
\begin{proof} $\beta$ is clearly an injective $*$-homomorphism. To see
  that $\beta$ is surjective, let $g
  \in C(K,M_2)$. Define $h : \Omega \to M_2$ such that
$$
h(x) = \begin{cases} g(x), & \ x \in K \\ ug(\sigma(x))u , & \ x \in
  \sigma(K). \end{cases}
$$
Then $h \in C(\Omega) \rtimes_{\sigma} \mathbb Z_2$ and $\beta(h) = g$.
\end{proof}
\end{lemma}

\begin{cor}\label{trivK1} Assume $\varphi$ is minimal. It follows that 
$$K_1\left(C(\Omega) \rtimes_{\varphi} \mathbb Z \rtimes_{\sigma} \mathbb
  Z_2\right) = 0.
$$
\begin{proof} It follows from Lemma \ref{iso} that $K_1\left(C(\Omega)
    \rtimes_{\sigma} \mathbb Z_2\right) = 0$. For the same reason we
  have also that $K_1\left(C(\Omega)
    \rtimes_{\varphi \sigma} \mathbb Z_2\right) = 0$. The conclusion
  follows now from the six terms exact sequence of \cite{N}.
\end{proof}
\end{cor}

Note that the diagram
\begin{equation}
\begin{xymatrix}{
K_0\left(C(\Omega) \rtimes_{\sigma} \mathbb Z_2\right) \ar[r]^-{Tr_*}
\ar[d]_-{\beta_*} &  \left\{ f \in C(\Omega, \mathbb Z): \ f \circ
  \sigma = f\right\} \ar[d]^-{f \mapsto f|_K} \\
K_0\left(C(K,M_2)\right) \ar[r]^-{Tr_*} &  C(K,\mathbb Z)}
\end{xymatrix}
\end{equation}
commutes. It is wellknown and easy to show that the lower $Tr_*$-map
is an isomorphism. It follows in this way from Lemma \ref{iso} that
$$
Tr_* : \ K_0\left(C(\Omega) \rtimes_{\sigma} \mathbb Z_2\right) \to
\left\{ f \in C(\Omega, \mathbb Z): \ f \circ  \sigma = f\right\}
$$
is an isomorphism. Similarly, 
$$
Tr_* : \ K_0\left(C(\Omega) \rtimes_{\varphi\sigma} \mathbb Z_2\right)
\to
 \left\{ f \in C(\Omega, \mathbb Z): \ f \circ  \varphi\sigma = f\right\}
$$
is an isomorphism. 

Set $G_{\sigma} = \left\{ f \in C(\Omega,\mathbb Z) : \ f \circ \sigma
  = f\right\}$ and $G_{\varphi \sigma} = \left\{ f \in C(\Omega,\mathbb Z) : \ f \circ \varphi\sigma
  = f\right\}$. It follows that $K_0\left(C(\Omega) \rtimes_{\varphi}
  \mathbb Z \rtimes_{\sigma} \mathbb
  Z_2\right)$ is isomorphic to the cokernel of the map
\begin{equation*}\label{themap2}
C(\Omega, \mathbb Z) \ni h \ \mapsto \ (h+ h\circ \sigma, h + h \circ
\varphi \sigma) \in G_{\sigma} \oplus G_{\varphi \sigma} .
\end{equation*}
Consider the embedding $\mathbb Z \ni z \mapsto (z,z) \in G_{\sigma}
\oplus G_{\varphi\sigma}$ obtained by considering $z$ as a constant
function on $\Omega$. If $h \in C(\Omega,\mathbb Z)$ is such that
\begin{equation}\label{eqeven}
\left( {i_1}_*(h),{i_2}_*(h)\right) = (h+h \circ \sigma, h + h \circ \varphi\sigma) = (z,z)
\end{equation} 
for some $z
\in \mathbb Z$ we have in particular that $h\circ \sigma = h\circ
\varphi \sigma = h \circ \sigma \varphi^{-1}$ which implies that
$h\circ \sigma$ and hence also $h$ is constant because $\varphi$ is
minimal. But then (\ref{eqeven}) implies that $z$ is even. This shows
that there is an injection 
\begin{equation*}
t : \mathbb Z_2 \to 
K_0\left(C(\Omega) \rtimes_{\varphi} \mathbb Z \rtimes_{\sigma} \mathbb
  Z_2\right).
\end{equation*}

By composing the map $G_{\sigma} \oplus G_{\varphi \sigma} \ni (f,g) \mapsto f-g
\in C(\Omega, \mathbb Z)$ with the quotient map $C(\Omega,\mathbb Z)
\to C(\Omega, \mathbb Z)/(1-\varphi_*)\left(C(\Omega, \mathbb
  Z)\right)$ we obtain a map which clearly annihilates elements of the
form $(h+h \circ \sigma, h + h \circ \varphi \sigma)$ and hence give
rise to a map 
\begin{equation*}
s: K_0\left(C(\Omega) \rtimes_{\varphi}
  \mathbb Z \rtimes_{\sigma} \mathbb
  Z_2\right) \to C(\Omega, \mathbb Z)/(1-\varphi_*)\left(C(\Omega, \mathbb
  Z)\right).
\end{equation*}

\begin{lemma}\label{extension} Assume that $\varphi$ is minimal. 
$$
0 \to \mathbb Z_2  \overset{t}\to K_0\left(C(\Omega) \rtimes_{\varphi}
  \mathbb Z \rtimes_{\sigma} \mathbb
  Z_2\right) \overset{s}\to (1 + \sigma_*)\left(C(\Omega, \mathbb Z)\right)/(1-\varphi_*)\left(C(\Omega, \mathbb
  Z)\right) \to 0
$$
is exact.
\begin{proof} It is obvious that $s \circ t = 0$ so consider an
  element $(f,g)\in G_{\sigma} \oplus G_{\alpha \sigma}$ such that
  $f-g$ represent $0$ in $C(\Omega, \mathbb Z)/(1-\varphi_*)\left(C(\Omega, \mathbb
  Z)\right)$. There is then a function $h \in C(\Omega,\mathbb Z)$
such that 
\begin{equation}\label{eq300}
f-g = h - \varphi(h).
\end{equation}
Since $f \in G_{\sigma}, g \in G_{\varphi\sigma}$ we find that
\begin{equation}\label{eq301}
f - \sigma(g) =
\sigma(h) - \sigma\varphi(h)
\end{equation}
and 
\begin{equation}\label{eq302}
\varphi(f) - g = \varphi\sigma(h) -
\sigma(h).
\end{equation}
Subtracting (\ref{eq302}) from (\ref{eq300}) yields that $f
-\varphi(f) = h-\varphi(h) - \varphi\sigma(h) + \sigma(h)$ or
\begin{equation}\label{eq306}
f  - \sigma(h) - h - \varphi\left( f - \sigma(h) -h\right) = 0 . 
\end{equation}
In the same way (\ref{eq301}) and (\ref{eq300}) imply that $-g +
\sigma(g) = h - \varphi(h) -\sigma(h) + \sigma\varphi(h)$ or (since
$\sigma(g) = \varphi^{-1}(g)$),
\begin{equation}\label{eq305}
-g + \sigma(h) + \varphi(h) - \varphi^{-1}\left( -g + \sigma(h) +
  \varphi(h)\right) = 0.
\end{equation}
Since $\varphi$ is minimal we conclude from (\ref{eq306}) and
(\ref{eq305}) that $f  - \sigma(h) - h$ and $-g + \sigma(h) +
\varphi(h) $ are constant. Inserting this conclusion into
(\ref{eq300}) we see that the two constants must add to zero, i.e. there
is an integer $z \in \mathbb Z$ such that i.e. $f - \sigma(h) - h =z$ and $-g
+\sigma(h) + \varphi(h) = -z$. Set $k = \sigma(h)$ and note that
$(f,g) = \left(f-k-\sigma(k), g - k - \varphi \sigma(k)\right) =
(z,z)$ in $K_0\left(C(\Omega) \rtimes_{\varphi}
  \mathbb Z \rtimes_{\sigma} \mathbb
  Z_2\right)$. It
follows that the element of $K_0\left(C(\Omega) \rtimes_{\varphi}
  \mathbb Z \rtimes_{\sigma} \mathbb
  Z_2\right)$ which $(f,g)$ represents is in $t\left(\mathbb Z_2\right)$.

It remains now only to show that
\begin{equation}\label{range}
s\left(K_0\left(C(\Omega) \rtimes_{\varphi}
  \mathbb Z \times_{\sigma} \mathbb
  Z_2\right)\right) = (1 + \sigma_*)\left(C(\Omega, \mathbb Z)\right)/(1-\varphi_*)\left(C(\Omega, \mathbb
  Z)\right) .
\end{equation}
Let $(f,g) \in G_{\sigma} \oplus G_{\varphi\sigma}$. Let $K \subseteq
\Omega$ be the compact open set of Lemma \ref{freeorder2}  and set
$$
h_1(x) = \begin{cases} f(x), & \ x \in K \\ 0, & \ x \in \sigma(K)
  . \end{cases}
$$
Then $h_1 \in C(\Omega,\mathbb Z)$ and $h_1 + \sigma(h_1) =
f$. Similarly there is an $h_2 \in C(\Omega,\mathbb Z)$ such that and $h_2 + \varphi\sigma(h_2) =
g$. Since $f-g = (1 + \sigma_*)\left(h_1 -h_2\right) + \sigma(h_2) -
\varphi\sigma(h_2)$, we see that $f-g \in  (1 +
\sigma_*)\left(C(\Omega, \mathbb Z)\right)$ modulo $(1-\varphi_*)\left(C(\Omega, \mathbb
  Z)\right)$. This proves the inclusion $\subseteq$ in
(\ref{range}). Since the reversed inclusion is trivial, the proof is complete. 
\end{proof}
\end{lemma}

\begin{thm}\label{directsummand} Assume that $\varphi$ is
  minimal. Then
$$
K_1\left(C(\Omega)\rtimes_{\alpha} D_{\infty}\right) = 0
$$
and
$$
K_0\left(C(\Omega)\rtimes_{\alpha} D_{\infty}\right) \ \simeq \
\mathbb Z_2 \ \oplus \ (1 + \sigma_*)\left(C(\Omega, \mathbb Z)\right)/(1-\varphi_*)\left(C(\Omega, \mathbb
  Z)\right).
$$
\begin{proof} In view of Corollary \ref{trivK1} and Lemma
  \ref{extension} it remains only to show that the $\mathbb
  Z_2$-subgroup from the extension of Lemma \ref{extension} is a
  direct summand. But this follows from a theorem of Kulikov,
  cf. Theorem 24.5 of \cite{F}, since the quotient of the extension, $(1 + \sigma_*)\left(C(\Omega, \mathbb Z)\right)/(1-\varphi_*)\left(C(\Omega, \mathbb
  Z)\right)$, is torsion free. 
\end{proof}
\end{thm}

\begin{remark}\label{remarkBEK} 

The torsion subgroup $\mathbb Z_2$ was overlooked in \cite{BEK}. Its
appearance shows that $C(\Omega) \rtimes_{\alpha} D_{\infty}$ can fail
to be
AF when neither $\sigma$ nor $\varphi\sigma$ has fixed points, a
question left open in \cite{BEK}. With
this remark we try to
give a better description of this torsion element and explain why it
disappears when $\sigma$ or $\varphi\sigma$ has fixed points. 

Observe first that the unitary $u$ occuring in (\ref{isocrossed}) is unitarily equivalent to $v = \left(
    \begin{smallmatrix} 1 & 0 \\ 0 & -1 \end{smallmatrix}
  \right)$. Specifically $wuw^* = v$ where
$$
w = \frac{1}{\sqrt{2}} \left( \begin{matrix} 1 & 1 \\ 1 & -1 \end{matrix}
\right).
$$
There is therefore a $*$-isomorphism $C(\Omega)
\rtimes_{\sigma} \mathbb Z_2  \simeq A$ where
\begin{equation}\label{io1}
A = \left\{ f \in C\left(\Omega,
    M_2\right) : \ f \circ \sigma = \Ad v (f) \right\} .
\end{equation}
Similarly, therefore is a $*$-isomorphism $C(\Omega)
\rtimes_{\varphi\sigma} \mathbb Z_2  \simeq B$ where
\begin{equation}\label{io2}
B = \left\{ f \in C\left(\Omega,
    M_2\right) : \ f \circ \varphi\sigma = \Ad v (f) \right\} .
\end{equation}
Under this identification $i_1$ and $i_2$ are given by 
$$
i_1(h) = \frac{1}{2} \left( \begin{matrix} h + h\circ \sigma    &  h -
    h\circ \sigma \\ h-h\circ \sigma  & h + h \circ \sigma \end{matrix}
  \right),
$$
and
$$
i_2(h) = \frac{1}{2} \left( \begin{matrix} h + h\circ \varphi\sigma    &  h -
    h\circ \varphi\sigma \\ h-h\circ \varphi\sigma  & h + h \circ \varphi\sigma \end{matrix}
  \right),
$$
respectively. The projection 
$$
p = \left(\begin{smallmatrix} 0 & 0 \\ 0 & 1 \end{smallmatrix} \right)
$$
is then an element of both $A$ and $B$ in the obvious way. The element $\left([p],[p]\right) \in
K_0(A)\oplus K_0(B)$ is not in
the range of $\left({i_1}_*,{i_2}_*\right)$, for if it was there would
have to be a function $h \in C(\Omega,\mathbb Z)$ such that 
\begin{equation}\label{hurra}
Tr_*\left([p]\right) = 1 = h + h\circ \sigma = h + h\circ
\varphi\sigma .
\end{equation}
The last identity of (\ref{hurra}) and minimality of $\varphi$ would force $h$ to be
constant which makes the second identity in (\ref{hurra})
impossible. This shows that $([p],[p])$ represents a non-zero element
of $K_0\left(C(\Omega)\rtimes_{\alpha} D_{\infty}\right)$ when
$\varphi$ is minimal, regardless of whether or not the action $\alpha$
is free. But when it \emph{is} free the element has order two,
because in this case       
 we have the element
$$ 
a = \left( \begin{matrix} 0 & 1_K - 1_{\sigma(K)} \\ 0 & 0
  \end{matrix} \right) 
$$
in $A$, where $K$ is the clopen subset of $\Omega$ from Lemma
\ref{freeorder2}). Since $aa^* =1_A -p$ and $a^*a = p$ this shows that
$[1_A] = 2[p]$ in $K_0(A)$. It follows in the same way that $[1_B] = 2[p]$ in
$K_0(B)$. Since $\left([1_A],[1_B]\right) =
\left({i_1}_*\left[1_{C(\Omega)}\right],{i_2}_*\left[1_{C(\Omega)}\right]\right)$
we see that $([p],[p])$ represents an element of order 2 in the
cokernel of $\left({i_1}_*,{i_2}_*\right)$. In contrast, when $\sigma$ or $\varphi\sigma$ has a fix point there is no
substitute for the element $a$ in $A$ or $B$ and the
element $([p],[p])$ has in fact infinite order in the
cokernel of $\left({i_1}_*,{i_2}_*\right)$ in this case, as shown in \cite{BEK}. 
\end{remark}

\begin{thm}\label{mainoncross} Let $\alpha$ be a free minimal action
  of the infinite dihedral group $D_{\infty}$ on the Cantor set $\Omega$, generated by the
  homeomorphisms $\varphi$ and $\sigma$ satisfying (\ref{order2}) and
  (\ref{conjcond}). Then 
$$C(\Omega)
  \rtimes_{\alpha} D_{\infty}
$$ 
is a simple unital AH-algebra with real
  rank zero and no dimension growth. It is an AT-algebra if and only
  if $\varphi$ is not minimal. 
\begin{proof} That $C(\Omega)
  \rtimes_{\alpha} D_{\infty}$ is AT when $\varphi$ is not minimal follows from Theorem
  \ref{easycase}. When $\varphi$ is 
  minimal we conclude from Corollary \ref{trrank0} and Corollary
  \ref{UCTforthecase} that $A$ has tracial rank zero and satisfies the
  UCT. It follows then from \cite{Lin4} that $A$ is AH with no
  dimension growth. That the real rank is zero follows from \cite{Lin1}
  and \cite{Lin2} because the
  tracial rank is zero. Finally, it follows from Theorem
  \ref{directsummand} that the algebra is not AT when $\varphi$ is minimal.
\end{proof}
\end{thm}

This ends our excursion into crossed products of the form $C(\Omega)
\rtimes_{\alpha} D_{\infty}$ because we have obtained all the
information we need for the study of the heteroclinic and homoclinic
$C^*$-algebras of 1-solenoids.

\subsection{The inductive limit decomposition of $B_{\overline{h}^{-1}}\left(\overline{\Gamma}\right)$.}

We can now put the pieces together and obtain the following

\begin{thm}\label{invmainthm} Let $({\Gamma},{h})$
  be a pre-solenoid and
   $ B_{\overline{h}^{-1}}\left(\overline{\Gamma}\right)$ the
  heteroclinic algebra of the inverse $\overline{h}^{-1} :
  \overline{\Gamma} \to \overline{\Gamma}$ of the 1-solenoid
  $\left(\overline{\Gamma}, \overline{h}\right)$. Then
  $B_{\overline{h}^{-1}}\left(\overline{\Gamma}\right) \simeq A
  \otimes \mathbb K$ where $A$ is a unital simple AH-algebra with real
  rank zero, no dimension growth and a unique trace state. 
\begin{proof} Combine Theorem \ref{1mainrev}
  with Proposition \ref{trace} and Theorem \ref{mainoncross}.
\end{proof}
\end{thm}

The $K$-theory of $
B_{\overline{h}^{-1}}\left(\overline{\Gamma}\right)$ is described in
Theorem \ref{directsummand} when $\varphi$ is minimal and in Theorem
\ref{easycase} when $\varphi$ is not minimal. By Lemma
\ref{minimalityofphi} and Theorem \ref{orientthm} $\varphi$ is minimal
if and only if $\left(\overline{\Gamma}, \overline{h}\right)$ is
orientable. This gives us the following conclusion in the same way we
obtained Theorem \ref{orient2}:

\begin{thm}\label{invorinet} Let $(\Gamma, h)$ be an oriented
  pre-solenoid. Then the heteroclinic algebra $
  B_{\overline{h}^{-1}}\left(\overline{\Gamma}\right)$ of
  $\overline{h}^{-1} : \overline{\Gamma} \to \overline{\Gamma}$ is a simple
  stable AT-algebra of real rank zero.
\end{thm}

\begin{example}\label{counterexample} 

The $K_0$-groups of
$K_0\left(B_{\overline{h}^{-1}}\left(\overline{\Gamma}\right)\right)$
agrees with the first \v Cech cohomology group of $\overline{\Gamma}$
when $(\Gamma,h)$ is oriented. This follows from \cite{Y2} and
\cite{Y4}. There is even a natural way to define an order on the first
\v Cech cohomology group which corresponds to the order on
$K_0\left(B_{\overline{h}^{-1}}\left(\overline{\Gamma}\right)\right)$.

I don't know if
$K_0\left(B_{\overline{h}}\left(\overline{\Gamma}\right)\right)$ has
an interpretation in terms of known invariants for 1-solenoids, even
in the orientable case. In many cases the two groups,
$K_0\left(B_{\overline{h}^{-1}}\left(\overline{\Gamma}\right)\right)$
and $K_0\left(B_{\overline{h}}\left(\overline{\Gamma}\right)\right)$
are isomorphic as dimension groups and it follows from the results in
this paper that $B_{\overline{h}^{-1}}\left(\overline{\Gamma}\right)$
and $B_{\overline{h}}\left(\overline{\Gamma}\right)$ are then
isomorphic. This happens for example for the 1-solenoids of Remark
\ref{exampl}. In fact, it is not so easy to come up with an example of
an
oriented 1-solenoid for which
$B_{\overline{h}^{-1}}\left(\overline{\Gamma}\right) \neq
B_{\overline{h}}\left(\overline{\Gamma}\right)$, but here is one: In
the setting of Remark \ref{exampl} let $h : \Gamma \to \Gamma$ be
given by the wrapping rule
$$
a \ \mapsto \ a^{65}b^7, \ b \ \mapsto \ a^{24}b^{67} .
$$
By Corollary 3.7 of \cite{Y4}
$K_0\left(B_{\overline{h}^{-1}}\left(\overline{\Gamma}\right)\right)$
is the stationary dimension group given by the matrix
$$
 A = \left( \begin{matrix} 65 & 7 \\ 24 & 67 \end{matrix} \right) .
$$
This group is not isomorphic to
$K_0\left(B_{\overline{h}}\left(\overline{\Gamma}\right)\right)$, even
when we ignore the order. Indeed, by using the description of
$B_{\overline{h}}\left(\overline{\Gamma}\right)$ given in Section
\ref{het} it can be shown that
$K_0\left(B_{\overline{h}}\left(\overline{\Gamma}\right)\right)$ is
the stationary dimension group corresponding to the \emph{transpose}
$A^t = \left( \begin{smallmatrix} 65 & 24 \\ 7 & 67 \end{smallmatrix}
\right)$ of $A$. These dimension groups are not even isomorphic as
groups, cf. Example 3.6 of \cite{BJKR} and Appendix C of \cite{Th1}.
\end{example}

\section{The homoclinic algebra of a 1-solenoid}\label{homo}

By using, for the first time, that 1-solenoids are Smale spaces we can
combine the preceding results with \cite{Pu} to obtain the following

\begin{thm}\label{homomain} Let $(\Gamma,h)$ be a pre-solenoid. The homoclinic $C^*$-algebra
  $A_{\overline{h}}\left(\overline{\Gamma}\right)$ of the 1-solenoid
  $\left(\overline{\Gamma}, \overline{h}\right)$ is a simple unital
  $AH$-algebra of real rank zero with no dimension growth and a unique
  trace state.
\begin{proof} By Theorem 3.1 of \cite{Pu} 
  $A_{\overline{h}}\left(\overline{\Gamma}\right)$ is stably
  isomorphic to $B_{\overline{h}}\left(\overline{\Gamma}\right) \otimes
  B_{\overline{h}^{-1}}\left(\overline{\Gamma}\right)$. It follows then from Theorem \ref{invmainthm} and Theorem
  \ref{main} that $A_{\overline{h}}\left(\overline{\Gamma}\right)$ is
  stably isomorphic to a tensor-product $D_1 \otimes D_2$ where $D_1$
  and $D_2$ are simple unital
  AH-algebras of real rank zero with no dimension growth and a unique
  trace state. It is standard to deduce from this that $D_1 \otimes
  D_2$ itself is a simple unital AH-algebra with no dimension growth
  and a unique trace state. To see that $D_1 \otimes D_2$ has real
  rank zero note that AH-algebras without dimension growth are
  inductive limits of recursive sub-homogeneous algebras in the sense
  of \cite{Ph3}. It follows then from Proposition 1.10 and Theorem 4.2
  of \cite{Ph3} that $D_1 \otimes D_2$ has real rank zero.
\end{proof}
\end{thm}

Since $A_{\overline{h}}\left(\overline{\Gamma}\right)$ is stably
isomorphic to the tensor product
$B_{\overline{h}^{-1}}\left(\overline{\Gamma}\right) \otimes
B_{\overline{h}}\left(\Gamma\right)$ we can find the K-groups of
$A_{\overline{h}}\left(\overline{\Gamma}\right)$ by using the results
from Section \ref{het} and Section \ref{invhet} in combination with
the K\" unneth theorem, \cite{Bl}.

By combining Theorem \ref{invorinet} and Theorem \ref{orient2} we
obtain the following

\begin{thm}\label{invorinet2} Let $(\Gamma, h)$ be an oriented
  pre-solenoid. Then the homoclinic algebra $
  A_{\overline{h}}\left(\overline{\Gamma}\right)$ is a simple unital
  AT-algebra of real rank zero with a unique trace state.
\end{thm}

\end{document}